\documentclass[preprint,onefignum,onetabnum]{siamart220329}

\usepackage{lipsum}
\usepackage{amsfonts}
\usepackage{graphicx}
\usepackage{epstopdf}
\usepackage{algorithmic}
\ifpdf
  \DeclareGraphicsExtensions{.eps,.pdf,.png,.jpg}
\else
  \DeclareGraphicsExtensions{.eps}
\fi

\newsiamremark{remark}{Remark}
\newsiamremark{hypothesis}{Hypothesis}
\crefname{hypothesis}{Hypothesis}{Hypotheses}
\newsiamthm{claim}{Claim}
\newsiamthm{example}{Example}

\headers{Structure-preserving scheme for  PNP equation}{F.  Tong and Y. Cai}

\title{Positivity preserving and mass conservative projection method for the Poisson-Nernst-Planck equation\thanks{Submitted to the editors DATE.
\funding{This work was partially supported by NSFC grants 12171041 and 11771036 (Y. Cai)}}}

\author{Fenghua Tong\thanks{Laboratory of Mathematics and Complex Systems and School of Mathematical Sciences, Beijing
Normal University, Beijing 100875, China
  (\email{fenghua.tong@mail.bnu.edu.cn}).}
\and Yongyong Cai\thanks{Laboratory of Mathematics and Complex Systems and School of Mathematical Sciences, Beijing
Normal University, Beijing 100875, China
  (\email{yongyong.cai@bnu.edu.cn}).}
}

\usepackage{amsopn}

\ifpdf
\hypersetup{
  pdftitle={Structure-preserving scheme for the PNP equations},
  pdfauthor={F.  Tong, and Y. Cai}
}
\fi

\begin{document}

\maketitle
\begin{abstract}
We propose and analyze a novel approach to construct structure preserving approximations for the Poisson-Nernst-Planck equations, focusing on  the  positivity preserving and mass conservation properties. 
The strategy consists of a standard time marching step with a projection (or correction) step to satisfy the desired physical constraints (positivity and mass conservation). Based on the $L^2$ projection, we construct a second order Crank-Nicolson type finite difference scheme, which is  linear (exclude the very efficient $L^2$ projection part), positivity preserving and mass conserving. Rigorous error estimates 
in $L^2$ norm are established, which are both second order accurate in space and time.
The other choice of projection, e.g. $H^1$ projection,  is discussed.  
 Numerical examples are presented to verify the theoretical results and demonstrate the efficiency of the proposed method.
\end{abstract}



\section{Introduction}
In this article, we  consider the two-component system of Poisson-Nernst-Planck  (PNP) equations \cite{NernstDieEW,PlanckUeberDE,2009Prohl,Mock1972,LU20106979} as
\begin{equation}\label{PDE}
\begin{aligned}
&p_t = \Delta p +\nabla\cdot( p \nabla \phi ), &\mathbf{x}\in\Omega, t>0,\\
&n_t = \Delta n-\nabla\cdot( n  \nabla \phi  ),&\mathbf{x}\in\Omega, t>0,\\
&-\Delta \phi  = p - n, &\mathbf{x}\in\Omega, t>0,\\
\end{aligned}
\end{equation}
 and the initial conditions are given as
\begin{equation}\label{ini}
p(\mathbf{x},0)=p_{0}(\mathbf{x})\ge0, \ \ n(\mathbf{x},0)=n_{0}(\mathbf{x})\ge0,   \quad\mathbf{x}\in\Omega,
\end{equation}
where $\Omega\subset\mathbb{R}^d$ ($d=1,2,3$)  is  an interval in one dimension (1D) or a rectangle in two dimensions (2D) or a cube in three dimensions (3D),  $\mathbf{x}\in\Omega$  is the spatial coordinate ($\mathbf{x}=x$ in 1D, $\mathbf{x}=(x,y)^T$ in 2D and $\mathbf{x}=(x,y,z)^T$ in 3D),  $t\ge0$ is time,   $\Delta$ is the Laplace operator, $\nabla$ is the gradient operator,  $\nabla\cdot$ is the divergence operator,  $p:=p(\mathbf{x},t)$ and $n:=n(\mathbf{x},t)$ represent the concentration of positively and negatively charged particles respectively, and $\phi:=\phi(\mathbf{x},t)$ denotes the electric potential.  Here,  periodic boundary conditions are imposed for $n$, $p$ and $\phi$ in \eqref{PDE}, and the case of homogeneous Neumann boundary conditions is similar.

Under the periodic boundary conditions, the PNP system \eqref{PDE} conserves the mass of each component, i.e.
\begin{equation}\label{eq:mass}
M_p=\int_{\Omega} p_0(\mathbf{x})d\mathbf{x} = \int_{\Omega} p(\mathbf{x},t)d\mathbf{x},\ \ M_n=\int_{\Omega} n_0(\mathbf{x})d\mathbf{x} = \int_{\Omega} n(\mathbf{x},t)d\mathbf{x},\ \ \forall t>0.
\end{equation}
For the consistency of periodic boundary conditions, we require  that
\begin{equation}\label{eq:consist}
M_p=M_n>0,
\end{equation}
which means the system is at the neutral state (net charge is zero).
In addition, in order to uniquely determine the electric potential $\phi$ in \eqref{PDE}, we set 
\begin{align}\label{interalphi}
   \int_{\Omega} \phi d\mathbf{x} =0.
\end{align}

For non-negative regular initial data \eqref{ini}, the PNP system \eqref{PDE} admits solutions with nonnegative ion densities, i.e. $p(\mathbf{x},t),n(\mathbf{x},t)\ge0$. In addition to the mass conservation \eqref{eq:mass},    \eqref{PDE} also enjoys energy laws.  

PNP system can be  viewed  as  a Wasserstein  gradient flow \cite{Luigi,Kinderlehrer2015AWG,2021Unconditionally}  driven by the total free energy  
\begin{align}\label{energytotal}
E(t)= \int_{\Omega} \left( p(\mathbf{x},t) \ln p(\mathbf{x},t) + n(\mathbf{x},t) \ln n(\mathbf{x},t) + \frac{1}{2} |\nabla \phi|^2  \right)d \mathbf{x},
\end{align}
which indicates energy dissipation property $\frac{\rm d}{\rm dt}E(t)\leq0$.
On the other hand,  the electric potential energy (the last term of total free energy \eqref{energytotal}) \cite{Biler1994TheDS,He2016AnEP, 2021Unconditionally}
\begin{align}\label{energyel}
E_{\phi}(t)=\frac{1}{2} \int_{\Omega}|\nabla \phi|^2 d \mathbf{x},  
\end{align}
 satisfies another energy identity as
$\frac{d}{d t}E_{\phi}(t)=-\int_{\Omega}\left((p-n)^2+(p+n)|\nabla \phi|^2\right) d \mathbf{x}\leq0$.

The PNP  system \eqref{PDE} is a macroscopic model  \cite{NernstDieEW,PlanckUeberDE} to describe  the potential difference in a galvanic cell.
Such a system and its extensions have been developed for various physical problems and applications,  including   the semiconductor theory (where they are known as  the  drift-diffusion  model  introduced  by Van Roosbroeck \cite{Gajewski1986OnTB,Semiconductor}),  biological systems \cite{Mock1972,LU20106979}, electrochemistry \cite{Kili2007Steric}.  In addition,  the PNP system can be coupled with the Navier-Stokes equation \cite{Schmuck2009ANALYSISOT,PNPcoupledNSLiuchun} or the convection equation  \cite{PNPcoupledConvection} to  address the interactions between electric field and complex flow field \cite{PhysRevE}.    

For the theoretical studies, there have been  many efforts devoted to the investigation of PNP system \eqref{PDE}, including the well-posedness \cite{Biler1994TheDS,YuExistence, Biler1994,Biler1993},  stationary states \cite{Cartailler2017A,Lyu2021NearAF,Gajewski1986OnTB,solexist}, asymptotic behavior \cite{Biler1994TheDS}, etc.
For numerical treatments on PNP \eqref{PDE}, Prohl and Schmuck \cite{2009Prohl}  developed a fully implicit finite element method that maintains  the mass conservation,  ion concentration positivity and electric potential energy law.  Flavell et al.  proposed a finite difference scheme  satisfying the conservation property \cite{2014A}.  
Based on a reformulation of Nernst Planck equation into a
nonlinear diffusive flux form,  e.g.  $\nabla\cdot( e^{\phi} \nabla(\frac{n}{e^{\phi}}))$ or $\nabla\cdot( n \nabla(\log(n)))$ in the Wasserstein gradient form, various numerical schemes have been developed for solving the
time dependent PNP equation, which admit the mass conservation, ion concentration
positivity and/or free energy dissipation \cite{DING2019108864,FU2022115031,Liu2020APE,LIU2022110777,He2019APP,METTI20161,ZHANG2022111086}.

Recently,  a new Lagrange multiplier approach has been introduced to construct positivity preserving  schemes for parabolic type systems in \cite{Vegt2019PositivityPL},  and to deal with additional constraints,  such as mass conservation and   bound preserving \cite{Cheng2020GlobalCP,ChengQing1,ChengQing2}. On the other hand, a natural way to ensure the physical constraint  numerically  is to project (or correct) the numerical solutions onto the constrained manifold, e.g. the projection methods \cite{Temam1968UneMD,shenNs,Cai2018ErrorEF} for incompressible fluid equations.    The cut off method \cite{libuyang,LU2013cutoff,yangjiang2022} for maximal bound preserving property of the Allen-Cahn equation can be viewed as such projection/correction methods.  The purpose of this work is to explore the projection/correction approach to numerically solve PNP \eqref{PDE}.

In this paper, we view the positivity of $(p,n)$ and the mass conservation \eqref{eq:mass} as physical constraints for PNP \eqref{PDE}.  Based on the projection (correction) strategy,   we follow the prediction-correction step, i.e. we first compute an intermediate numerical solution,   by  using a backward Euler scheme or Crank–Nicolson scheme with explicit treatment on the nonlinear terms in time, then project (correct) the intermediate solution onto the constrained manifold with positivity and mass conservation (we set such constraints correspondingly at the discrete level after finite difference discretization in space).  The projection step can be defined through a convex $L^2$ minimization problem to ensure the positivity and the mass conservation, where efficient solvers through Lagrange multipliers exist \cite{Cheng2020GlobalCP,ChengQing1,ChengQing2}.   The desired Crank-Nicolson type finite difference method with $L^2$ projection is shown to be second order convergent in both space and time in $L^2$ norm. In addition, different projection strategy can be employed, e.g. $H^1$ projection onto the constrained manifold. 
Moreover, $L^2$ error estimates (for $(n,p)$) are natural for the $L^2$ projection, while $H^1$ estimates (for $(n,p)$) are natural with $H^1$ projection.

The rest of the paper is organized as follows. In section \ref{sec:L2discrete}, we present the fully discrete Crank-Nicolson type finite difference scheme with $L^2$ projection strategy  for the PNP equation \eqref{PDE},  which enforces the mass conservation and positivity preservation at the discrete level. 
The error estimates are presented in section \ref{sec:L2erroestimate}.   Section \ref{sec:H1discrete} is devoted to the study of the $H^1$ projection.   In  section \ref{sec:example},  we present several numerical experiments to validate the accuracy and efficiency of our scheme.  Some conclusion are drawn in section \ref{sec:conclusion}.

\section{Finite difference scheme with $L^2$ projection}
\label{sec:L2discrete}
In this section, we introduce the finite difference scheme for solving the PNP system \eqref{PDE}. 
For the simplicity of presentation, we focus on the 2D case, i.e.  $\Omega$ is a  rectangular domain with $\Omega=[a_l, a_r] \times[b_l, b_r]$. The numerical scheme can be easily extended to 1D and 3D cases for tensor grids and the results remain the same.

\subsection{Finite difference  discretization}
Choosing $N_x$ and $N_y$ two positive integers, the domain $\Omega$ is uniformly partitioned with mesh sizes $h_x=\frac{a_r-a_l}{N_x}, h_y=\frac{b_r-b_l}{N_y}$, and the numerical domain $\Omega_h$ is
\begin{equation*}
\Omega_h=\left\{\left(x_i, y_j\right) \mid x_i=a_l+i h_x, y_j=b_l+ j  h_y, 0 \leq i \leq N_x,  0  \leq j \leq N_y\right\}.
\end{equation*}
We write the discrete grid function $u_h(\mathbf{x}) \in \mathbb R$ ($\mathbf{x}\in \Omega_h$), which can be regarded as a vector consisting of  $u_{i,j}=u_h(x_i,y_j) $($(x_i,y_j)\in\Omega_h$). For the discrete periodic boundary conditions,   we introduce the periodic discrete  function space for $u_h$ (or $u_{i,j}$) as
\begin{equation}\label{periobc}
X=\left\{u_h:\Omega_h\to\mathbb{R}|u_{0, j}=u_{N_x,j}, \; u_{i, 0}=u_{i,  N_y},\; \forall 0\leq i\leq N_x,\, 0\leq j\leq N_y\right\},
\end{equation}
and the mean zero subspace
\begin{equation}
X_0=\left\{u_h\in X\bigg| \sum\limits_{i,j=1}^Nu_{i,j}=0\right\}.
\end{equation}
 For the simplicity of presentation, we assume  $a_r-a_l=b_r-b_l$, $h_x=h_y=h, N_x=N_y=N$ in the subsequent discussion, i.e. a uniform tensor grids on a 2D square. 
For $u_h\in X$, we denote the  averaging operator $A_{x}(A_y): X\to X$ as
\begin{equation*}
\begin{aligned}
\mathcal A_x u_{i,j}=\frac{u_{i, j}+u_{i-1, j}}{2}, \quad \mathcal A_y u_{i, j}=\frac{u_{i, j}+u_{i, j-1}}{2},
\end{aligned}
\end{equation*}
 and the forward/backward finite difference operators $D^{\pm}_x: X\to X$ in $x$ direction 
\begin{equation*}
D^-_x u_{i, j} = \frac{u_{i, j}-u_{i-1, j}}{h}, \quad D^+_x u_{i, j} = \frac{u_{i+1, j}-u_{i, j}}{h}, 
\end{equation*}
where similar difference operators $D^{\pm}_y: X\to X$ can be introduced as
$D^-_y u_{i, j} = \frac{u_{i, j}-u_{i, j-1}}{h},  D^+_y u_{i, j} = \frac{u_{i, j+1}-u_{i, j}}{h}$. 
To derive a second order finite difference discretization of the PNP system \eqref{PDE}, for $u_h,v_h\in X$, 
we introduce the notations
\begin{equation*}
\mathcal Av_{i,j} =
      \begin{pmatrix}
      \mathcal A_xv_{i,j} & 0 \\
      0 & \mathcal A_yv_{i,j} \\
     \end{pmatrix},
\quad 
    \nabla_h^- u_{i, j}=
      \begin{pmatrix}
     D^-_x  \\
      D^-_y   \\
     \end{pmatrix}
u_{i, j},
\quad 
    \nabla_h^+ u_{i, j}=
      \begin{pmatrix}
     D^+_x  \\
      D^+_y   \\
     \end{pmatrix}
u_{i, j},
 \end{equation*}
 and the   2D discrete elliptic operator with variable coefficients can be written as
\begin{equation*}
\begin{aligned}
\nabla^+_h\cdot (( \mathcal A v)_{i,j} \nabla^-_h u_{i,j}) = \begin{pmatrix}
    D^+_x \quad
     D^+_y   
     \end{pmatrix}
      \begin{pmatrix}
      \mathcal A_x v_{i,j} & 0 \\
      0 & \mathcal A_y v_{i,j} \\
     \end{pmatrix}
      \begin{pmatrix}
     D^{-}_x \\
      D^{-}_y   \\
     \end{pmatrix}
 u_{i, j},
\end{aligned}
 \end{equation*}
while the standard 2D discrete Laplacian operator is
\begin{equation*}
\begin{aligned}
\Delta_h  u_{i, j}: &=\nabla_h^+\cdot( \nabla^-_h u_{i,j})=\frac{1}{h^2}\left( u_{i+1, j}+ u_{i-1, j}+ u_{i, j+1}+ u_{i, j-1}-4  u_{i, j}\right) .
\end{aligned}
\end{equation*}

For $u_h,v_h\in X$, the discrete $L^2$  inner products, the $L^2$ norm and the $L^p$ ($1\leq p<\infty$) norm  are  given by
\begin{equation*}
\langle u_h, v_h\rangle=h^2\sum_{i,j=1}^{N}  u_{i, j}  v_{i, j} , \quad\|u_h\|^2=\langle u_h, u_h\rangle,\quad
\|u_h\|_p^p=h^2\sum_{i,j=1}^{N} |u_{i, j}|^p,
\end{equation*}
and the $L^\infty$ norm is defined as $\|u\|_\infty=\max_{1\leq i,j\leq N} |u_{i,j}|$.
Without confusion, we shall use the same notation for the vector case, i.e. $\langle U, V\rangle=h^2\sum_{i,j=1}^{N}  U_{i, j}\cdot  V_{i, j}$, $\|U\|^2=\langle U, U\rangle$, $\|U_h\|_p^p=h^2\sum_{i,j=1}^{N} \|U_{i, j}|^p$ and $\|U_h\|_\infty=\max_{1\leq i,j\leq N} |U_{i,j}|$ ($|U_{i,j}|$ stands for the usual Euclidean norm) for $U_h,V_h\in X\times X$. 
The discrete semi-$H^1$ norm and the $H^1$ norm for $u_h\in X$ are 
\begin{equation*}
\|\nabla_h^- u_h\|^2=\langle\nabla_h^- u_h, \nabla_h^- u_h\rangle,\quad \|u_h\|^2_{H^1}=\|u_h\|^2+\|\nabla_h^- u_h\|^2.
\end{equation*}

\begin{lemma}\label{lem:summation}(Summation by parts)  For any grid functions $u_h, v_h,w_h \in X$,  the following identities are valid
\begin{align*}
\langle \nabla_h^- u_h, \nabla_h^- v_h\rangle=-\langle \Delta_h u_h, v_h\rangle,\quad \langle (\mathcal{A}w_h)\nabla_h ^-u_h, \nabla_h^- v_h\rangle=-\langle \nabla_h^+\cdot((\mathcal{A}w_h)\nabla^-_h u_h), v_h\rangle.
\end{align*}
 \end{lemma}
\subsection{Finite difference methods with positivity preservation and mass conservation}
Now, we are ready to construct our finite difference schemes based on the prediction-correction strategy. Choose the time step size $\tau>0$,  and the time steps are $t_k:=k\tau$ for $k=0,1,\ldots $.  Let $({p}^{k}_h(\mathbf{x}),{n}^{k}_h(\mathbf{x}),{\phi}^{k}_h(\mathbf{x}) )\in X\times X\times X_0$ ($\mathbf{x} \in \Omega_{h}$) be the numerical approximation of  the exact solution $(p,n,\phi)$ of PNP \eqref{PDE} on $\Omega_h$ at time $t_k$.  For any continuous periodic function $v(\mathbf{x})$ on $\Omega$, we define the interpolation operator $I_h$  on the numerical domain $\Omega_h$ as
\begin{align*}
I_h v \in X \quad \text{with} \quad {I_h v}(x_i,y_j) = v(x_i,y_j),  \quad \forall   \mathbf{x}_{i,j}=(x_i,y_j)\in\Omega_h.
\end{align*}

We assume the initialization is compatible, i.e. $\langle p_h^0,1\rangle=\langle n_h^0,1\rangle$. For the initial data $(p_0,n_0)$ given in \eqref{ini},  we assume $\min\left\{\langle I_hp_0,1\rangle,\langle I_hn_0,1\rangle\right\}>0$ (valid for sufficiently small $h$) and set $M_0=\max\left\{\langle I_hp_0,1\rangle,\langle I_hn_0,1\rangle\right\}>0$,
\begin{equation}\label{eq:ini}
p_h^0(x_i,y_j)=\varepsilon_1p_0(x_i,y_j),\quad n_h^0(x_i,y_j)=\varepsilon_2n_0(x_i,y_j),\quad -\Delta_h\phi_h^0=p_h^0-n_h^0,
\end{equation}
where $\varepsilon_1=M_0/\langle I_hp_0,1\rangle$, $\varepsilon_2=M_0/\langle I_hn_0,1\rangle$ and  $\phi_h^0\in X_0$ is uniquely determined.

For $k=0,1,2,\cdots$,  we follow the prediction-projection strategy as:

\textbf{Step 1}: intermediate solutions $\tilde{p}^{k+1}_{h}, \tilde{n}^{k+1}_{h}\in X$ are obtained from   the Crank-Nicolson finite difference scheme (CNFD) with the Adams–Bashforth strategy for the nonlinear part:
\begin{align}
\label{pseca}&\frac{ \tilde{p}^{k+1}_h  - p^{k}_h }{\tau}=\frac{1}{2}\Delta_h ( \tilde p^{k+1}_h + p^{k}_h)   + \nabla_h^+\cdot \left((\mathcal{A}p_{h}^{k+\frac{1}{2},*}) \nabla_h^- \phi_{h}^{k+\frac{1}{2},*}\right),\\
\label{psecb}&\frac{ \tilde{n}^{k+1}_h - n^{k}_h }{\tau}=\frac{1}{2}\Delta_h ( \tilde n^{k+1}_h + n^{k}_h)   - \nabla_h^+ \cdot \left((\mathcal{A}n_{h}^{k+\frac{1}{2},*}) \nabla_h^- \phi_{h}^{k+\frac{1}{2},*}\right),
\end{align}
where
\begin{equation*}
\begin{aligned}
& p^{k+\frac{1}{2},*}_h=\frac{3}{2}p^{k}_h -\frac{1}{2} p^{k-1}_h,\qquad  n^{k+\frac{1}{2},*}_h=\frac{3}{2}n^{k}_h -\frac{1}{2} n^{k-1}_h,\qquad  {\phi}^{k+\frac{1}{2},*}_h=\frac{3}{2}{\phi}^{k}_h -\frac{1}{2} {\phi}^{k-1}_h.\\
\end{aligned}
\end{equation*}

 \textbf{Step 2}: the correction step.
In general, the intermediate numerical solutions $\tilde{p}^{k+1}_h$ and $\tilde{n}^{k+1}_h$ computed by the linear semi-implicit schemes fail the positivity preservation \cite{libuyang}. 
Here,  we treat the positivity and mass conservation at the discrete level as the constraints for the numerical solution at $t_{k+1}$ as $(p_h^{k+1},n_h^{k+1})$ with $\langle {p}^{k+1}_h,1\rangle=\langle{p}^{0}_h,1\rangle$, $\langle {n}^{k+1}_h,1\rangle=\langle{n}^{0}_h,1\rangle$ and $p_h^{k+1},n_h^{k+1}\ge0$.
A most natural way to obtain such $(p_h^{k+1},n_h^{k+1})$ from  $(\tilde{p}^{k+1}_{h}, \tilde{n}^{k+1}_{h})$ is to
 project  the nodal vectors $\tilde{p}^{k+1}_h$ and $\tilde{n}^{k+1}_h$  to the constrained manifold with mass conservation and positivity preservation.
We adopt the $L^2$ projection here to enforce the positivity  and mass conservation, which reads
\begin{equation}
\begin{aligned} \label{ps2}
&\min_{{p}^{k+1}_h,{n}^{k+1}_h\in X} \quad \frac{1}{2}\left(\|{p}^{k+1}_h-\tilde{p}^{k+1}_h\|^2+\|{n}^{k+1}_h-\tilde{n}^{k+1}_h\|^2\right)\\
&\begin{array}{r@{\quad}l@{}l@{\quad}l}
\qquad \text{s.t.}    & {p}^{k+1}_h\geq 0,  {n}^{k+1}_h\geq 0,\,\langle {p}^{k+1}_h,1\rangle=\langle{p}^{0}_h,1\rangle,\quad \langle {n}^{k+1}_h,1\rangle=\langle{n}^{0}_h,1\rangle. 
\end{array}
\end{aligned}
\end{equation}
This is a convex minimization problem and can be solved efficiently. In particular, a simple semi-smooth Newton solver  through the following Karush-Kuhn-Tucker (KKT) conditions is available (see Appendix):
\begin{align}
\label{pproj}  &{p}^{k+1}_h = \tilde{p}^{k+1}_h +\lambda^{k+1}_h -\xi^{k+1},\quad {n}^{k+1}_h = \tilde{n}^{k+1}_h +\eta^{k+1}_h -\gamma^{k+1},\\
\label{projnonzero} & \lambda_{h}^{k+1}{p}^{k+1}_h=0,\quad \eta_{h}^{k+1}{n}^{k+1}_h=0,\quad \lambda_{h}^{k+1}\geq 0,\quad \eta_{h}^{k+1}\geq 0,\\
\label{projmass} & \quad \langle {p}^{k+1}_h,1\rangle=\langle {p}^{0}_h,1\rangle,\quad  \langle {n}^{k+1}_h,1\rangle=\langle {n}^{0}_h,1\rangle,
\end{align}
where $\xi^{k+1},\gamma^{k+1}\in \mathbb{R}$ are the Lagrange multipliers for the mass conservation, and $\lambda^{k+1}_h,\eta^{k+1}_h\in X$ are the Lagrange multipliers for the positivity preservation.  The idea is to use the semi-smooth Newton method  for solving the equations of the Lagrange multipliers $\xi^{k+1}$ and $\gamma^{k+1}$, where $p_h^{k+1}$ and $n_h^{k+1}$ can be obtained from \eqref{psolu}-\eqref{nsolu}. The detail is shown in Appendix.

\textbf{Step 3}: Update $\phi_{h}^{k+1}\in X_0$ by solving
\begin{equation}\label{pphi}
 -\Delta_h \phi^{k+1}_h= p^{k+1}_h- n^{k+1}_h.
\end{equation}
By initialization and the mass conservation property, we have $ p^{k+1}_h- n^{k+1}_h\in X_0$ and $\phi_h^{k+1}\in X_0$ can be uniquely determined.

Now, \eqref{pseca}-\eqref{pphi} complete the {\bf  CNFD projection} ({\bf CNFDP}) scheme.  Since \eqref{pseca}-\eqref{psecb} and \eqref{pphi} are linear  and the convex minimization problem \eqref{ps2} admits unique solutions, we find CNFDP is uniquely solvable at each time step. The semi-smooth Newton method (see Appendix) for solving \eqref{ps2} is efficient, and extensive numerical results (see Section 5) show that only one or two Newton iterations are required for the optimization problem \eqref{ps2} with  $O(N^2)$ complexity. Utilizing fast Fourier transform (FFT) to solve \eqref{pseca}-\eqref{psecb} and \eqref{pphi}, the overall memory cost of CNFDP is $O(N^2)$ and the computational cost is $O(N^2\ln N$), which is very efficient.

\begin{remark}
Since \eqref{pseca}-\eqref{psecb} is a three-level scheme,  for the first step $k=0$, we use the  first order scheme instead
\begin{equation}
\label{ps1a}\frac{\tilde{p}^{1}_{h} - p^{0}_{h} }{\tau}=\Delta_h \tilde{p}^{1}_{h}
+\nabla_h^+ \cdot((\mathcal{A}p_{h}^{0}) \nabla_h^- \phi_{h}^{0}),\;
\frac{ \tilde{n}^{1}_{h} - n^{0}_{h}  }{\tau}=\Delta_h \tilde{n}^{1}_{h} -\nabla_h^+ \cdot ((\mathcal{A}n_{h}^{0}) \nabla_h^- \phi_{h}^{0}).
\end{equation}
\end{remark}
\begin{remark}
It is not necessary to  explicitly  compute the value of $\lambda^{k+1}_h$ and  $\eta^{k+1}_h$,  since we can use the complementary slackness property in  KKT conditions \eqref{projnonzero} to determine the solution (see appendix and \cite{ChengQing1}). Projection in other norms rather than the $L^2$ norm in \eqref{ps2} is possible and we shall discuss this issue in section \ref{sec:H1discrete}.
\end{remark}
Concerning the projection part \eqref{pproj}-\eqref{projmass}, we have the following results.
\begin{lemma}
\label{lem:lampoistive}
For CNFDP \eqref{eq:ini}-\eqref{pphi},   the two Lagrange multipliers $\xi^{k+1}, \gamma^{k+1}\in\mathbb{R}$ satisfy
\begin{equation*}
\begin{aligned}
 \xi^{k+1} \geq0,  \gamma^{k+1} \geq 0,\quad k=0,1,2,\cdots.
\end{aligned}
\end{equation*}
\end{lemma}
\textbf{Proof:}
For $k\ge1$, combining \eqref{pseca} and \eqref{pproj}, we have
\begin{equation}\label{psecsum}
\frac{p^{k+1}_h -  p^{k}_h}{\tau}-\Delta_h \tilde{p}^{k+\frac{1}{2}}_h -  \nabla^+_h\cdot (\mathcal{A}{p}^{k+\frac{1}{2},*}_h  \nabla^-_h \phi^{k+\frac{1}{2},*}_h) =\frac{\lambda^{k+1}_h - \xi^{k+1} }{\tau}.
\end{equation}
Taking  inner products of  \eqref{psecsum} with 1 on both sides, noticing $\lambda^{k+1}_h\ge0$, we get
\begin{equation*}\label{psf1}
\langle\lambda^{k+1}_h- \xi^{k+1},1\rangle=0\Longrightarrow \xi^{k+1}=\frac{\langle \lambda^{k+1}_h,1\rangle}{|\Omega|}\geq 0.
\end{equation*}
The $k=0$ case is similar and $\gamma^{k+1}\geq0$ can be deduced analogously.  
$\hfill\Box$
\begin{lemma}\label{lem:solucomp}
For the projection part \eqref{pproj}-\eqref{projmass}, the following estimates hold,
\begin{equation*}
\begin{aligned}
\|p_h^{k+1}\|^2+\|p_h^{k+1}-\tilde{p}_h^{k+1}\|^2\leq \|\tilde p_h^{k+1}\|^2, \,  \|n_h^{k+1}\|^2+\|n_h^{k+1}-\tilde{n}_h^{k+1}\|^2\leq \|\tilde n_h^{k+1}\|^2,\, k\ge0.
\end{aligned}
\end{equation*}
\end{lemma}
\textbf{Proof.}
Testing both sides in \eqref{pproj} with $p_h^{k+1}$ and $n_h^{k+1}$, respectively,  we obtain
\begin{equation*}
\begin{aligned}
\frac{\|p_h^{k+1}\|^2+\|p_h^{k+1}-\tilde{p}_h^{k+1}\|^2-\|\tilde{p}_h^{k+1}\|^2}{2}=\langle\lambda_h^{k+1},p_h^{k+1}\rangle-\xi^{k+1}\langle p_h^{k+1},1\rangle,\\
\frac{\|n_h^{k+1}\|^2+\|n_h^{k+1}-\tilde{n}_h^{k+1}\|^2-\|\tilde{n}_h^{k+1}\|^2}{2}=\langle\eta_h^{k+1},n_h^{k+1}\rangle-\gamma^{k+1}\langle n_h^{k+1},1\rangle,
\end{aligned}
\end{equation*}
According to \cref{lem:lampoistive} and the KKT conditions,  we find both right-hand-sides in the above equations are non-positive, and  the conclusions in the lemma follow. $\hfill\Box$

\section{Error estimates}
\label{sec:L2erroestimate}
In this section, we carry out the error analysis for CNFDP \eqref{pseca}-\eqref{pphi} with \eqref{eq:ini} and \eqref{ps1a}. 

\subsection{Main results}
Let  $T>0$ be a fixed time, and $(p(\mathbf{x},t)\ge0,n(\mathbf{x},t)\ge0,\phi(\mathbf{x},t))$ be the exact solution of \eqref{PDE}.
Based on the theoretical results, we make the following assumptions,
\begin{equation*}
\text{(A)} \quad p(t):=p(\mathbf{x},t),\, n(t):=n(\mathbf{x},t), \,\phi(t):=\phi(\mathbf{x},t)\in C^3([0, T];C^4_p(\Omega) ),
\end{equation*}  
where $C^4_p(\Omega)=\{u\in C^4(\Omega)| \partial_{x}^k\partial_y^lu\;\text{is periodic on}\, \Omega,\, \forall k,l\ge0,\, k+l\leq 4\}$. The following error  bounds can be established.

\begin{theorem}\label{thm:main}
Under the Assumption (A), let $(p_h^k,n_h^k,\phi_h^k)\in X\times X\times X_0$ be obtained by CNFDP  \eqref{pseca}-\eqref{pphi} with \eqref{eq:ini} and \eqref{ps1a}, then there exists $\tau_0,h_0>0$, such that for $0<\tau<\tau_0$ and $0<h<h_0$ satisfying a mild CFL type condition $\tau \leq C_0h$ ($C_0>0$),  the following error  estimates hold
\begin{equation*}
\begin{aligned}
\| I_hp(t_k)-p_h^{k} \|+ \| I_hn(t_k)-n_h^{k} \|+\|\nabla_h^-(I_h\phi(t_k)-\phi_h^k)\| \leq  C(\tau^2+h^2),\, 0\leq k\leq T/\tau, 
\end{aligned}
\end{equation*}
where $C>0$ is a constant  independent of $h$, $\tau$ and $k$.  
\end{theorem}

From initialization \eqref{eq:ini},  $\langle p_h^0,1 \rangle=\langle n_h^0,1\rangle =\max\{\langle I_hp_0,1 \rangle, \langle I_hn_0,1 \rangle\}=M_0$.  It is easy to see that $I_hp(t)$ and $I_hn(t)$ do not preserve the discrete mass $M_0$ in general. For the purpose of numerical analysis, we introduce the `biased' error functions
$e_p^k,e_n^k,\tilde{e}_p^k,\tilde{e}_n^k\in X, e_\phi^k\in X_0$ ($k\ge0$):
\begin{equation}\label{eq:error}
\begin{aligned}
  &e_p^{k}=p_h(t_{k})-p^{k}_{h},\quad \tilde{e}_p^{k}=p_h(t_{k})-\tilde{p}^{k}_{h},\quad e_n^{k}=n_h(t_{k})-n^{k}_{h},\\
  & \tilde{e}_n^{k}=n_h(t_{k})-\tilde{n}^{k}_{h},   \quad e_{\phi}^{k}=\phi_h(t_{k})-\phi^{k}_{h},
\end{aligned}
\end{equation}
where $\tilde{p}_h^0=p_h^0$, $\tilde{n}_h^0=n_h^0$ and 
\begin{equation}\label{eq:ref}
p_h(t_k)=(1+\epsilon_p^k)I_hp(t_k),\quad n_h(t_k)=(1+\epsilon_n^k)I_hn(t_k),\quad \phi_h(t_k)=I_h\phi(t_k)-\epsilon_\phi^k,
\end{equation}
with
\begin{equation}
\epsilon_p^k=\frac{M_0-\langle I_hp(t_k),1\rangle}{\langle I_hp(t_k),1\rangle},\quad
\epsilon_n^k=\frac{M_0-\langle I_hn(t_k),1\rangle}{\langle I_hn(t_k),1\rangle},\quad
\epsilon_\phi^k=\frac{1}{N^2h^2}\langle I_h\phi(t_k),1\rangle.
\end{equation}
It is straightforward to check $\langle p_h(t_k),1\rangle=\langle n_h(t_k),1\rangle=M_0$ and $\phi_h(t_k)\in X_0$, which immediately implies  $\langle e_p^{k},1\rangle=\langle e_n^{k},1\rangle=\langle \tilde{e}_p^{k},1\rangle=\langle \tilde{e}_n^{k},1\rangle=0$ ($k\ge0$).  Under the assumption (A), for sufficiently small $h$, we have $p_h(t_k),n_h(t_k)\ge0$ and
\begin{equation}\label{eq:factor}
|\epsilon_p^k|+|\epsilon_n^k|+|\epsilon_\phi^k|\leq Ch^2, \quad|\epsilon_p^{k+1}-\epsilon_p^k|+|\epsilon_n^{k+1}-\epsilon_n^k|\leq C\tau h^2,\quad 0\leq k\leq T/\tau,
\end{equation}
where $C>0$ depends on $(p(t),n(t),\phi(t))$. The first part in \eqref{eq:factor} is a direct consequence of the second order trapezoidal rule (it is indeed of spectral accuracy for periodic functions). For the second part in \eqref{eq:factor}, using the error formula for the trapezoidal rule, we can get
\begin{align*}
\int_{\Omega}p(\mathbf{x},t_k)\,d\mathbf{x}-\langle I_hp(t_k),1\rangle=&\sum_{i=1}^{N}\int_{y_0}^{y_N}\int_{x_{i-1}}^{x_{i}}
\frac{(x_i-x)(x_{i-1}-x)}{2}\partial_{xx}p(x,y,t_k)\,dxdy\\
&+h\sum_{i,j=1}^N\int_{y_{j-1}}^{y_j}\frac{(y_j-y)(y_{j-1}-y)}{2}\partial_{yy}p(x_i,y,t_k)\,dy.
\end{align*}
Since $ \int_{\Omega}p(\mathbf{x},t)\,d\mathbf{x}$ is a constant in time, we have from the above error formula that $\epsilon_p^{k+1}-\epsilon_p^k=M_0\frac{\langle I_hp(t_k),1\rangle-\langle I_hp(t_{k+1}),1\rangle}{\langle I_hp(t_k),1\rangle\,\langle I_hp(t_{k+1}),1\rangle}=O(\tau h^2(\partial_{xxt}p+\partial_{yyt}p ))$ and we omit the details here for brevity.

Now, we are ready to introduce the `local truncation errors' $R_p^k,R_n^k,R_\phi^k\in X$ ($k\ge0$) as
\begin{align}\label{eq:lc1}
  R_{\phi}^{k}=&-\Delta_h \phi_h(t_{k}) - p_h(t_{k})+n_h(t_{k}),\quad k\ge0,\\ \label{eq:lc2}
R_p^0=&\frac{p_{h}(\tau) - p_{h}(0) }{\tau}-\Delta_h p_{h}(\tau)
-\nabla_h^+ \cdot((\mathcal{A}p_{h}(0)) \nabla_h^- \phi_{h}(0)),\\
R_n^0=&\frac{ n_{h}(\tau) - n_{h}(0)  }{\tau}-\Delta_h n_{h}(\tau) +\nabla_h^+ \cdot ((\mathcal{A}n_{h}(0)) \nabla_h^- \phi_{h}(0)),\label{eq:lc3}
\end{align}
and for $k\ge1$
\begin{align}\label{eq:lc4}
 R_p^{k} &=\frac{p_h(t_{k+1})-p_h(t_k)}{\tau} -\frac{\Delta_h \left(p_h(t_{k+1})+p_h(t_{k})\right)}{2}- \nabla^{+}_h \cdot \left((\mathcal{A}p_{h}^{t_k,*})\nabla^{-}_h \phi_h^{t_k,*}\right),\\
 R_n^{k} &=\frac{n_h(t_{k+1})-n_h(t_k)}{\tau} -\frac{\Delta_h \left(n_h(t_{k+1})+n_h(t_{k})\right)}{2}
+\nabla^{+}_h \cdot ((\mathcal{A}n_h^{t_{k},*}) \nabla^{-}_h \phi_h^{t_{k},*}),\label{eq:lc5}
\end{align}
where $ p_h^{t_{k},*}=\frac{3}{2}p_h(t_{k}) -\frac{1}{2} p_h(t_{k-1}), \,  n_h^{t_{k},*}=\frac{3}{2}n_h(t_{k}) -\frac{1}{2} n_h(t_{k-1}),\,
  \phi_h^{t_{k},*}=\frac{3}{2}\phi_h(t_{k}) -\frac{1}{2} \phi_h(t_{k-1})$ ($k\ge1$).
By Taylor expansion, we can obtain the following estimates on the local error.
\begin{lemma}\label{lemma:local}
For the local truncation errors $R_p^k,R_n^k,R_\phi^k\in X$ ($k\ge0$) defined in \eqref{eq:lc1}-\eqref{eq:lc5}, under the assumption (A), we have the estimates for $0\leq k\leq\frac{T}{\tau}-1$
\begin{equation} \label{eq:bd}
\|R_p^{k+1}\|+ \|R_n^{k+1}\| \leq C(\tau^2 + h^2),\,\|R_\phi^{k}\|  \leq Ch^2,\,\|R_p^0\|+\|R_n^0\|\leq C(\tau+h^2),
\end{equation}
where $C$ is independent of $\tau$, $h$ and $k$.
\end{lemma}
\textbf{Proof}.
By classic Taylor expansion approach, it is not difficult to check the bounds in \eqref{eq:bd} 
if replacing $(p_h(t),n_h(t),\phi_h(t))$ in \eqref{eq:lc1}-\eqref{eq:lc5} by $(I_hp(t),I_hn(t),I_h\phi(t))$.
When  $(p_h(t),n_h(t),\phi_h(t))$ \eqref{eq:ref} are taken into consideration, recalling \eqref{eq:factor}, we have
\begin{align*}
R_p^0=&O(\tau+h^2)+O(\frac{\epsilon_p^1-\epsilon_p^0}{\tau}p(\cdot)+\epsilon_p^1\partial_tp(\cdot)+\epsilon_p^1(\Delta p(\cdot)
+(|\epsilon_p^0|+|\epsilon_\phi^0|)\nabla\cdot(p(\cdot)\nabla\phi(\cdot)))\\
=&O(\tau+h^2),
\end{align*}
and the rest estimates can be derived similarly. $\hfill\Box$

The following lemma due to the $L^2$ projection adopted in CNFDP is useful and is the main reason to introduce the error functions as in \eqref{eq:error}. 
\begin{lemma}\label{errlo}
For the errors defined in \eqref{eq:error}, there hold
\begin{align*}
\|e_p^{k}\|^2+\|e_p^{k}-\tilde{e}_p^{k}\|^2\leq \|\tilde{e}_p^{k}\|^2, \, \|e_n^{k}\|^2+\|e_n^{k}-\tilde{e}_n^{k}\|^2\leq \|\tilde{e}_n^{k}\|^2,\,   0\leq k\leq \frac{T}{\tau}.
\end{align*}
\end{lemma}
\textbf{Proof}. $k=0$ is trivial. For $k\ge1$, we have
from \eqref{pproj} that
\begin{equation*}
e_p^{k}= \tilde{e}_p^{k}-\lambda^{k}_h+\xi^{k}.
\end{equation*}
Testing  both sides with $e_p^k$, we have
\begin{equation*}
\frac{1}{2}(\|e_p^{k}\|^2+\|e_p^{k}-\tilde{e}_p^{k}\|^2-\|\tilde{e}_p^{k}\|^2)=-\langle\lambda_h^k,e_p^{k} \rangle,
\end{equation*}
where we have used the fact that $\langle e_p^{k},\xi^{k}\rangle=0$  due to the precise mass conservation.  Using the KKT conditions and $\lambda_h^k\ge0$, we have $-\langle\lambda_h^k,e_p^{k} \rangle=-\langle\lambda_h^k,p_h(t_k)\rangle\leq0$ ($p_h(t)\ge0$) and the estimate on $e_p^k$ follows.
The case for $e_n^k$ is the same and is omitted here for brevity.
\hfill$\square$

 Subtracting  \eqref{pseca}-\eqref{psecb} and  \eqref{pphi} from \eqref{eq:lc4}, \eqref{eq:lc5} and \eqref{eq:lc1}, respectively, we obtain the error equations for $k\ge1$ as
\begin{equation}\label{Errfunsec}
\begin{aligned}
&\frac{\tilde{e}_p^{k+1}-e_p^{k}}{\tau}=\frac{1}{2}\Delta_h (\tilde{e}_p^{k+1}+{e}_p^{k} ) +T_1^k+R_p^{k},\\
&\frac{\tilde{e}_n^{k+1}-e_n^{k}}{\tau}=\frac{1}{2}\Delta_h (\tilde{e}_n^{k+1}+{e}_n^{k} ) +T_2^k+R_n^{k},\\
& -\Delta_h {e}_{\phi}^{k+1}= {e}_p^{k+1}- {e}_n^{k+1}+R_{\phi}^{k+1},
\end{aligned}
\end{equation}
where $T_1^k,T_2^k\in X$ ($k\ge1$ ) are defined as
\begin{equation}\label{eq:non}
\begin{aligned}
T_1^k& =   \nabla^{+}_h  \cdot \left((\mathcal{A}p_h^{t_k,*}) \nabla^{-}_h \phi_h^{t_k,*}-(\mathcal{A}p^{k+\frac{1}{2},*}_{h}) \nabla^{-}_h \phi^{k+\frac{1}{2},*}_{h}\right),\\
T_2^k& =- \nabla^{+}_h \cdot \left((\mathcal{A}n_h^{t_k,*})\nabla^{-}_h \phi_h^{t_k,*}-(\mathcal{A}n^{k+\frac{1}{2},*}_{h}) \nabla^{-}_h \phi^{k+\frac{1}{2},*}_{h}\right).\\
\end{aligned}
\end{equation}
For the nonlinear part \eqref{eq:non}, denoting
 \begin{equation*}
B=\max_{0\leq k\leq T/\tau}\{\|p_h(t_k)\|_{\infty}+\|n_h(t_k)\|_{\infty}+\|\nabla_h^-\phi_h(t_k)\|_{\infty}\},
\end{equation*}
where $B>0$ is well-defined under the assumption (A) and sufficiently small $0<h<h_1$ ($h_1>0$), we have the estimates below.
\begin{lemma}\label{lemma:non} 
 Assuming $\|p^m_h\|_{\infty}+\|n^m_h\|_{\infty}+\|\nabla^{-}_h\phi^m_h\|_{\infty} \leq B+1$  ($m=k,k-1$, $k\ge1$), under the assumption (A), for $T_1^k,T_2^k$ given in \eqref{eq:non} and any $f_h\in X$, we have
\begin{align*}
&|\langle T_1^k,f_h\rangle|\leq C_B(\|e_p^k\|+\|e_p^{k-1}\|+\|\nabla_h^- e_\phi^{k}\|+\|\nabla_h^- e_\phi^{k-1}\|)\|\nabla_h^-f_h\|,
\\
&|\langle T_2^k,f_h\rangle|\leq C_B(\|e_n^k\|+\|e_n^{k-1}\|+\|\nabla_h^- e_\phi^{k}\|+\|\nabla_h^- e_\phi^{k-1}\|)\|\nabla_h^-f_h\|,
\end{align*}
where $C_B$ is a constant depending only on $B$.
\end{lemma}
{\bf Proof}. We shall  show the case for $T_1^k$ only, since the proof for $T_2^k$ is the same. Recalling \eqref{eq:lc4}-\eqref{eq:lc5}, we have $p_h^{t_k,*}-p_h^{k+1/2,*}=\frac{3}{2}e_p^k-\frac{1}{2}e_p^{k-1}$, $\phi_h^{t_k,*}-\phi_h^{k+1/2,*}=\frac{3}{2}e_\phi^k-\frac{1}{2}e_\phi^{k-1}$ and
\begin{align*}
(\mathcal{A}p_h^{t_k,*}) \nabla^{-}_h \phi_h^{t_k,*}-(\mathcal{A}p^{k+\frac{1}{2},*}_{h}) \nabla^{-}_h \phi^{k+\frac{1}{2},*}_{h}=&\left(\mathcal{A}\left(\frac{3}{2}e_p^k-\frac{1}{2}e_p^{k-1}\right)\right)\nabla^-_h\phi^{k+\frac{1}{2},*}\\
&+\left(\mathcal{A}p_h^{t_k,*}\right)\nabla^-_h\left(\frac{3}{2}e_\phi^k-\frac{1}{2}e_\phi^{k-1}\right).
\end{align*}
Under the assumptions of Lemma \ref{lemma:non}, $\|\nabla^{-}_h \phi^{k+\frac{1}{2},*}_{h}\|_\infty\leq 2(B+1)$ and
$\|p_h^{t_k,*}\|_\infty\leq 2B$. The estimate on $\langle T_1^k,f_h\rangle$ then follows from the summation by parts formula as presented in Lemma \ref{lem:summation}.
$\hfill\Box$

Now, we proceed to  show the main error estimates.

{\it Proof of Theorem \ref{thm:main}.}
 We shall prove by induction that for  sufficiently small $h$ and $\tau$ satisfying $\tau\leq C_0h$ ($C_0>0$ is a constant),
\begin{equation}\label{eq:ind:1}
\|p^k_h\|_{\infty}+\|n^k_h\|_{\infty}+\|\nabla^{-}_h\phi^k_h\|_{\infty} \leq B+1, \quad 0 \leq k \leq T/\tau,
\end{equation}
and
\begin{equation}
\label{eq:ind:2}
\|e_p^k\|^2+\|e_n^k\|^2+\|\nabla_h^-e_\phi^k\|^2+\|\tilde{e}_p^k\|^2+\|\tilde{e}_n^k\|^2\leq Ce^{Ck\tau}(\tau^2+h^2)^2,\quad 0 \leq k \leq T/\tau,
\end{equation}
where $C$ (to be determined later) is independent of $\tau$, $h$ and $k$. Noticing \eqref{eq:ref} and \eqref{eq:factor},  it is easy to check that \eqref{eq:ind:2} implies the estimates in Theorem \ref{thm:main}.

{\bf Step 1}.  For $k=0$, due to the initialization  \eqref{eq:ini} and error definition \eqref{eq:error} , we  have $e_p^k=e_n^k=\tilde{e}_p^k=\tilde{e}_n^k=0$.
 And subtracting the third equation of \eqref{eq:ini} from \eqref{eq:lc1} with $k=0$,  we can get $-\Delta_he_\phi^0=R_\phi^0,$

which implies $\|\nabla_h^-e_\phi^0\|^2=\langle-\Delta_he_\phi^0,e_\phi^0\rangle=\langle R_\phi^0,e_\phi^0\rangle\leq \|R_\phi^0\|\,\|e_\phi^0\|$. The discrete Sobolev inequality $\|e_\phi\|\leq C_{\Omega}\|\nabla_h^-e_\phi^0\|$  ($C_\Omega$ only depends on $\Omega$) for $e_\phi^0\in X_0$ leads to $\|\nabla_h^-e_\phi^0\|^2\leq \|R_\phi^0\|^2\leq C_1h^4$ ($C_1>0$), and \eqref{eq:ind:2} holds.
On the other hand,  $\|\nabla_h^-e_\phi^0\|_\infty\leq h^{-1}\|\nabla^{-}_h e_\phi^0\|\leq \sqrt{C_1}h$ in view of the  inverse inequality. Then for some sufficiently small $h_2>0$, when $0<h<h_2$, we have \eqref{eq:ind:1}.

{\bf Step 2}. For $k=1$, subtracting   \eqref{pphi} and  \eqref{ps1a} from \eqref{eq:lc1}-\eqref{eq:lc2}, respectively, we have the error equations
\begin{align}\label{eq:k1}
&\frac{\tilde{e}_p^{1}-e_p^{0}}{\tau}=\Delta_h \tilde{e}_p^{1}+\nabla_h^+ \cdot \left((\mathcal{A}p_{h}(0)) \nabla_h^- e_\phi^{0}\right)+R_p^{0},\\
&\frac{\tilde{e}_n^{1}-e_n^{0}}{\tau}=\Delta_h \tilde{e}_n^{1}-\nabla_h^+ \cdot \left((\mathcal{A}n_{h}(0)) \nabla_h^- e_\phi^{0}\right)+R_n^{0},\label{eq:k2}\\
& -\Delta_h {e}_{\phi}^{1}= {e}_p^{1}- {e}_n^{1}+R_{\phi}^{1},\label{eq:k3}
\end{align}
where we have used the fact that $p_h^0=p_h(0)$ and $n_h^0=n_h(0)$.
Taking the $L^2$ inner product of \eqref{eq:k1} with $ \tilde{e}_p^{1}$, using Cauchy inequality, under the assumption (A), we have
\begin{align*}
 &\| \tilde e_p^{1} \|^2 - \|  e_p^{0} \|^2+\| \tilde e_p^{1}-  e_p^{0} \|^2  +  \tau \| \nabla_h^-\tilde e_p^{1} \|^2  =\tau \langle\tilde e_p^{1}, R_p^0 \rangle-\tau\langle(\mathcal{A}p_{h}(0)) \nabla_h^- e_\phi^{0},\nabla_h^-\tilde{e}_p^1\rangle\\
 \leq&\frac{1}{2}\|\tilde{e}_p^1\|^2+\frac{\tau^2}{2}\|R_p^0\|^2+\frac{\tau}{2}\|\nabla_h^-\tilde{e}_p^1\|^2+\frac{\tilde{C}_1}{2}\tau\|\nabla_h^-e_\phi^0\|^2,
\end{align*}
where $\tilde{C}_1>0$ is a constant depending on $p(\mathbf{x},t)$.  From Lemma \ref{lemma:local}, we get $\|\tilde{e}_p^1\|^2+\tau\|\nabla_h^-\tilde{e}_p^1\|^2\leq \tau^2\|R_p^0\|^2+\tilde{C}_1\tau\|\nabla_h^-e_\phi^0\|^2=O(\tau^2(\tau+h^2)^2+\tau h^4)$, which leads to $\|\tilde{e}_p^1\|=O(\tau^2+h^2)$. Similarly, we have $\|\tilde{e}_n^1\|=O(\tau^2+h^2)$. Lemma \ref{errlo} implies for some constant $C_2>0$
\begin{equation}\label{eq:k=1}
\|\tilde{e}_p^1\|^2+\|\tilde{e}_n^1\|^2+\|e_p^1\|^2+\|e_n^1\|^2\leq C_2 (\tau^2+h^2)^2.
\end{equation}
Following the derivation for $k=0$, the discrete elliptic equation \eqref{eq:k3} gives that
\begin{equation}
\|\nabla_h^-e_\phi^1\|\leq C_\Omega \left(\|e_p^1\|+\|e_n^1\|+\|R_\phi^1\|\right)=O(\tau^2+h^2).
\end{equation}
Combing the above results, we find for some constant $C_3>0$ that $\|\tilde{e}_p^1\|^2+\|\tilde{e}_n^1\|^2+\|e_p^1\|^2+\|e_n^1\|^2+\|\nabla_h^-e_\phi^1\|^2\leq C_3 (\tau^2+h^2)^2$, i.e. \eqref{eq:ind:2} holds. Again, the inverse inequality ensures that for $\tau\leq C_0 h$ and $0<h<h_3$ ($h_3$ is sufficiently small), we have \eqref{eq:ind:1} for $k=1$.

{\bf Step 3}. Now, assuming \eqref{eq:ind:1} and \eqref{eq:ind:2} hold for $k\leq m$ ($1\leq m\leq T/\tau-1$), we are going to prove the case  $k=m+1$.

Taking the  inner products of \eqref{Errfunsec} with $ \tilde{e}_p^{k+1}+{e}_p^{k}\in X$, $\tilde{e}_n^{k+1}+{e}_n^{k}\in X$ and ${e}_{\phi}^{k+1}\in X_0$, respectively, applying Sobolev inequality, we have the $L^2$ error growth for $k\ge1$,
\begin{align*}
&\frac{1}{\tau}\left( \| \tilde e_p^{k+1} \|^2 - \|  e_p^{k} \|^2\right) +\frac{1}{2}  \| \nabla_h^- (\tilde e_p^{k+1}+e_p^{k})  \|^2 = \langle T_1^k,  \tilde{e}_p^{k+1}+{e}_p^{k} \rangle+ \langle\tilde{e}_p^{k+1}+{e}_p^{k}, R_p^{k}\rangle,\\
& \frac{1}{\tau}\left( \| \tilde e_n^{k+1} \|^2 - \|  e_n^{k} \|^2\right) +  \frac{1}{2} \| \nabla_h^- (\tilde e_n^{k+1}+e_n^{k})  \|^2 = \langle T_2^k,  \tilde{e}_n^{k+1}+{e}_n^{k} \rangle + \langle\tilde{e}_n^{k+1}+{e}_n^{k}, R_n^{k} \rangle,
\end{align*}
and
\begin{equation}\label{eq:phik}
 \| \nabla_h^- e_{\phi}^{k+1} \| \leq C_\Omega \left(\|  e_{p}^{k+1} \| +\|  e_{n}^{k+1} \| +\|R_{\phi}^{k+1}\|\right),\quad k\ge0.
\end{equation}
Recalling $\langle \tilde{e}_p^{k+1},1\rangle=\langle e^k_p,1\rangle=0$, $\|\tilde{e}_p^{k+1}+e_p^{k}\|\leq C_\Omega \| \nabla_h^- (\tilde e_p^{k+1}+e_p^{k})  \|$ ($k\ge1$) by Sobolev inequality,  and we get from Cauchy inequality,
\begin{equation}\label{eq:1}
 |\langle\tilde{e}_p^{k+1}+{e}_p^{k}, R_p^{k}\rangle|\leq \| \tilde{e}_p^{k+1}+e_p^{k}\|\|R_p^k\|\leq\frac{1}{8}\| \nabla_h^- (\tilde e_p^{k+1}+e_p^{k}) \|^2
 +2C_\Omega^2 \|R_p^k\|^2.
\end{equation}
In view of Lemma \ref{lemma:non}, under the induction hypothesis, applying Cauchy inequality and \eqref{eq:phik}, we obtain for $k\ge1$
 \begin{align}\label{eq:2}
 |\langle T_1^{k}, \tilde{e}_p^{k+1}+{e}_p^{k} \rangle|\leq& C_B\| \nabla_h^-(\tilde{e}_p^{k+1}+e_p^{k})\| \left(\|e_p^k\|+\|e_p^{k-1}\|+\|\nabla_h^- e_\phi^k\|+\|\nabla_h^-e_\phi^{k-1}\|\right)\nonumber\\
 \leq&\frac{1}{8}\| \nabla_h^-(\tilde{e}_p^{k+1}+e_p^{k})\|^2+C_{B,\Omega}\sum\limits_{l=k,k-1}(\|e_p^l\|^2+\|e_n^l\|^2+\|R_\phi^l\|^2),
 \end{align}
 where $C_{B,\Omega}$ is a constant depending on $B$ and $\Omega$. Combining \eqref{eq:1} and \eqref{eq:2}, we obtain from the $L^2$ error growth that
 \begin{align}\label{eq:3}
 &\frac{1}{\tau}\left( \| \tilde e_p^{k+1} \|^2 - \|  e_p^{k} \|^2\right) +\frac{1}{4}  \| \nabla_h^- (\tilde e_p^{k+1}+e_p^{k})  \|^2\\
 &\leq C_3(\|R_p^k\|^2+\sum\limits_{l=k,k-1}(\|e_p^l\|^2+\|e_n^l\|^2+\|R_\phi^l\|^2),\nonumber
 \end{align}
where $C_3=\max\{2C_\Omega^2,C_{B,\Omega}\}$. Similarly, we have for $\tilde{e}_n^{k+1}$ ($k\ge1$) that
 \begin{align}\label{eq:4}
 &\frac{1}{\tau}\left( \| \tilde e_n^{k+1} \|^2 - \|  e_n^{k} \|^2\right) +\frac{1}{4}  \| \nabla_h^- (\tilde e_n^{k+1}+e_n^{k})  \|^2\\
 &\leq C_3(\|R_n^k\|^2+\sum\limits_{l=k,k-1}(\|e_p^l\|^2+\|e_n^l\|^2+\|R_\phi^l\|^2).\nonumber
 \end{align}
 Denoting $S^k= \| \tilde{e}_p^{k+1} \|^2+\|\tilde{e}_n^{k+1}\|^2+\|e_p^{k+1}\|^2+\|e_n^{k+1}\|^2$ ($k\ge1$), in view of \eqref{eq:3} and \eqref{eq:4}, recalling Lemma \ref{errlo} where $\|e_p^{l}\|^2\leq\|  \tilde{e}_p^{l} \|^2$ and $\|e_n^{l}\|^2\leq\| \tilde{e}_n^{l} \|^2$ ($l=k,k-1$), we have for $k\ge1$
  \begin{align}\label{eq:5}
 &\frac{1}{\tau}\left( S^{k+1}-S^k\right) +\frac{1}{2}  \| \nabla_h^- (\tilde e_p^{k+1}+e_p^{k})  \|^2+
 \frac{1}{2}  \| \nabla_h^- (\tilde e_n^{k+1}+e_n^{k})  \|^2\\
 &\leq 4C_3(\|R_p^k\|^2+\|R_n^k\|^2+\sum\limits_{l=k,k-1}(S^l+\|R_\phi^l\|^2).\nonumber
 \end{align}
 Summing \eqref{eq:5} together for $1,2,\cdots,k$, using the local error in Lemma \ref{lemma:local} and the estimates \eqref{eq:k=1} at the first step, for $1\leq k\leq m$, we arrive at
\begin{align}
&S^{k+1}+  \frac{\tau}{2}\sum^{k-1}_{l=1}\left( \| \nabla_h^-(\tilde e_p^{l+1}+e_p^{l})  \|^2 +   \| \nabla_h^- (\tilde e_n^{l+1}+e_n^{l})  \|^2 \right) \\
& \leq S^1+ 8C_3\tau\sum\limits_{l=1}^k(\|R_p^k\|^2+\|R_n^l\|^2+\|R_\phi^l\|^2+S^l)+4\tau C_3\|R_\phi^0\|^2\nonumber
\\& \leq 8C_3\tau\sum\limits_{l=1}^k S^l+C_4(\tau^2+h^2)^2,\nonumber
\end{align}
where $C_4>0$ is a constant independent of $\tau$, $h$ and $k$.
Discrete Gronwall inequality \cite{DisGronwalllemma1} yields for some $0<\tau<\tau_0$ ($\tau_0>0$),
\begin{equation*}
S^{k+1} \leq \exp(8C_3(k+1)\tau)C_4(\tau^2+h^2)^2,\quad 1\leq k\leq m.
\end{equation*}
From \eqref{eq:phik}, we have $\|\nabla_h^-e_\phi^{k+1}\|^2\leq C_5(S^{k+1}+h^4)$ ($C_5>0$) and
\begin{equation*}
S^{k+1}+ \|\nabla_h^-e_\phi^{k+1}\|^2 \leq (C_5+1)(C_4+1)\exp(8C_3(k+1)\tau)C_4(\tau^2+h^2)^2,\quad 1\leq k\leq m,
\end{equation*}
and \eqref{eq:ind:2} holds at $k=m+1$, if we set $C=\max\{8C_3,(C_5+1)(C_4+1),C_4\}$. It is easy to check that the constant $C$ is independent of $h$, $\tau$ and $k$. We are left with \eqref{eq:ind:1} for $k=m+1$, which can be derived similarly as the $k=0$ case by the inverse inequality. More precisely, for $k=m+1$, \eqref{eq:ind:1} implies  for $\tau\leq C_0h$ with $h<h_4$ ($h_4>0$)
\begin{equation*}
\begin{aligned}
&\|e^{m+1}_p\|_{\infty}+\|e^{m+1}_n\|_{\infty} +\|\nabla_h^- e^{m+1}_\phi\|_{\infty}\\
\leq&
h^{-1}(\|e^{m+1}_p\|+\|e^{m+1}_n\|+\| \nabla_h^-e_\phi^{m+1} \|)
\leq 3\sqrt{Ce^{CT}}(h^{-1}\tau^2+h)\leq 1,
\end{aligned}
\end{equation*}
  Triangle inequality  implies \eqref{eq:ind:1} at $k=m+1$. Choosing $h_0=\min\{h_1,h_2,h_3,h_4\}$, we then complete the induction process for proving Theorem \ref{thm:main}. $\hfill\Box$

\begin{remark}From the proof, we have the following observations.
\begin{itemize}
\item One of the key part of the error analysis is from \eqref{eq:3}-\eqref{eq:4} to \eqref{eq:5}, which is guaranteed by the $L^2$ projection \eqref{ps2} (Lemma \ref{errlo}). Lemma \ref{errlo} enables the error analysis without tracking the Lagrange multipliers. In a word, the projection step does not ruin the convergence rates.
\item For the electric potential, the error satisfies $\|\Delta_h(I_h\phi(t_k)-\phi_h^k)\|\leq C(\tau^2+h^2)$, and the $l^\infty$ error estimates hold in view of the discrete Sobolev inequality $\|I_h\phi(t_k)-\phi_h^k\|_{\infty}\leq C(\tau^2+h^2)$.
\item  For the intermediate solutions,  $\|I_hp(t_k)-\tilde{p}_h^k\|+\|I_hn(t_k)-\tilde{n}_h^k\|\leq C(\tau^2+h^2)$.
\item For the 3D case, using the  inverse inequality in 3D ($\|u_h\|_\infty\leq h^{-3/2}\|u_h\|$), the analysis works  for $ \tau \leq C_0 h $ and sufficiently small $h$, since the scheme is of second order in both space and time.
\end{itemize}
\end{remark}

\section{Different projections}
\label{sec:H1discrete}
Apart from the $L^2$ projection employed in CNFDP \eqref{pseca}-\eqref{pphi},  the correction step \eqref{ps2} can be performed in other suitable norms, e.g. the discrete $H_1$ norm. Here, we briefly discuss the $H^1$ projection strategy, which would easily yield optimal convergence in $H^1$ norm.  In detail, the projection step \eqref{ps2} is replaced by the following $H^1$ projection: given $\tilde{p}_h^{k+1},\tilde{n}_h^{k+1}\in X$, find $p_h^{k+1},n_h^{k+1}\in X$ such that
\begin{equation}
\begin{aligned} \label{h1min}
&\min_{{p}^{k+1}_h,{n}^{k+1}_h} \quad \frac{1}{2} \left( \|{p}^{k+1}_h-\tilde{p}^{k+1}_h\|_{H_1}^2+\|{n}^{k+1}_h-\tilde{n}^{k+1}_h\|_{H_1}^2 \right) \\
&\begin{array}{r@{\quad}l@{}l@{\quad}l}
\qquad \text{s.t.}  & {p}^{k+1}_h\geq 0,  {n}^{k+1}_h\geq 0,\,\langle {p}^{k+1}_h,1\rangle=\langle{p}^{0}_h,1\rangle, \,\langle {n}^{k+1}_h,1\rangle=\langle{n}^{0}_h,1\rangle.   
\end{array}
\end{aligned}
\end{equation}

The convex minimization problem \eqref{h1min} can be similarly solved as \eqref{ps2}.
The KKT conditions for \eqref{h1min} become
\begin{equation}\label{hproj}
\begin{aligned}
  &(I-\Delta_h){p}^{k+1}_h= (I-\Delta_h)\tilde{p}^{k+1}_h+\lambda^{k+1}_h-\xi^{k+1},\\
   & (I-\Delta_h){n}^{k+1}_h= (I-\Delta_h)\tilde{n}^{k+1}_h+\eta^{k+1}_h-\gamma^{k+1},
\end{aligned}
\end{equation}
where  the Lagrange multiplies $\xi^{k+1},\gamma^{k+1}\in \mathbb{R}$, $\lambda^{k+1}_h,\eta^{k+1}_h\in X$ satisfy
\begin{align}
&\label{hpproj}\lambda_{h}^{k+1}{p}^{k+1}_h=0,\quad \eta_{h}^{k+1}{n}^{k+1}_h=0,\quad \lambda_{h}^{k+1}\geq 0,\quad \eta_{h}^{k+1}\geq 0,
\end{align}
and $I$ is the identity operator.
 In particular, an efficient approach based on the semi-smooth Newton method to solve \eqref{hproj}-\eqref{hpproj} is available (see Appendix).  Extensive numerical results (see Section 5) indicate that   the optimization problem \eqref{h1min} can be solved by the semi-smooth Newton method within several Newton iterations (serval matrix-vector products are required per iteration step). Therefore, the above $H^1$ projection approach is less efficient than the $L^2$ projection in terms of the computational complexity.

We shall denote the corresponding numerical scheme \eqref{eq:ini}-\eqref{psecb} with the $H^1$ correction step \eqref{h1min}-\eqref{hpproj} and the electric potential step \eqref{pphi} as {\bf CNFDP2}.
Similar to CNFDP with the $L^2$ projection, we also have the following lemmas regarding CNFDP2 with the $H^1$ projection.
\begin{lemma}
For  CNFDP2 \eqref{eq:ini}-\eqref{psecb} with  \eqref{h1min}-\eqref{hpproj} and\eqref{pphi}, we have for the $H^1$ projection strategy \eqref{hproj},  
\begin{equation*}
\begin{aligned}
 \xi^{k+1}\geq0,  \gamma^{k+1}\geq 0,\quad k\ge0.
\end{aligned}
\end{equation*}
\end{lemma}
By testing both sides of \eqref{hproj} by $p_h^{k+1}$ and $n_h^{k+1}$ respectively, we can obtain the following estimates from the KKT conditions \eqref{hpproj}.
\begin{lemma}
For the $H^1$ projection step \eqref{hproj} in CNFDP2, there holds for $k\ge0$,
\begin{equation*}
\|p_h^{k+1}\|_{H^1}^2+\|p_h^{k+1}-\tilde{p}_h^{k+1}\|_{H^1}^2\leq \|\tilde p_h^{k+1}\|_{H^1}^2, \quad  \|n_h^{k+1}\|_{H^1}^2+\|n_h^{k+1}-\tilde{n}_h^{k+1}\|_{H^1}^2\leq \|\tilde n_h^{k+1}\|_{H^1}^2.
\end{equation*}
\end{lemma}
Following the same notations used in CNFDP,  we can obtain the following  key properties of the $H^1$ projection. 
\begin{lemma}
For  CNFDP2 \eqref{eq:ini}-\eqref{psecb} with  \eqref{h1min}-\eqref{hpproj} and \eqref{pphi}, we have for $k\ge1$,
\begin{align}\label{hep}
\|e_p^{k}\|_{H^1}^2+\|e_p^k-\tilde{e}_p^{k}\|_{H^1}^2  \leq \|\tilde{e}_p^{k}\|_{H^1}^2, \quad \|e_n^{k}\|_{H^1}^2 + \|e_n^{k}-\tilde{e}_n^k\|_{H^1}^2   \leq \|\tilde{e}_n^{k}\|_{H^1}^2.
\end{align}
\end{lemma}
Following the similar procedure  for the error analysis of CNFDP, we can obtain the estimates for CNFDP2 in $H^1$ norm.
\begin{theorem} Assuming the exact solution $p(t), n(t), \phi(t)\in C^3([0, T];C^5_p(\Omega) )$,
 there exist $\tau_0>0$ and $h_0>0$ such that when $0<\tau<\tau_0$ and $0<h<h_0$ satisfying $\tau\leq C_0h$ ($C_0>0$), the following error  estimates hold for CNFDP2,
\begin{equation*}
\begin{aligned}
\| I_hp(t_k)-p_h^{k} \|_{H^1}+ \| I_hn(t_k)-n_h^{k} \|_{H^1}+\|\nabla_h^-(I_h\phi(t_k)-\phi_h^k)\|_{H^1} \leq  C(\tau^2+h^2), 
\end{aligned}
\end{equation*}
where $0\leq k\leq \frac{T}{\tau}$, $C>0$ is a constant  independent of $h$, $\tau$ and $k$.  
\end{theorem}
The $H^1$ projection property  \eqref{hep} is crucial in the error analysis, and we shall omit the details of the proof for the above theorem. 

 \section{Numerical experiments}
\label{sec:example}
In this section, we present some numerical examples to verify the convergence and show the performance of CNFDP/CNFDP2. The fast Fourier transform (FFT) approach is utilized to improve the computational efficiency \cite{He2019APP}.
\subsection{Convergence test}
The numerical schemes CNFDP/CNFDP2 can be easily extended to deal with the problems with mass conservation and positivity preserving properties. For the convergence tests, we choose the following  problem with available exact solutions.
\begin{example}\label{example1}
We consider the following  PNP system \eqref{PDE} with source terms on domain $\Omega=[0,1]^2$:
\begin{equation*}
\begin{aligned}
&\frac{\partial p_e}{\partial t}= \Delta p_e +\nabla\cdot( p_e \nabla \phi_e )+f_p, \,
\frac{\partial n_e}{\partial t} = \Delta n_e-\nabla\cdot( n_e \nabla \phi_e )+f_n,\,-\Delta \phi_e  = p_e - n_e,
\end{aligned}
\end{equation*}
which admits the exact solution as
\begin{equation*}
\begin{aligned}
&p_e(x,y,t) =  \cos^2(\pi t +\pi  x + \pi y),\quad n_e(x,y,t) = \cos^2(\pi t + \pi x - \pi y),\\
&\phi_e(x,y,t)  =\left( \cos(2\pi t +2\pi x + 2\pi y)-\cos(2\pi t +2\pi x - 2\pi y)\right)/(16\pi^2) .
\end{aligned}
\end{equation*}
The periodic boundary conditions are assumed and the source terms ($f_n$, $f_p$) can be computed from the above exact solution.
\end{example}
Let $p_{\tau,h}(\cdot, \cdot, T), n_{\tau,h}(\cdot, \cdot, T),\phi_{\tau,h}(\cdot, \cdot, T)$ be the numerical solution at time $t=T$ with time step size $\tau$ and mesh size $h$, and the error function is denoted by, e.g.  $p_{\tau,h}-p_e=p_{\tau,h}(\cdot,\cdot,T)-I_hp_e(\cdot,\cdot,T)$.  
For the temporal convergence test, we choose a sufficiently small $h$ such that the spatial error can be neglected,
while for the spatial convergence  test,  we choose a sufficiently small $\tau$ such that the temporal error is negligible.
The final time is $T=2$.  

For the temporal error analysis, a very fine mesh size $h=1/512$  is adopted, and the  errors are displayed  in Tabs.
 \ref{tab:periodictime} and \ref{tab:periodicHnormtime} for  CNFDP and CNFDP2, respectively.

 For the spatial error analysis, we take a small time step size  $ \tau= 2\times 10^{-6}$ for  Tabs. \ref{tab:periodicspace}   
\& \ref{tab:periodicHnormspace1},  depicting the spatial convergence for CNFDP and CNFDP2, respectively.
 Second order convergence rates are clearly observed in both space and time  for CNFDP ($L^2$ norm) and  for CNFDP2 ($H^1$ norm). 

\begin{figure}[!h]
 \begin{minipage}{0.3\linewidth}
  \centering
  \includegraphics[width=0.9\linewidth]{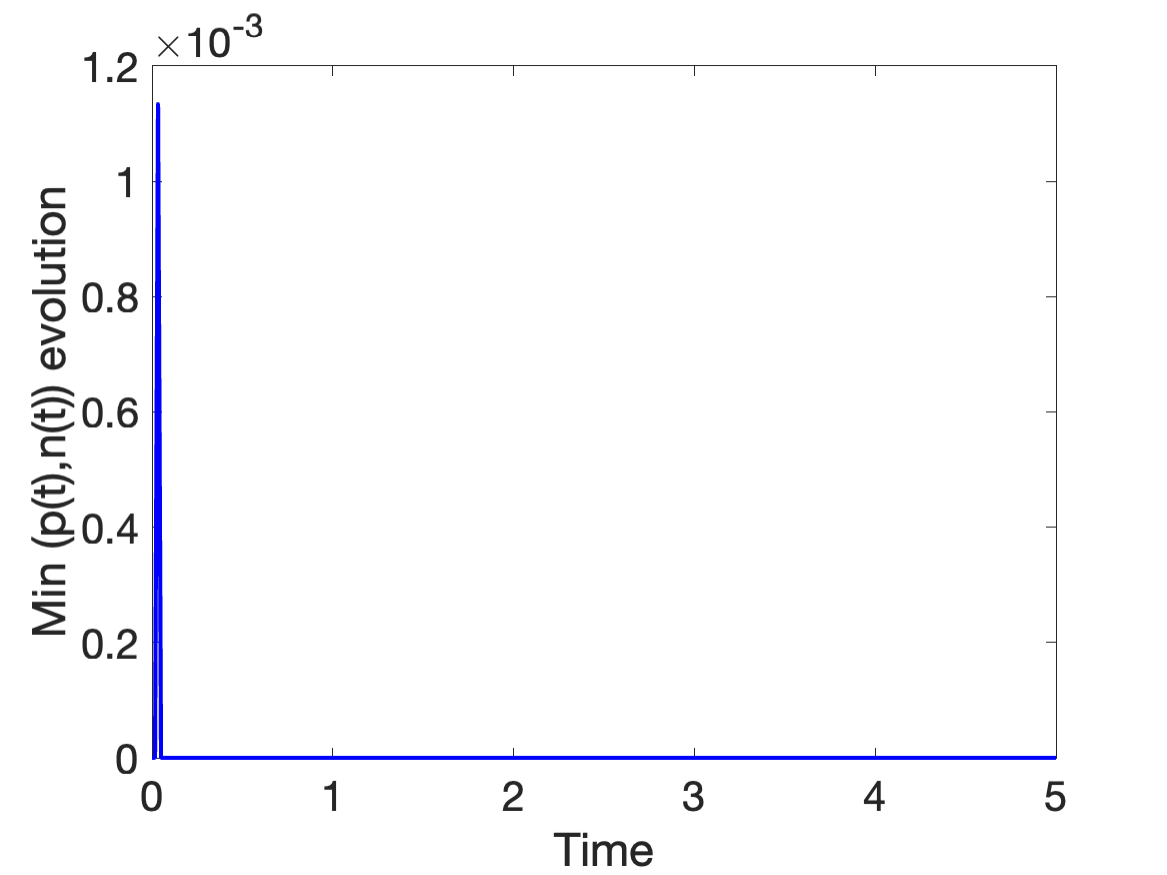}
 \end{minipage}
 \begin{minipage}{0.3\linewidth}
  \centering
  \includegraphics[width=0.9\linewidth]{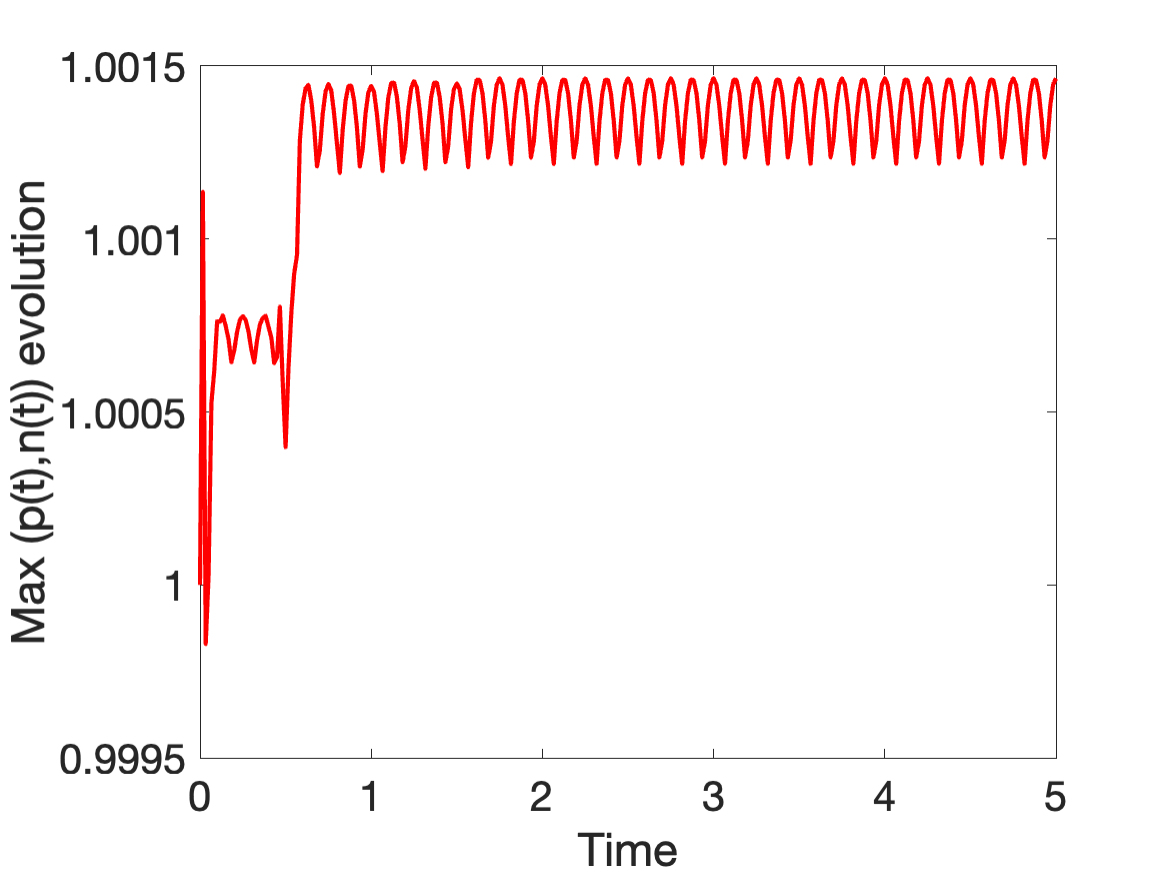}
 \end{minipage}
 \begin{minipage}{0.3\linewidth}
  \centering
  \includegraphics[width=0.9\linewidth]{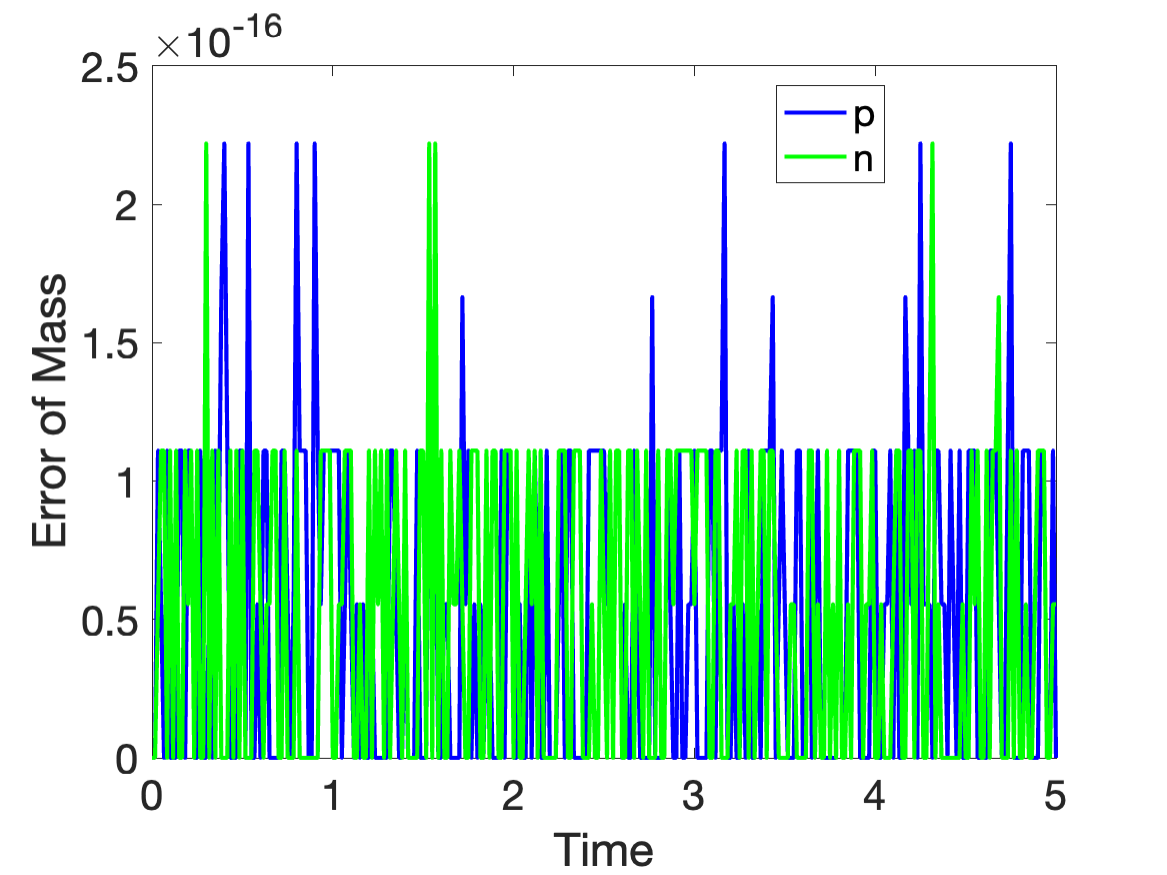}
 \end{minipage}
   \caption{ (\cref{example1} with CNFDP) Left: Lower  bounds of $(n,p)$.  Middle: Upper bounds of $(n,p)$.   Right: Conservation of discrete  masses  for $p$ and $n$. }
   \label{fig:ex1CNFDP}
\end{figure}
\begin{figure}[!h]
 \begin{minipage}{0.3\linewidth}
  \centering
  \includegraphics[width=0.9\linewidth]{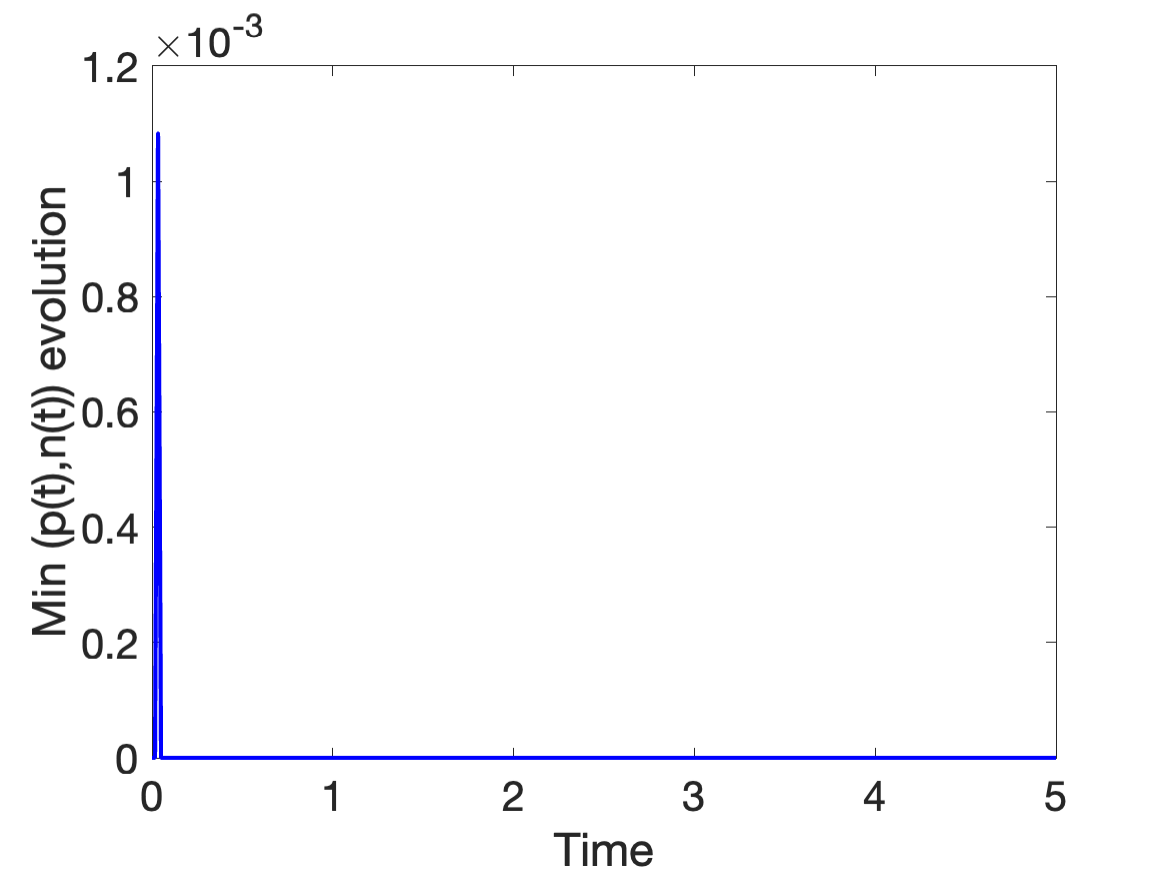}
 \end{minipage}
 \begin{minipage}{0.3\linewidth}
  \centering
  \includegraphics[width=0.9\linewidth]{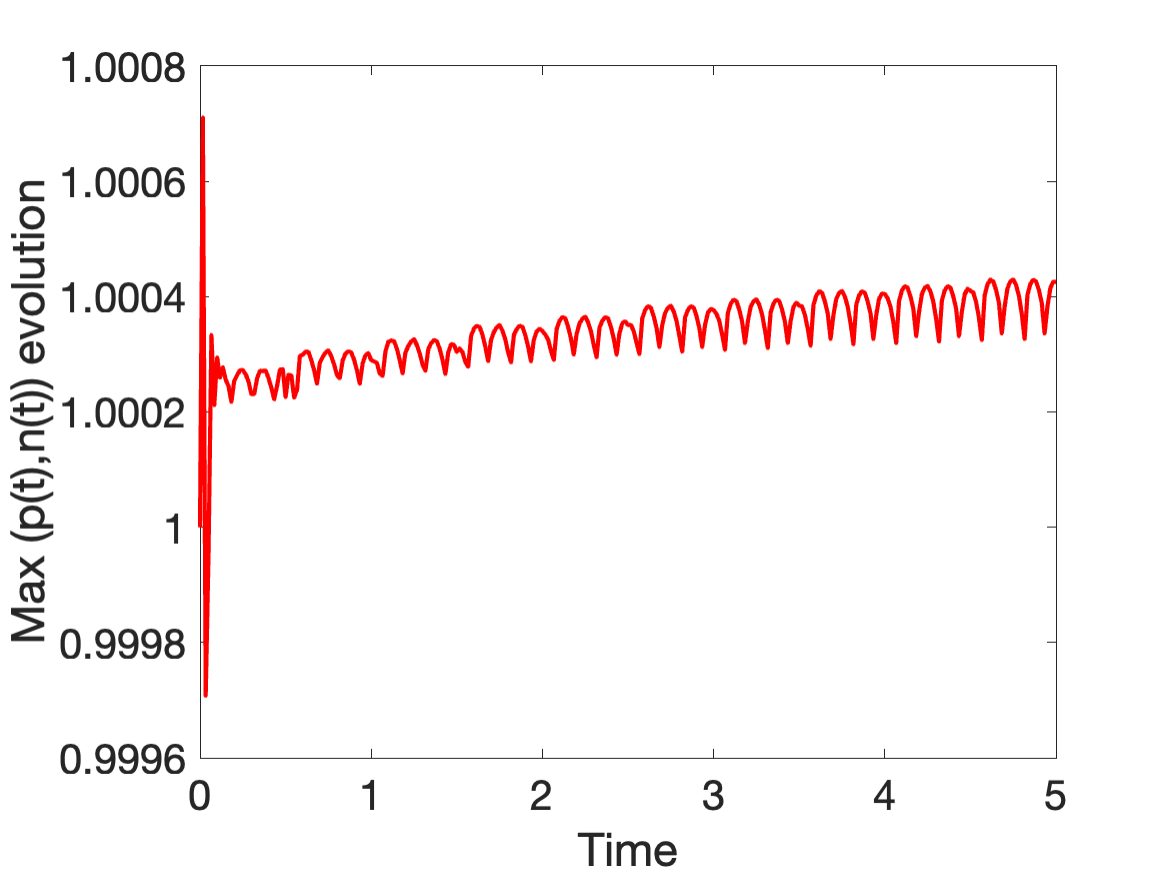}
 \end{minipage}
 \begin{minipage}{0.3\linewidth}
  \centering
  \includegraphics[width=0.9\linewidth]{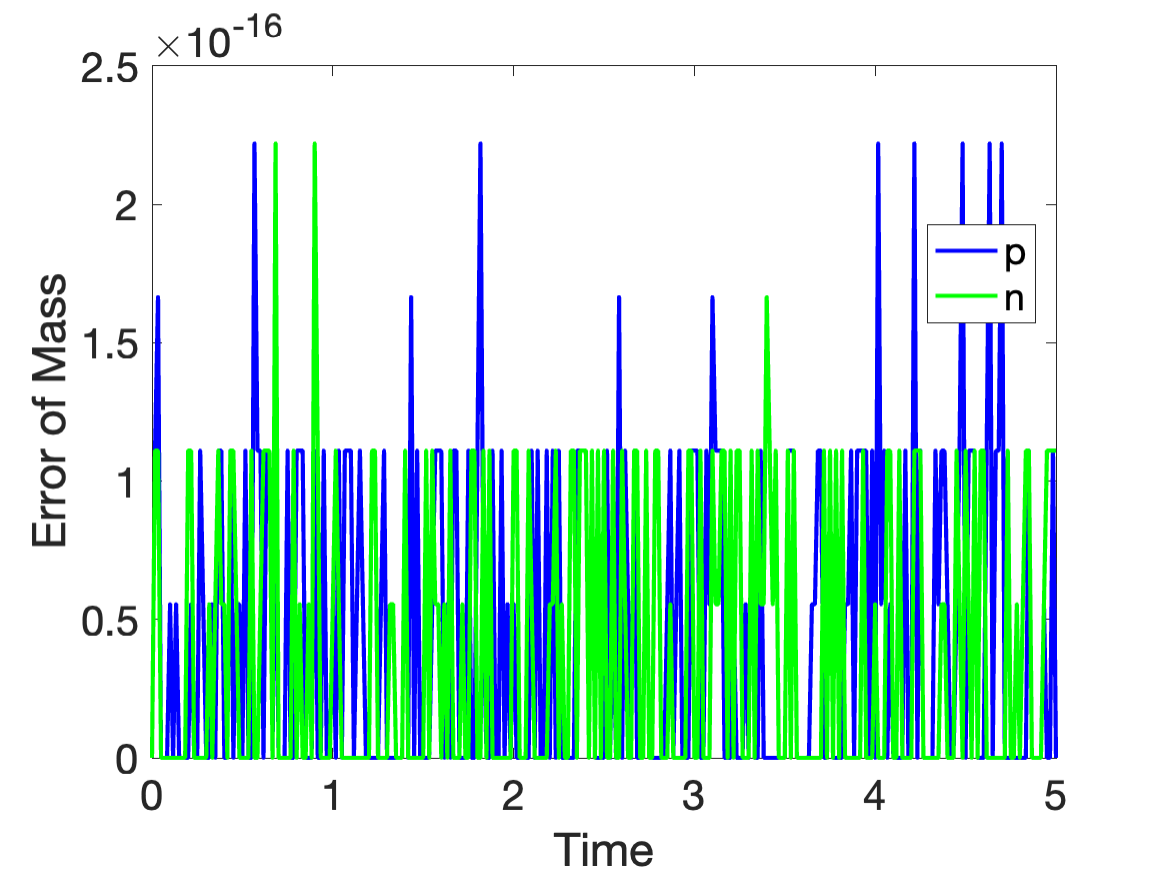}
 \end{minipage}
   \caption{ (\cref{example1} with CNFDP2) Left: Lower  bounds of $(n,p)$.  Middle: Upper bounds of $(n,p)$.   Right: Conservation of discrete  masses  for $p$ and $n$. }
   \label{fig:ex1CNFDP2}
\end{figure}
\begin{figure}[!htp]
 \centering
 \begin{minipage}{0.3\linewidth}
  \centering
  \includegraphics[width=0.9\linewidth]{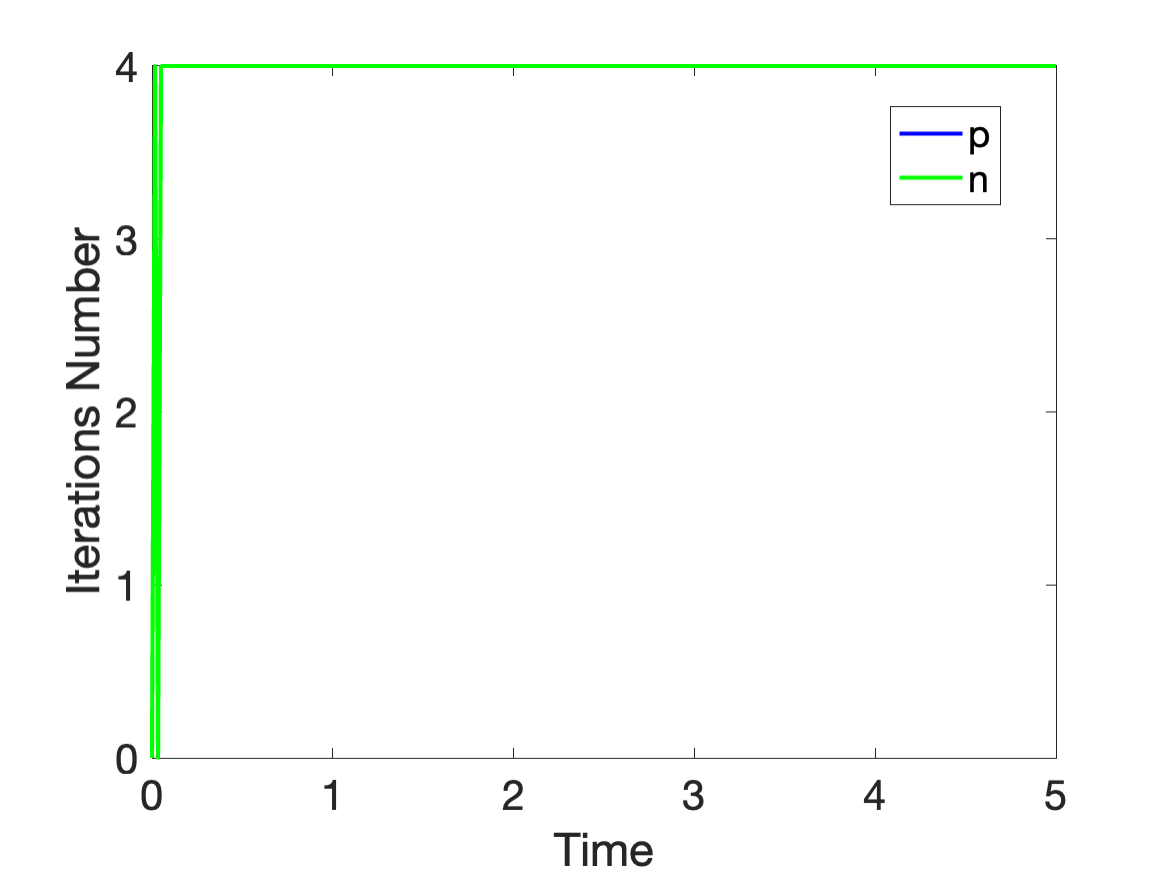}
 \end{minipage}
 \begin{minipage}{0.3\linewidth}
  \centering
  \includegraphics[width=0.9\linewidth]{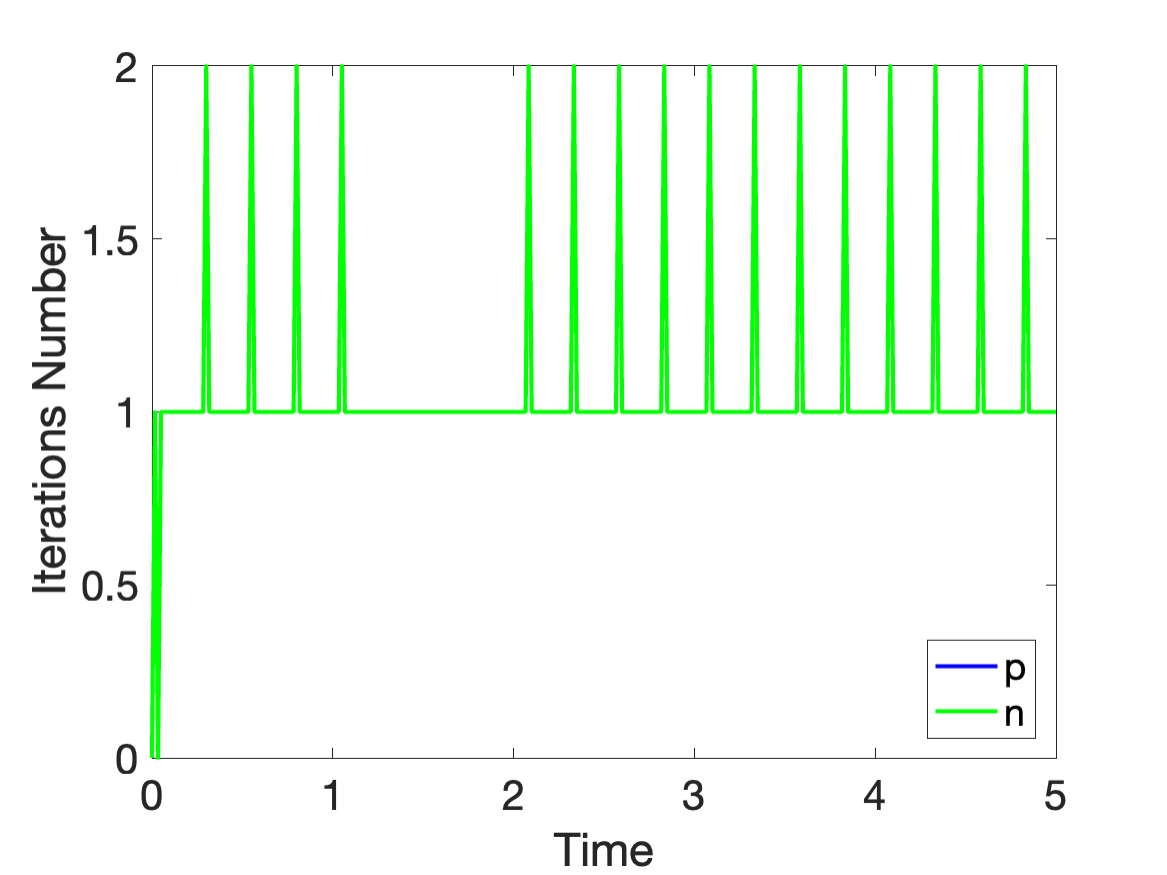}
 \end{minipage}
 \begin{minipage}{0.3\linewidth}
  \centering
  \includegraphics[width=0.9\linewidth]{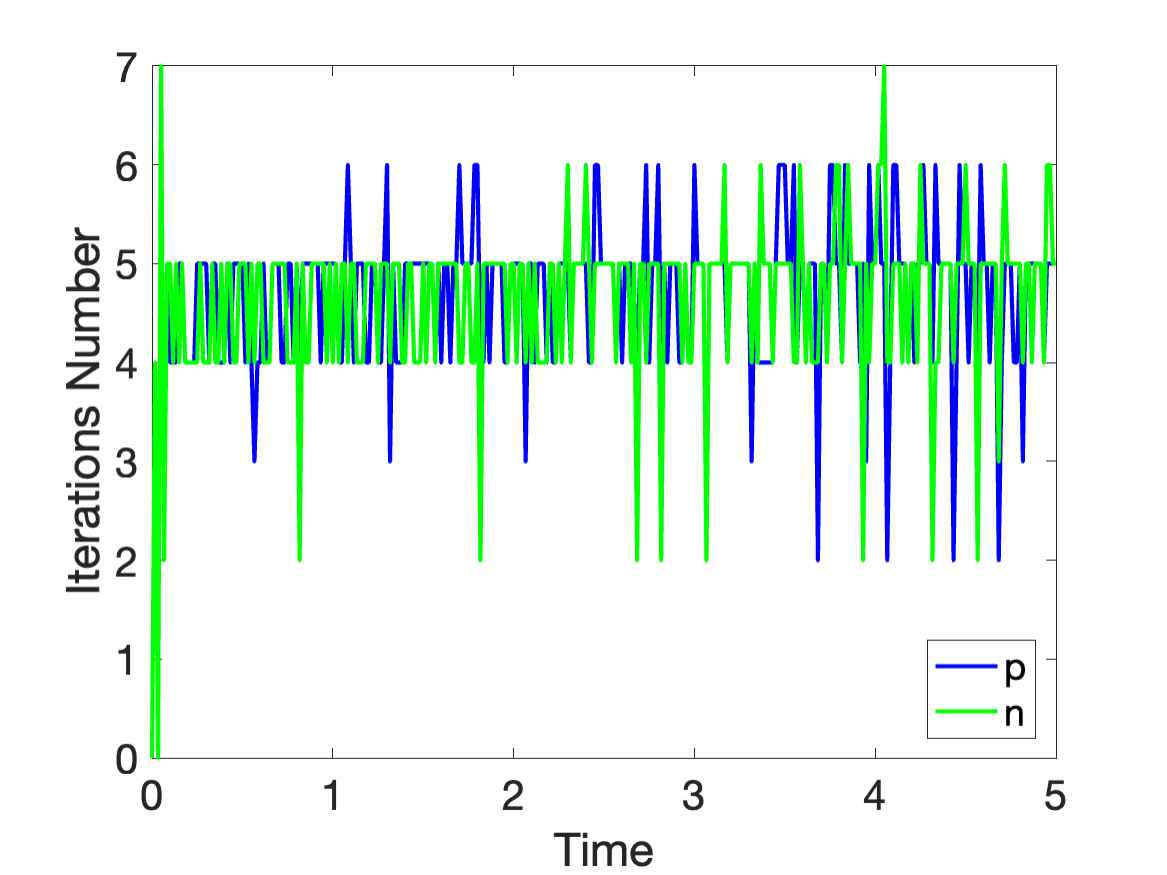}
 \end{minipage}
   \caption{Number of iterations in the  KKT condition solvers for Example \ref{example1}. Left: the secant method with CNFDP. Middle: the semi smooth Newton method with CNFDP.  Right:  the semi smooth Newton method with CNFDP2.}
  \label{fig:EX2Diter}
\end{figure}

 Figs. \ref{fig:ex1CNFDP} and \ref{fig:ex1CNFDP2}   show the  point-wise bounds of $(p,n)$ and the discrete masses, simulated  by CNFDP and CNFDP2 separately ($\tau=5/300$, $h=2^{-7}$). 
 From these figures, we observe that the numerical results   are consistent with the theoretical analysis and satisfy the desired physical properties, i.e. positivity preserving and mass conservation.  Fig. \ref{fig:EX2Diter}  presents the number of iterations of the secant method and the semi smooth Newton method for solving the mass preserving Lagrange multiplier $\xi$ and $\gamma$, which demonstrates the effectiveness of the proposed numerical scheme.  

\begin{table}[!htp]
\footnotesize
\caption{Temporal error analysis of CNFDP for  \cref{example1}. }\label{tab:periodictime}
\begin{center}
\begin{tabular}{c|cc|cc|cc}\hline
$\tau$ &  $\| p_{\tau,h} -p_{e}\|$ & Rate &$\|n_{\tau,h} -n_{e}\|$ & Rate &$\|\phi_{\tau,h} -\phi _{e}\|$ & Rate\\
\hline
 $T/2^4$   & 4.402e-02   & & 4.402e-02   & & 4.868e-04   &\\ 
 $T/2^5$  & 9.767e-03   &2.17 & 9.767e-03   &2.17 & 1.172e-04   &2.05\\ 
$T/2^6$   & 2.302e-03   &2.08  & 2.302e-03   &2.08 & 2.895e-05   &2.02\\ 
$T/2^7$  & 5.735e-04   &2.01  & 5.735e-04   &2.01 & 7.350e-06   &1.98\\ 
$T/2^8$ & 1.477e-04   &1.96 & 1.477e-04   &1.96 & 1.955e-06   &1.91\\
\hline
\end{tabular}
\end{center}
\end{table}
\begin{table}[!htp]
\footnotesize
\caption{Spatial error analysis for \cref{example1} with CNFDP. }\label{tab:periodicspace}
\begin{center}
\begin{tabular}{c|cc|cc|cc}\hline
$h$ &  $\|p_{\tau,h} -p_{e}\|$ & Rate &$\|n_{\tau,h} -n_{e}\|$ & Rate &$\|\phi_{\tau,h} -\phi_{e}\|$ & Rate\\
\hline
 $1/2^4 $ & 2.200e-03   &   & 2.200e-03   &  & 1.174e-04   & \\ 
  $1/2^5 $ & 6.107e-04   &1.85  & 6.107e-04   &1.85 & 3.030e-05   &1.95\\ 
   $1/2^6 $  & 1.622e-04   &1.91  & 1.622e-04   &1.91 & 7.756e-06   &1.97\\ 
  $1/2^7 $  & 4.081e-05   &1.99  & 4.081e-05   &1.99 & 1.953e-06   &1.99\\ 
 $1/2^8 $  & 1.004e-05   &2.02  & 1.004e-05   &2.02 & 4.875e-07   &2.00\\ 
\hline
\end{tabular}
\end{center}
\end{table}
\begin{table}[!h]
\footnotesize
\caption{Temporal  error analysis for \cref{example1} with CNFDP2. }\label{tab:periodicHnormtime}
\begin{center}
\begin{tabular}{c|cc|cc|cc}\hline
$\tau$ &  $\|p_{\tau,h} -p_{e}\|_{H^1}$ & Rate &$\|n_{\tau,h} -n_{e}\|_{H^1}$ & Rate &$\|\phi_{\tau,h} -\phi_{e}\|_{H^1}$ & Rate\\
\hline
 $T/2^4$  & 6.777e-01   &   & 6.776e-01   &  & 3.153e-03   &\\ 
  $T/2^5$  & 1.694e-01   &2.00  & 1.693e-01   &2.00 & 5.471e-04   &2.53\\ 
  $T/2^6$   & 4.007e-02   &2.08  & 4.007e-02   &2.08 & 9.723e-05   &2.49\\ 
 $T/2^7$   & 8.116e-03   &2.30  & 8.116e-03   &2.30 & 1.705e-05   &2.51\\ 
 $T/2^8$   & 1.422e-03   &2.51  & 1.422e-03   &2.51 & 3.714e-06   &2.20\\
\hline
\end{tabular}
\end{center}
\end{table}
\begin{table}[!htp]
\footnotesize
\caption{ Spatial error analysis for \cref{example1} with CNFDP2. }\label{tab:periodicHnormspace1}
\begin{center}
\begin{tabular}{c|cc|cc|cc}\hline
$h$ &  $\|p_{\tau,h} -p_{e}\|_{H^1}$ & Rate &$\|n_{\tau,h} -n_{e}\|_{H^1}$ & Rate &$\|\phi_{\tau,h} -\phi_{e}\|_{H^1}$ & Rate\\
\hline 
$ 1/2^4 $  & 7.971e-03   &  & 7.971e-03   & & 8.082e-04   & \\ 
  $ 1/2^5 $ & 2.118e-03   &1.91  & 2.118e-03   &1.91 & 2.025e-04   &2.00\\ 
$ 1/2^6 $  & 5.839e-04   &1.86  & 5.839e-04   &1.86 & 5.121e-05   &1.98\\ 
  $ 1/2^7 $  & 1.436e-04   &2.02  & 1.436e-04   &2.02 & 1.282e-05   &2.00\\
$ 1/2^8 $  & 3.350e-05   &2.10  & 3.350e-05   &2.10 & 3.205e-06   &2.00\\ 
\hline
\end{tabular}
\end{center}
\end{table}


\subsection{Extension to homogeneous Neumann boundary conditions}
\begin{example}\label{example2}
 We consider PNP \eqref{PDE}  posed on $\Omega=[-2,2]^2$ with homogeneous Neumann boundary conditions and initial conditions
\begin{equation*}
(p(x,y,0),n(x,y,0)) =\left\{\begin{array}{ll}
(1,0.5), & (x-0.5)^2+(y-0.5)^2 \leq 0.25,\\
(0.5,1), & (x+0.5)^2+(y+0.5)^2 \leq 0.25,\\
(0,0), & \text { otherwise}.\end{array}\right.
\end{equation*}
\end{example}
Here, CNFDP  is extended to investigate Example \ref{example2} with $\tau=4/200$ and $h=4/128$, where the finite difference treatment on the Neumann boundary conditions is done following \cite{CAI2022,He2019APP}. 

\begin{figure}[!htp]
\centering{
  \includegraphics[width=0.25\textwidth]{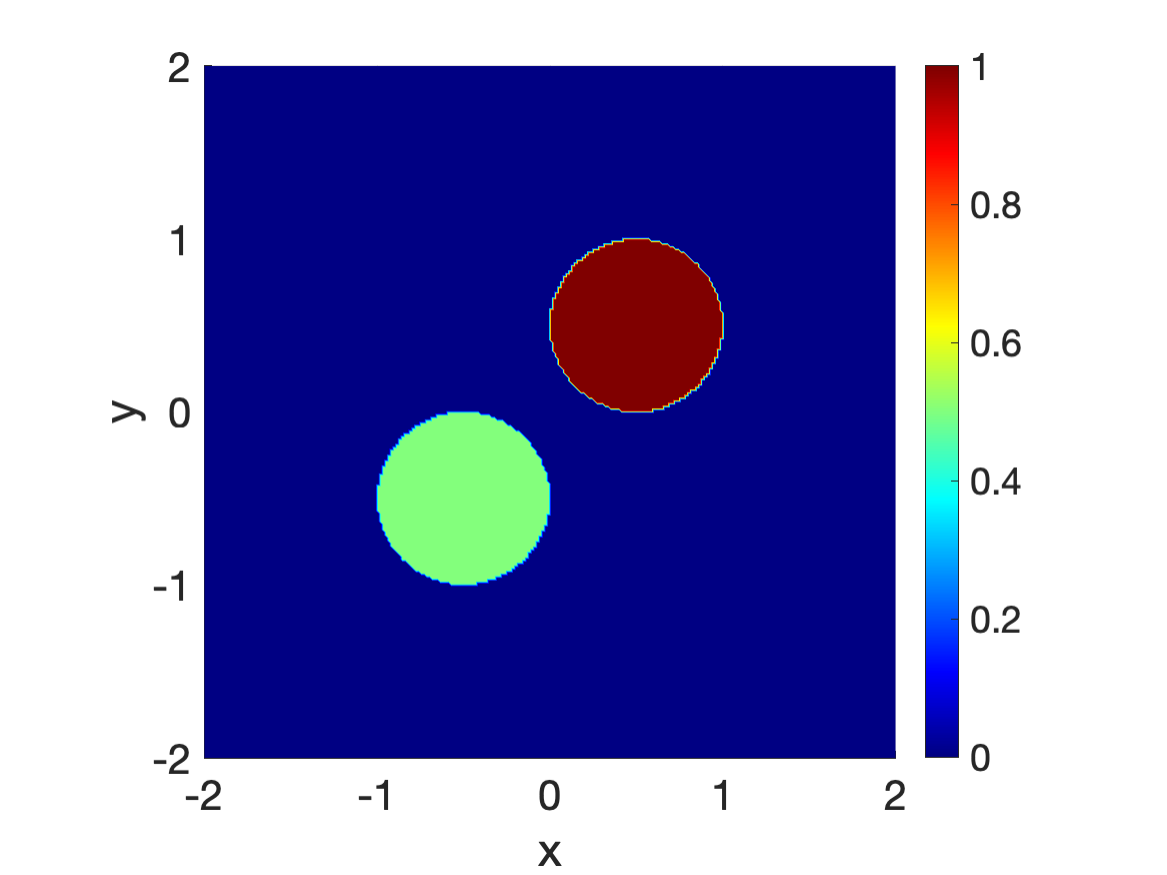}\hspace{-12.5pt}
 \includegraphics[width=0.25\textwidth]{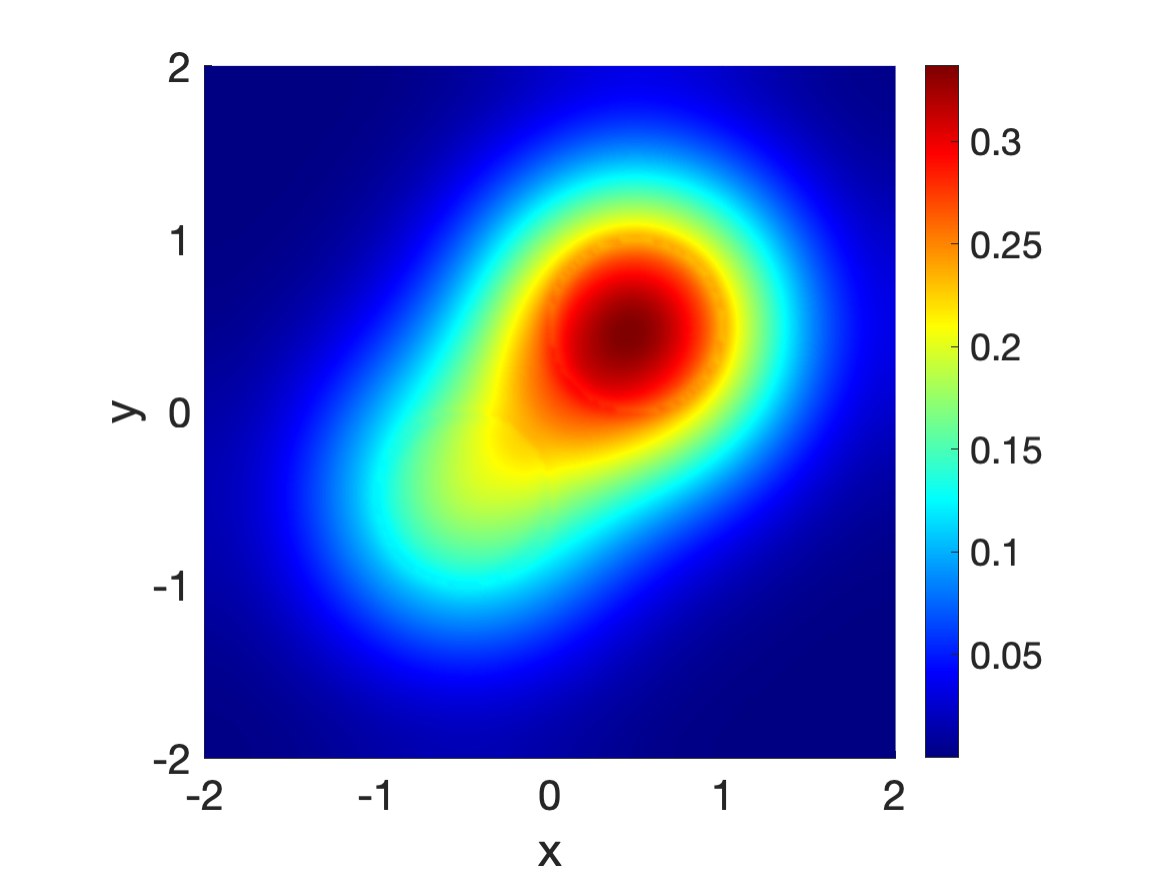}\hspace{-8.5pt}
 \includegraphics[width=0.25\textwidth]{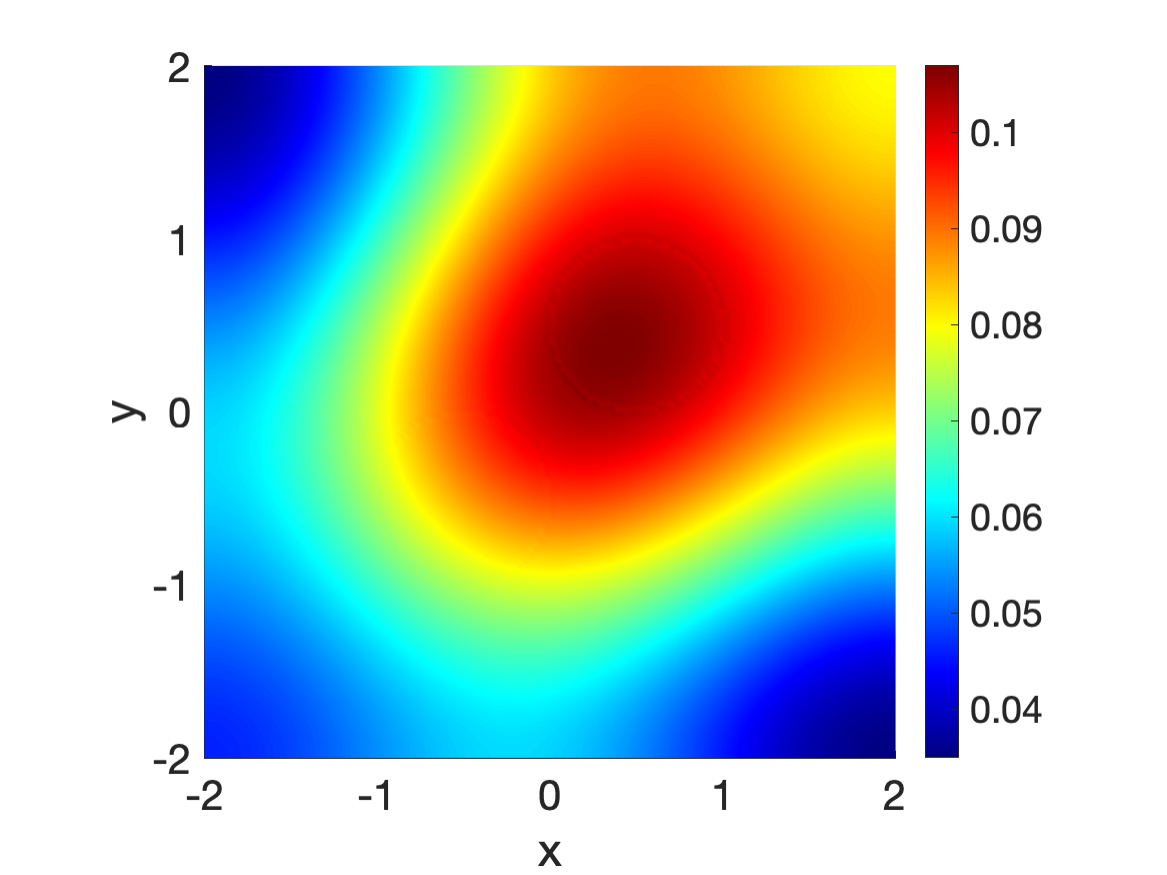}\hspace{-7.5pt}
 \includegraphics[width=0.25\textwidth]{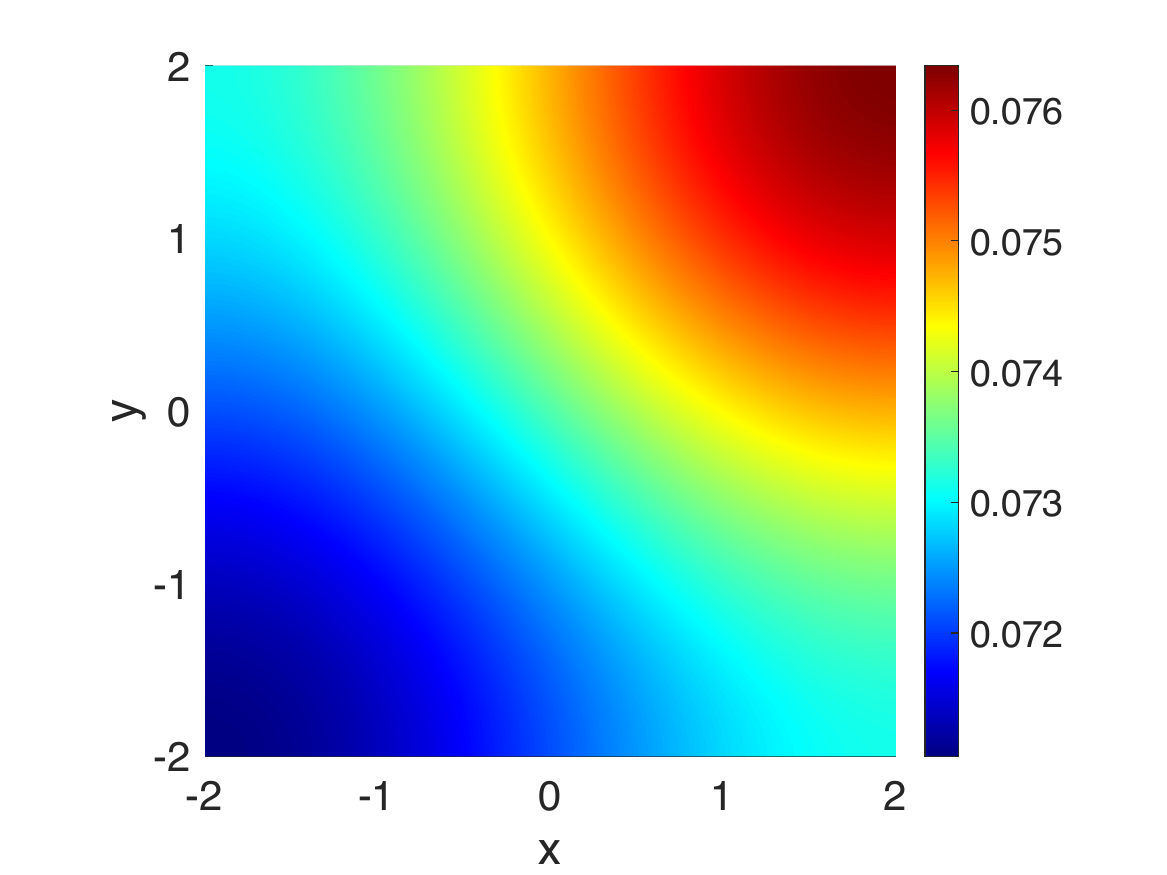}}
\centering{
  \includegraphics[width=0.25\textwidth]{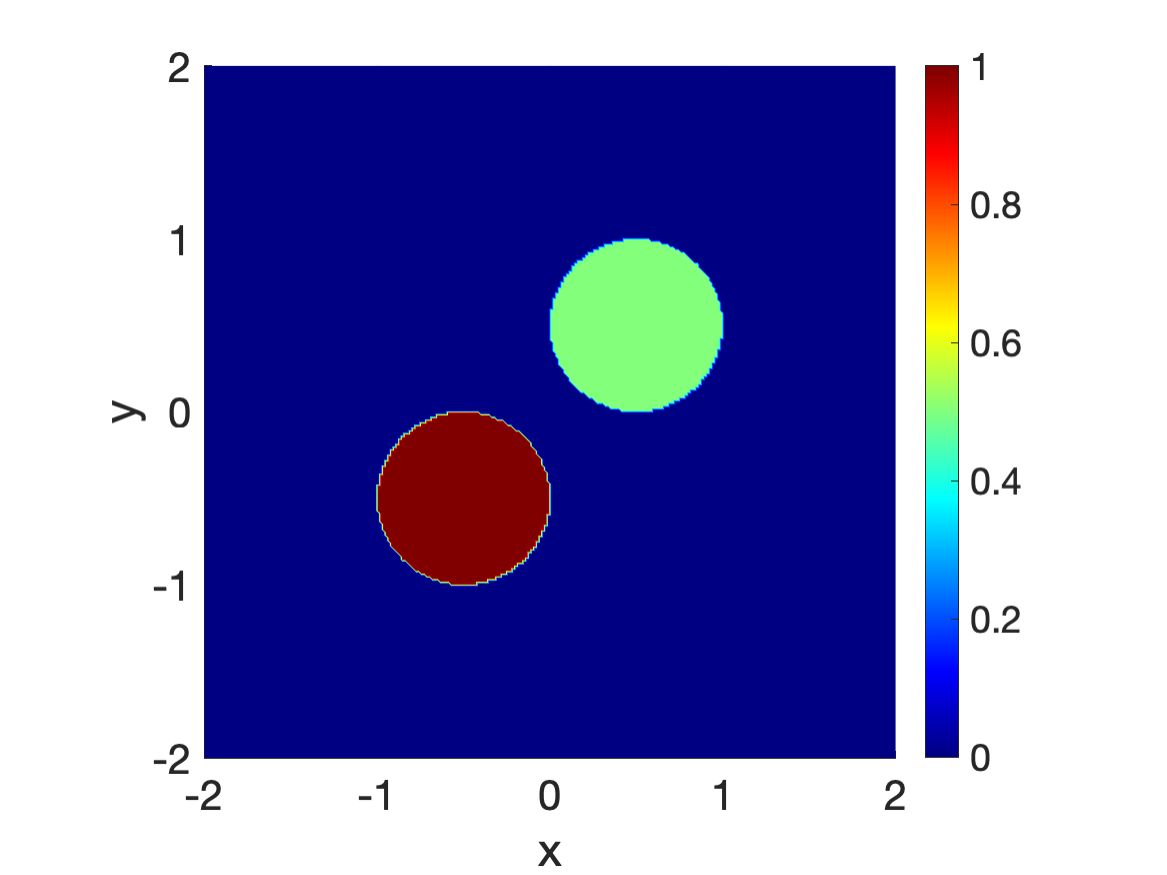}\hspace{-12.5pt}
 \includegraphics[width=0.25\textwidth]{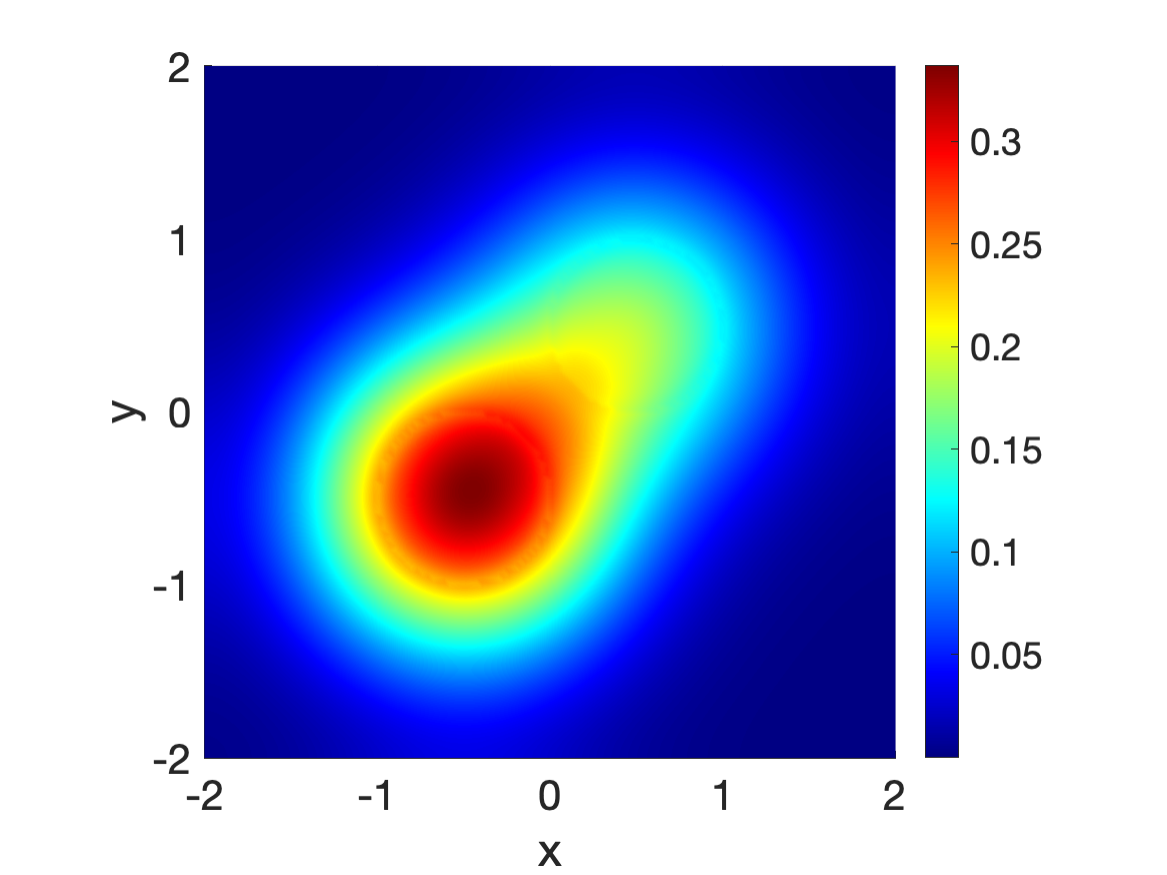}\hspace{-8.5pt}
 \includegraphics[width=0.25\textwidth]{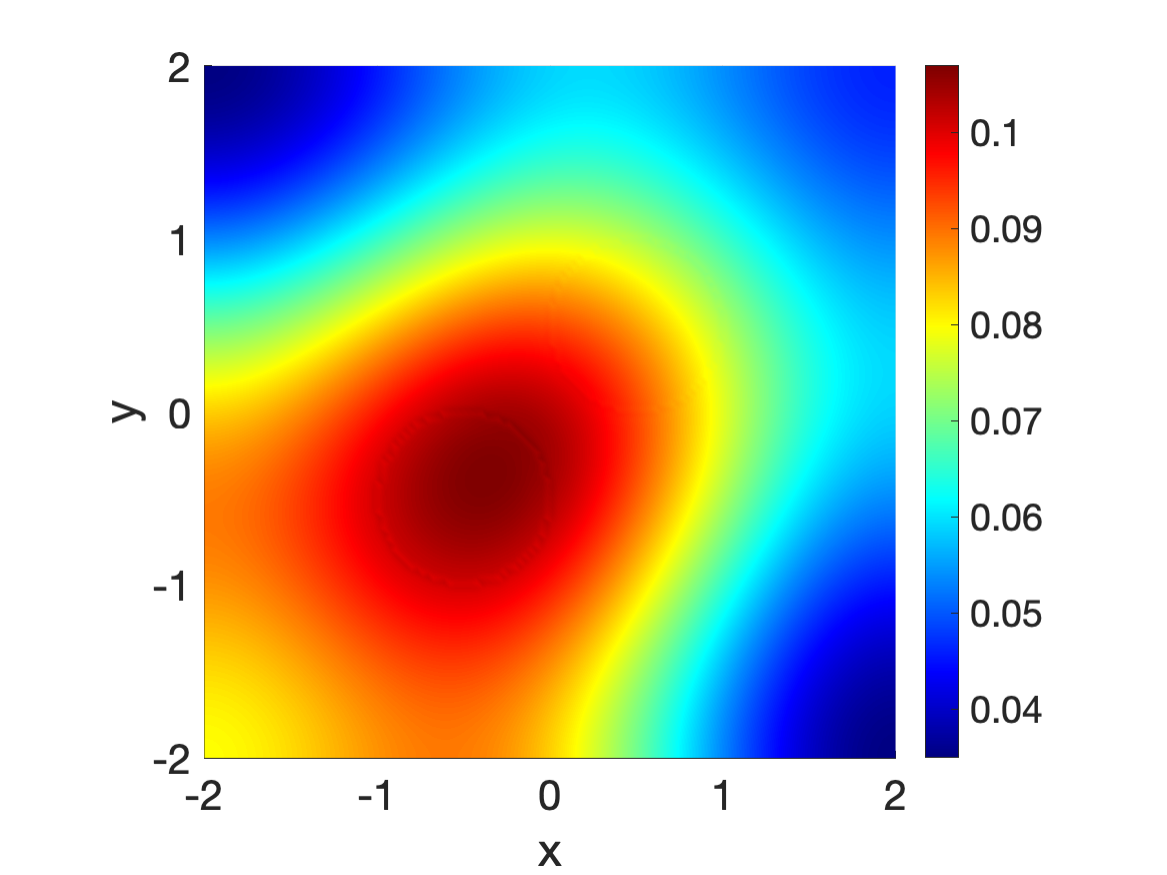}\hspace{-7.5pt}
 \includegraphics[width=0.25\textwidth]{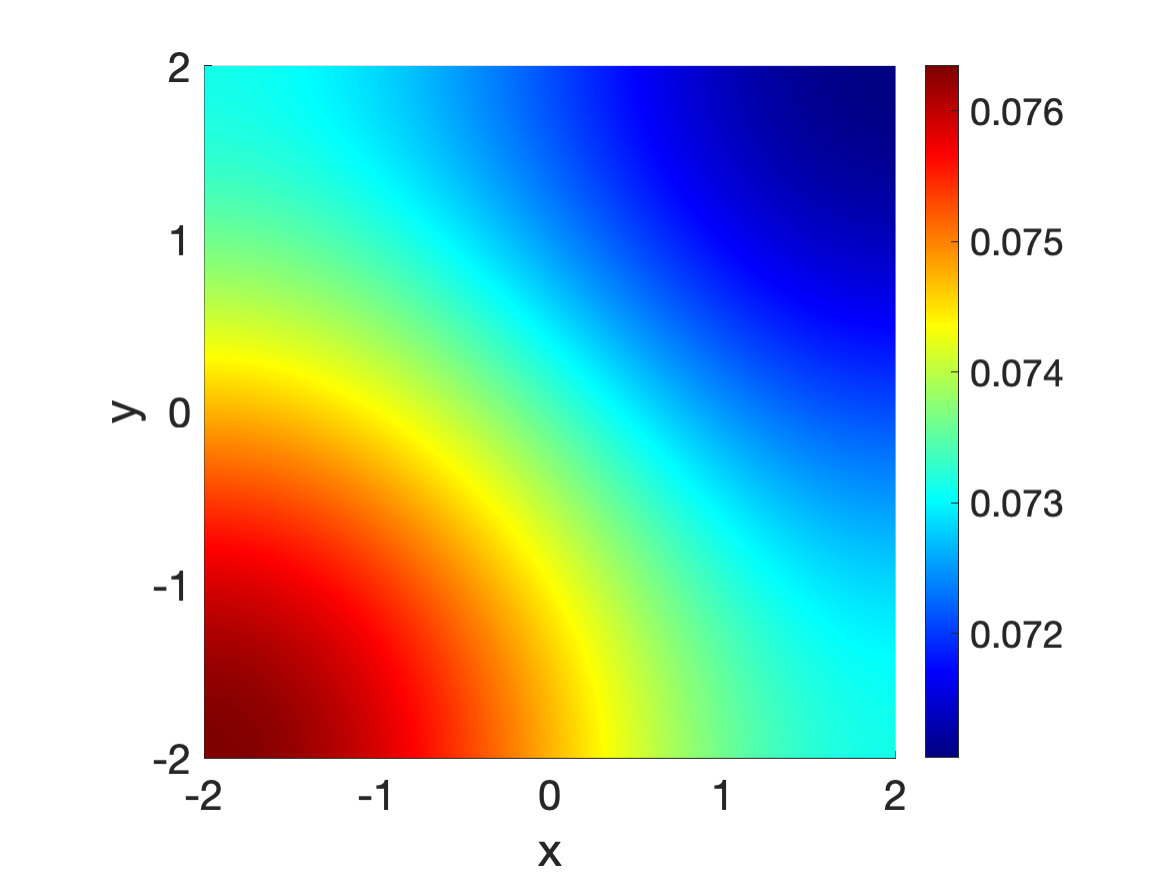}}
\centering{
  \includegraphics[width=0.25\textwidth]{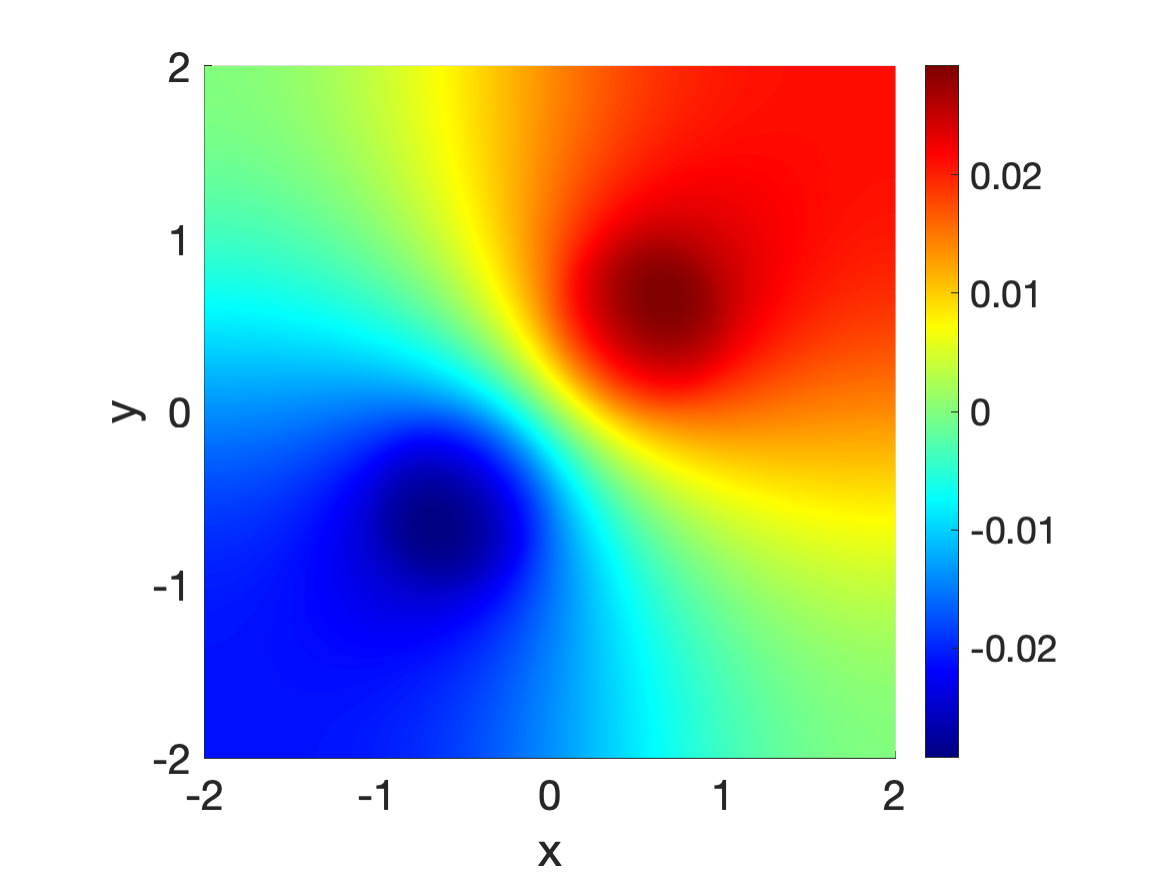}\hspace{-11pt}
 \includegraphics[width=0.25\textwidth]{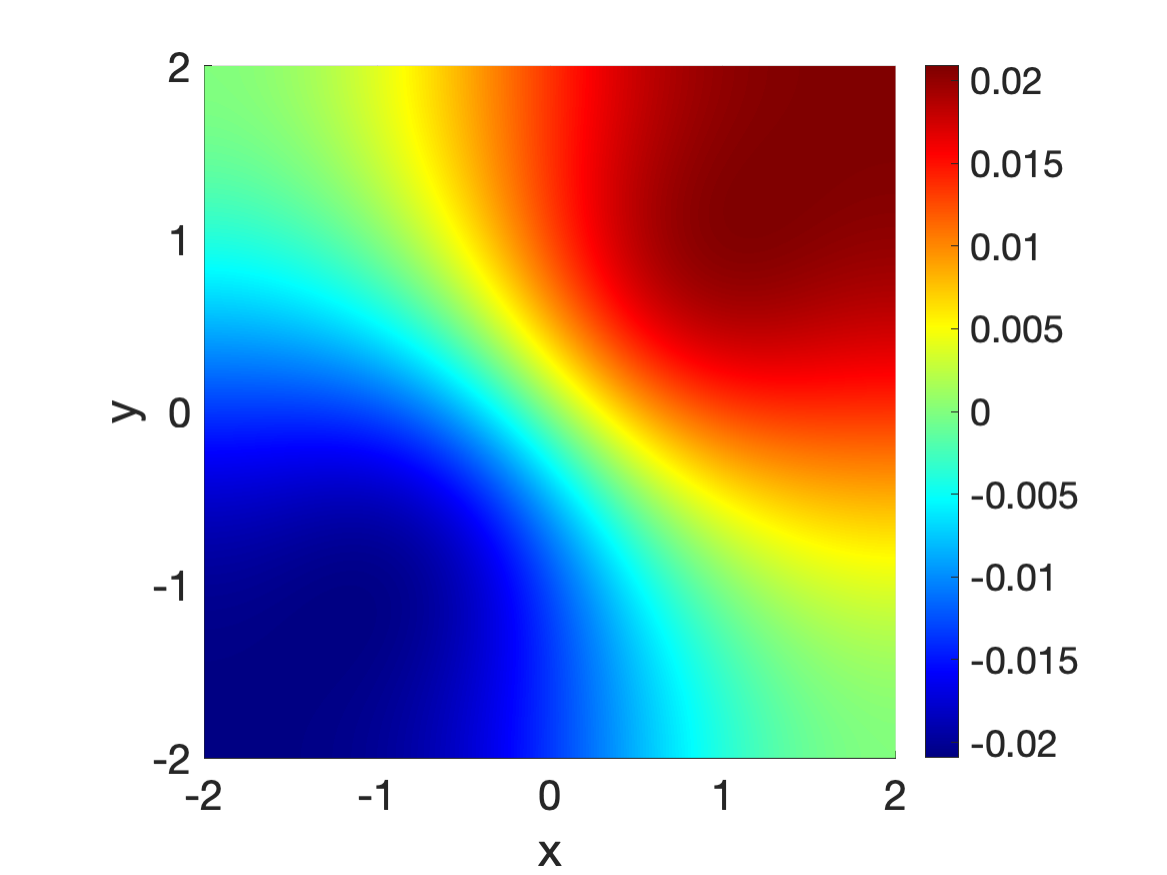}\hspace{-9.5pt}
 \includegraphics[width=0.25\textwidth]{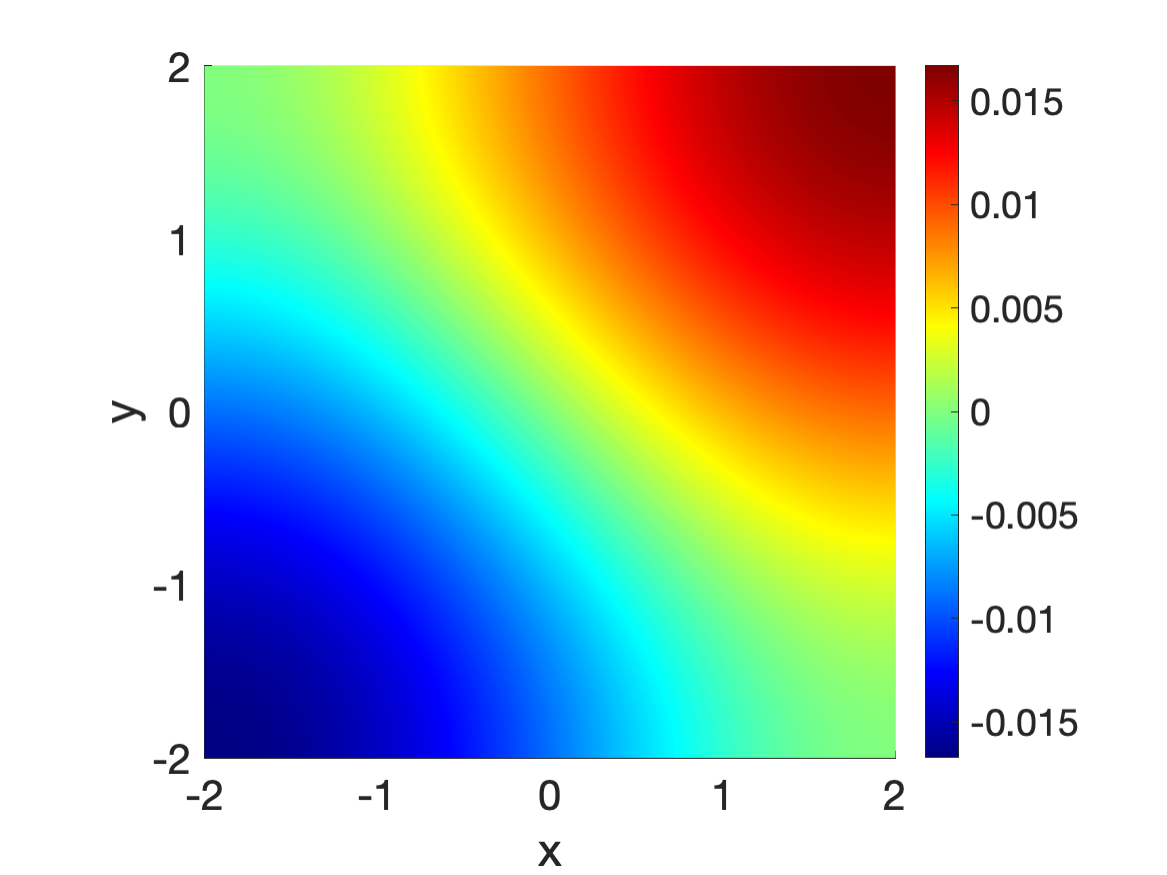}\hspace{-9.0pt}
 \includegraphics[width=0.25\textwidth]{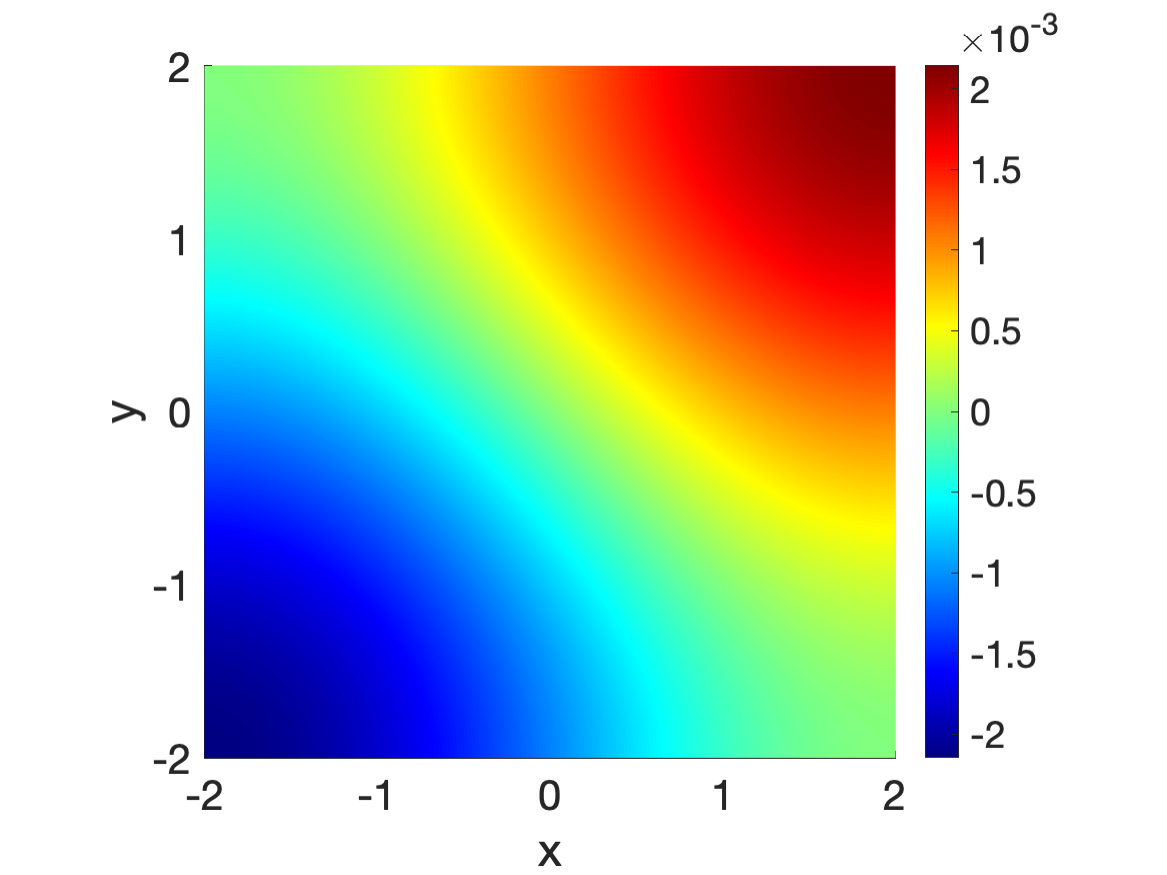}}
\caption{Example \ref{example2} with CNFDP.  Snapshots of $p$, $n$ and $\phi$ (from top to bottom) at $t = 0, 0.04, 0.2, 1$ (from left to right), respectively.}  
 \label{fig:Ex2p}
\end{figure}

\begin{figure}[!htp]
 \begin{minipage}{0.3\linewidth}
  \centering
  \includegraphics[width=0.9\linewidth]{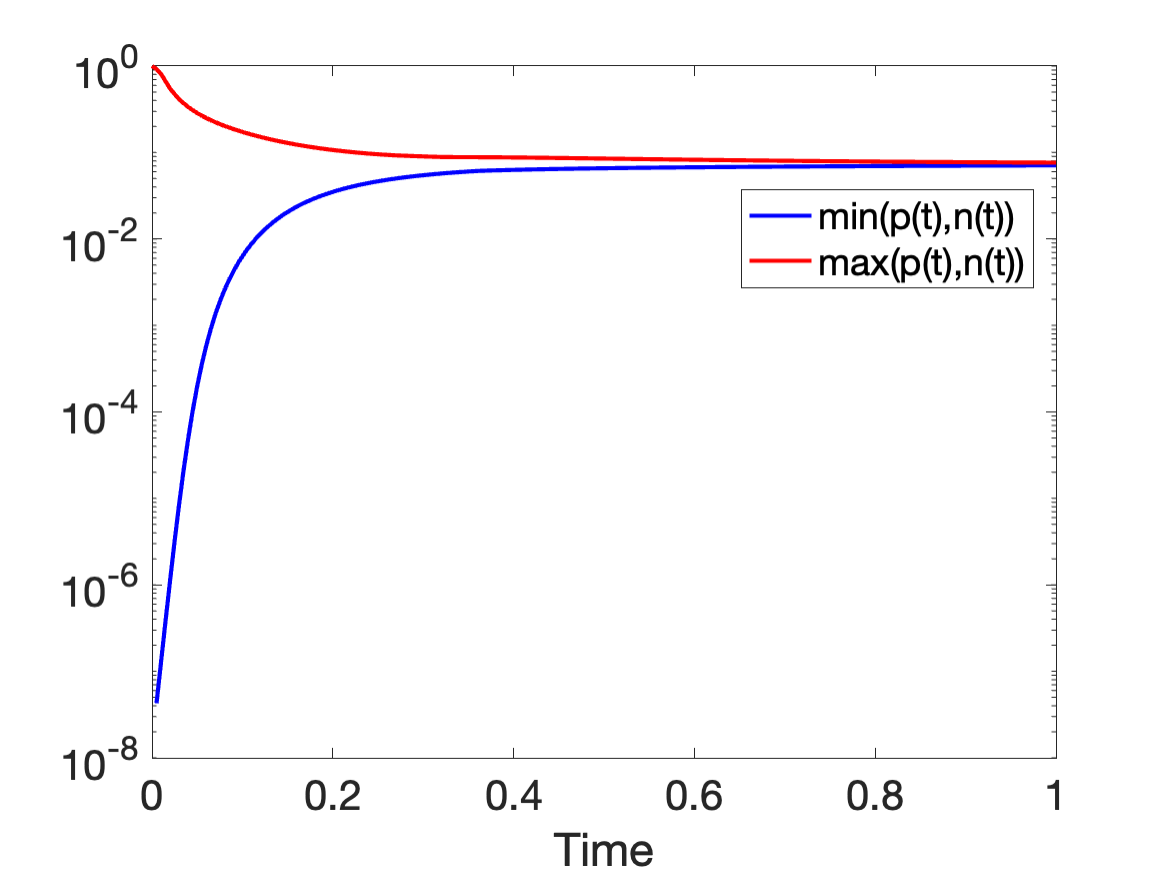}
 \end{minipage}
 \begin{minipage}{0.3\linewidth}
  \centering
  \includegraphics[width=0.9\linewidth]{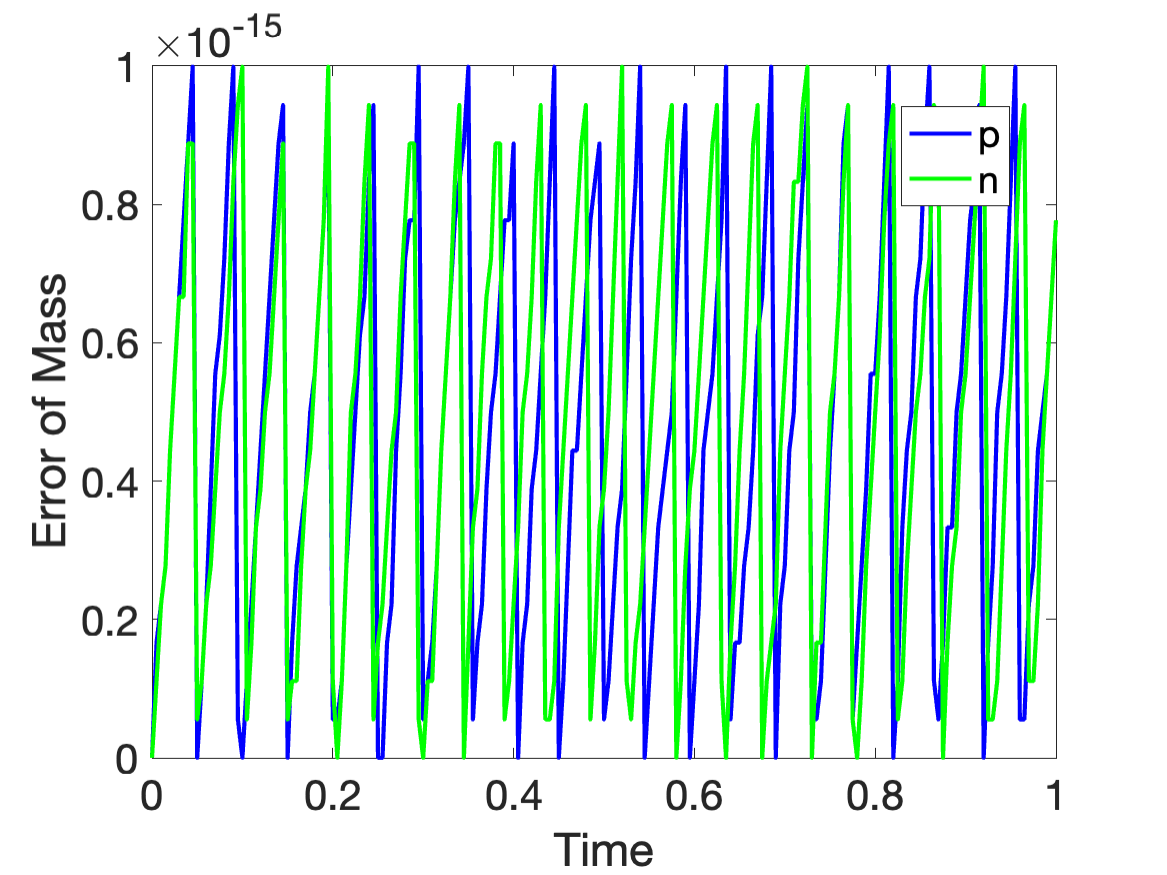}
 \end{minipage}
 \begin{minipage}{0.3\linewidth}
  \centering
  \includegraphics[width=0.9\linewidth]{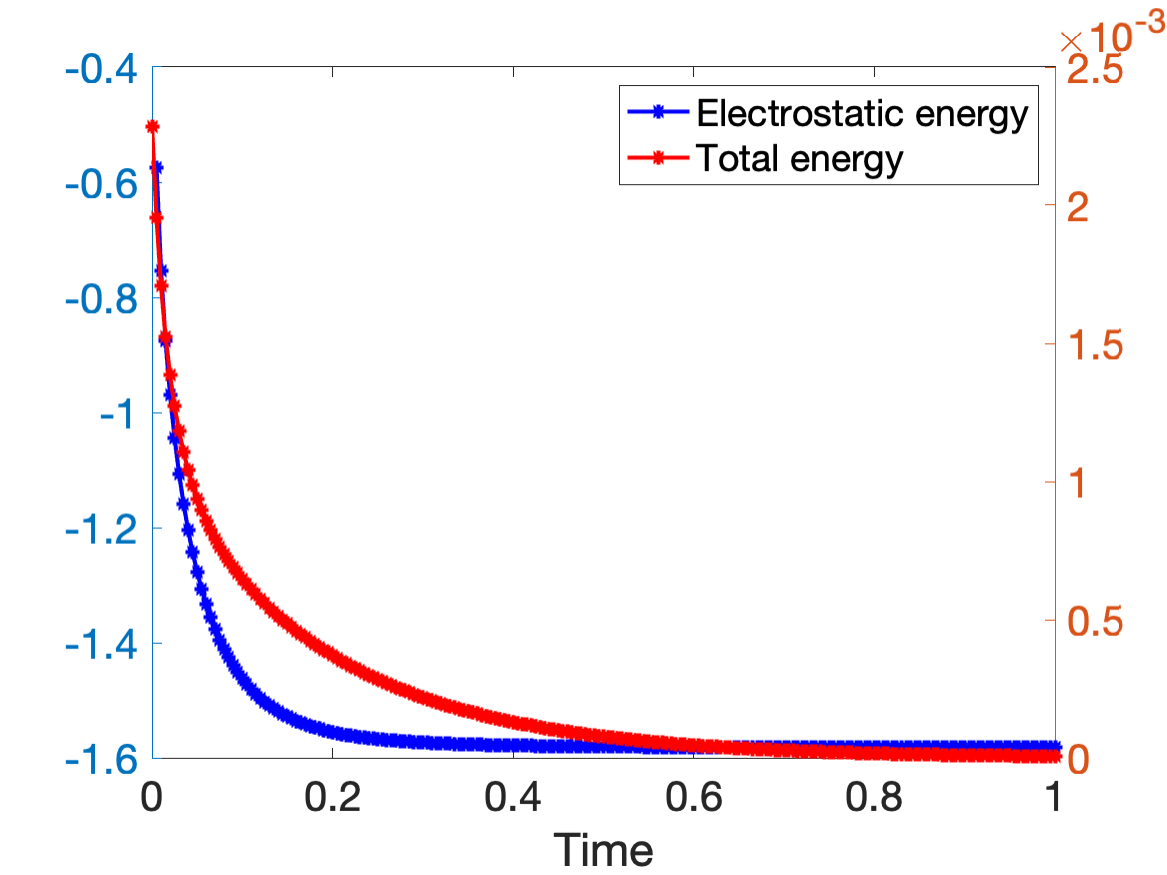}
 \end{minipage}
    \caption{ (Example \ref{example2} with CNFDP) Left: lower and upper bounds of $(n,p)$.  Middle:  error of discrete mass for    $p$ and $n$.   Right: evolution of total energy \eqref{energytotal} and the electric potential energy \eqref{energyel}. }
   \label{fig:EX2mass}
\end{figure}

 Figs. \ref{fig:Ex2p} and \ref{fig:EX2mass} depicts the evolution  of   $p$,  $n$ and $\phi$ along with time, where   positivity preserving and mass conservation for $p$ and $n$ are clearly observed. In addition, the total energy \eqref{energytotal} and the electric energy \eqref{energyel} are decreasing in time.

\subsection{3D case}
\begin{figure}[!htp]
\centering{
  \includegraphics[width=0.25\textwidth]{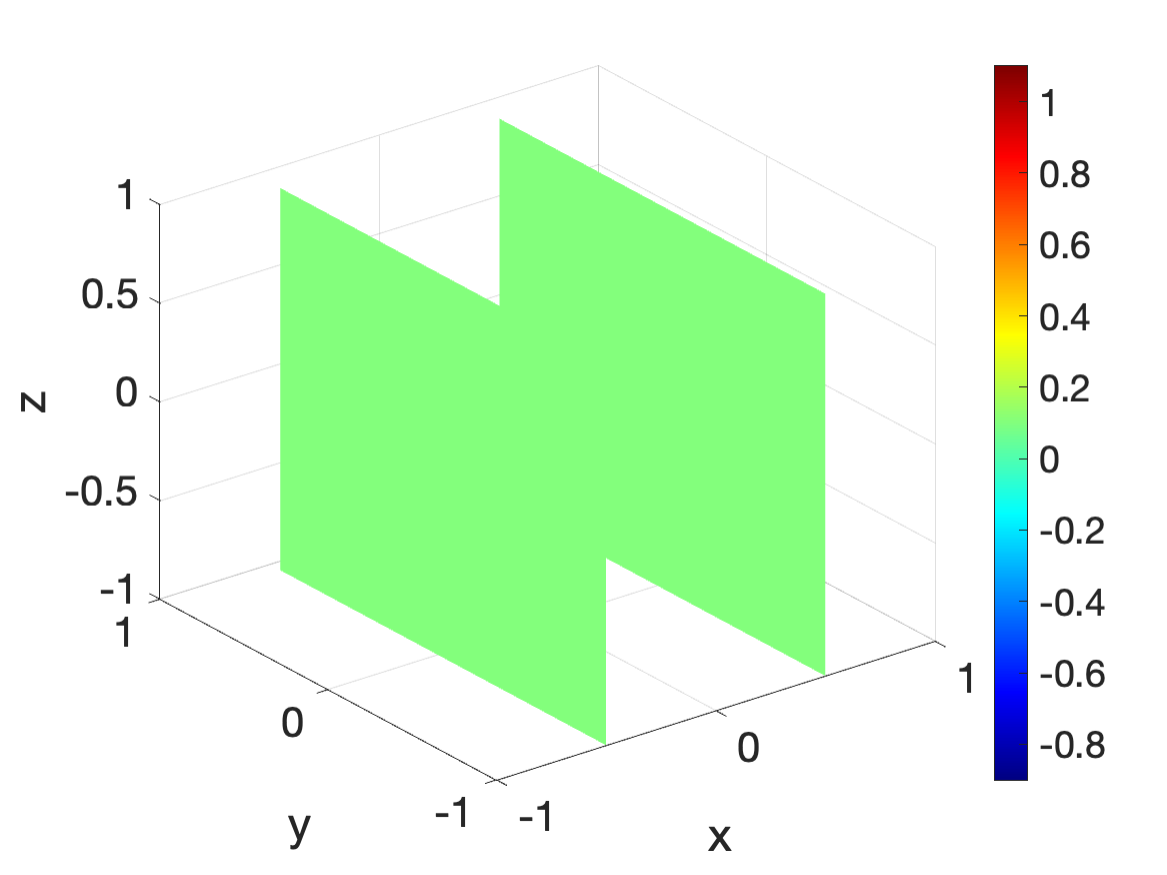}\hspace{-7pt}
 \includegraphics[width=0.25\textwidth]{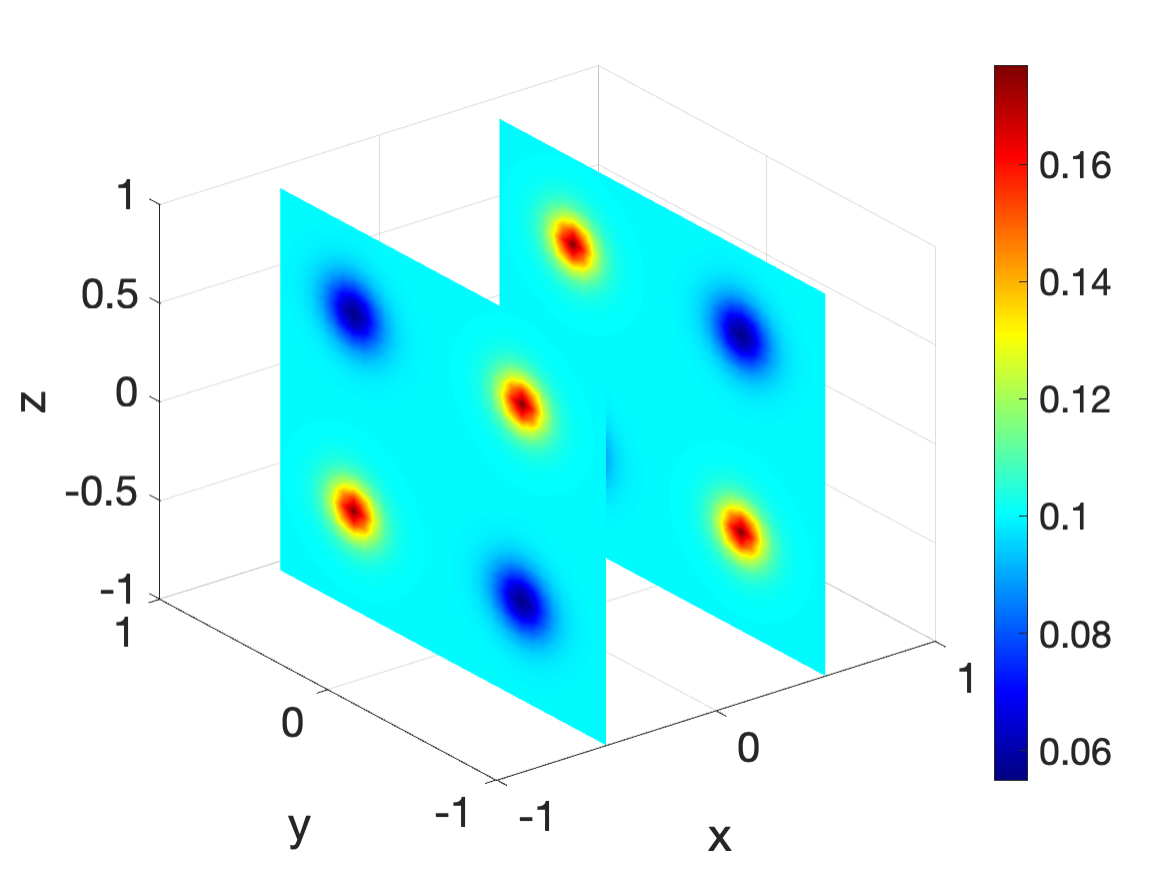}\hspace{-7pt}
 \includegraphics[width=0.25\textwidth]{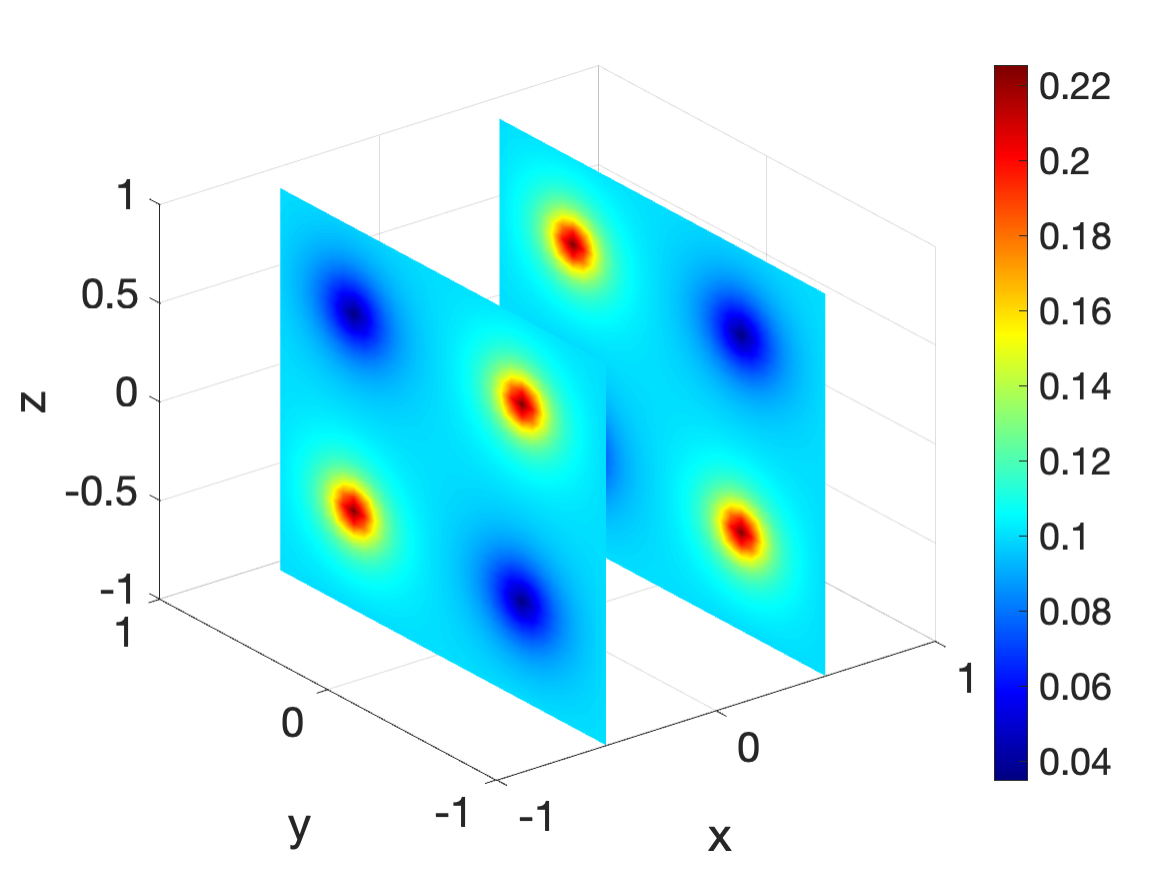}\hspace{-6pt}
 \includegraphics[width=0.25\textwidth]{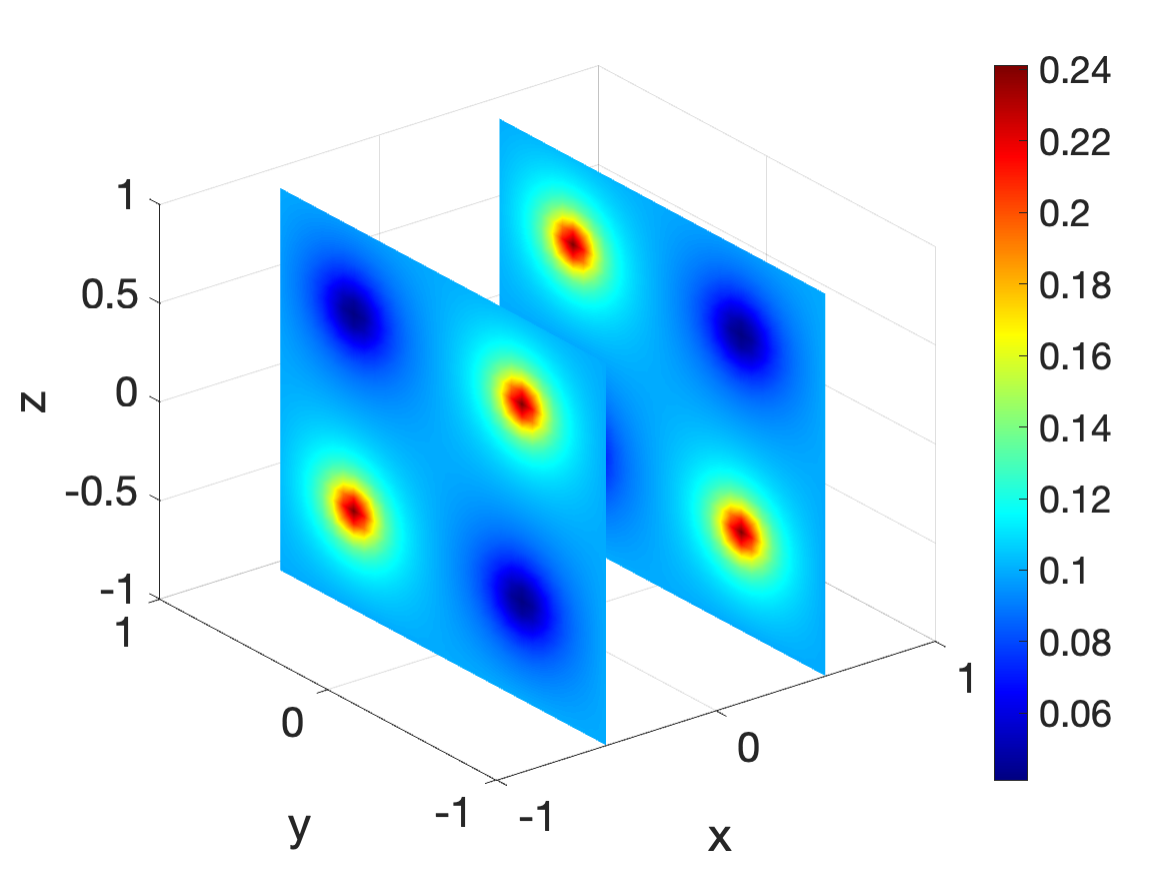}}
\centering{
  \includegraphics[width=0.25\textwidth]{Figures/e3p00.pdf}\hspace{-7pt}
 \includegraphics[width=0.25\textwidth]{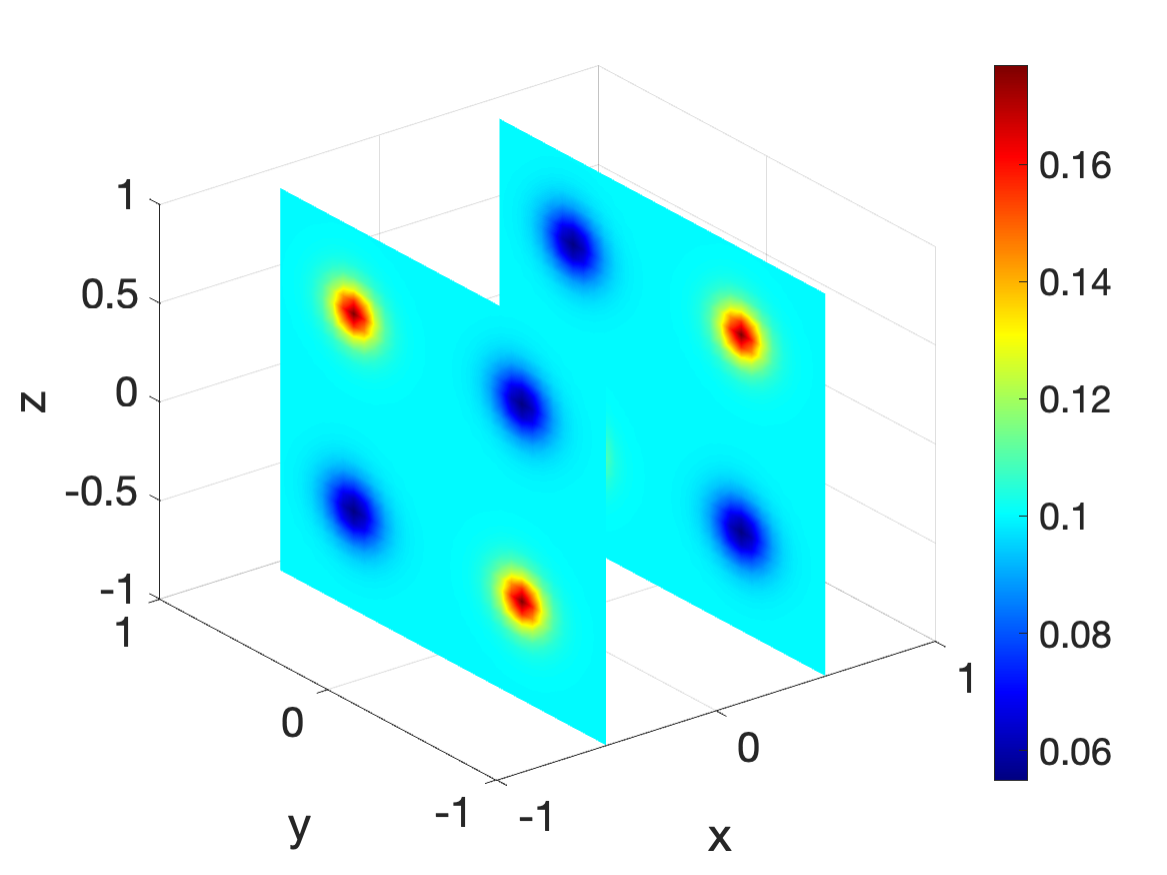}\hspace{-7pt}
 \includegraphics[width=0.25\textwidth]{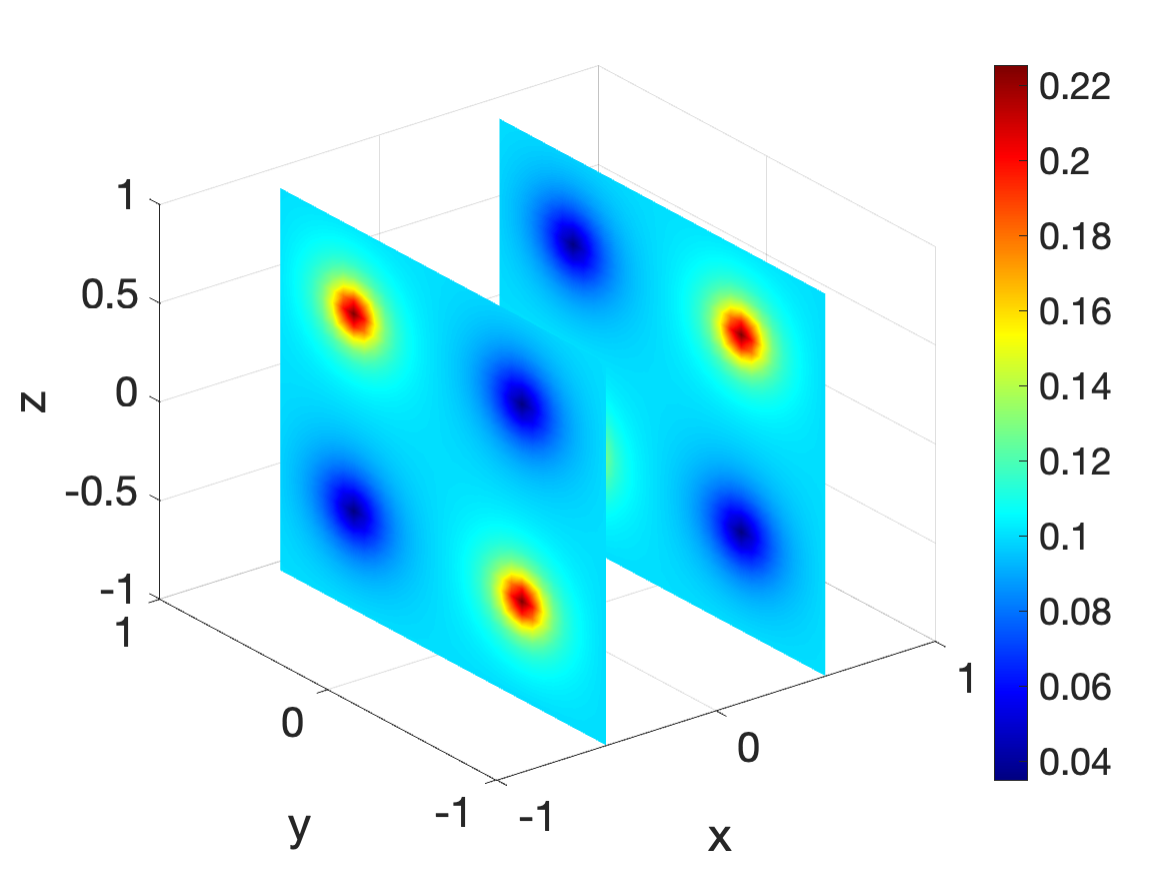}\hspace{-6pt}
 \includegraphics[width=0.25\textwidth]{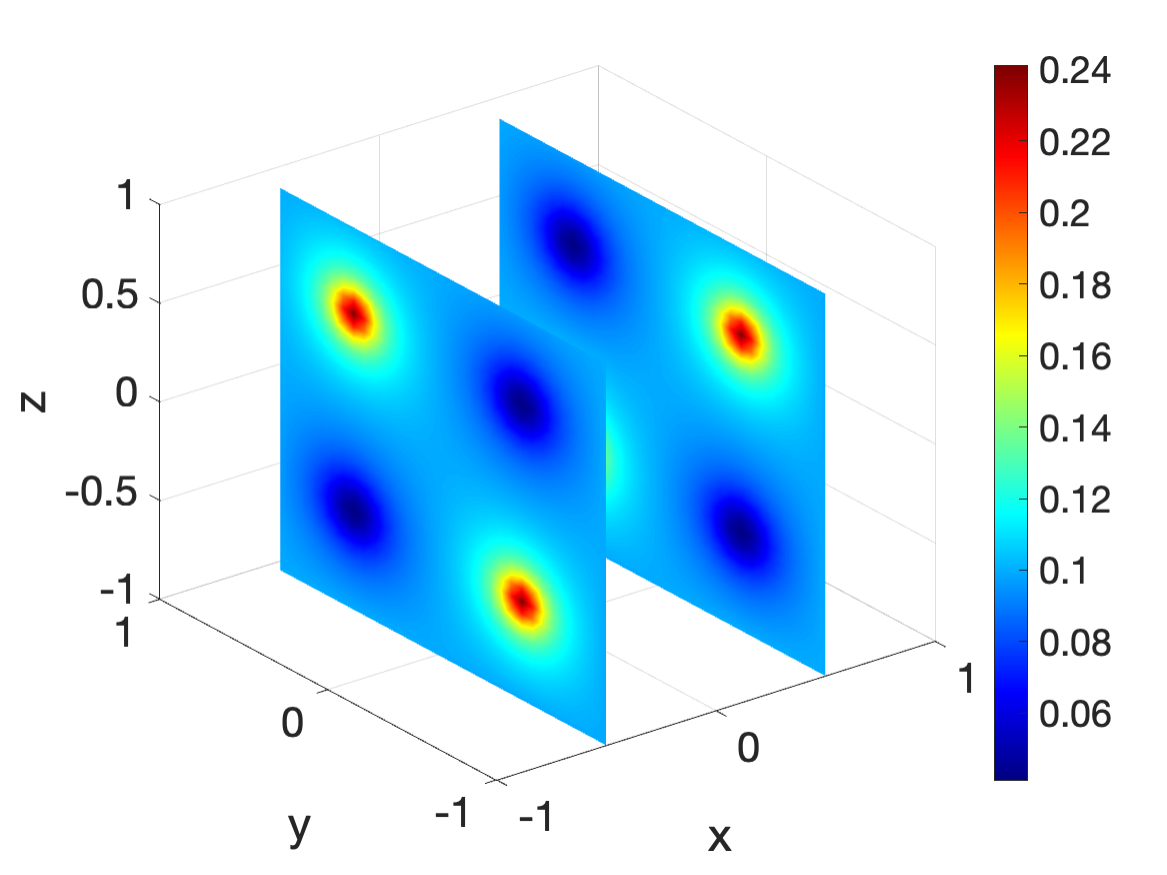}}
\centering{
  \includegraphics[width=0.25\textwidth]{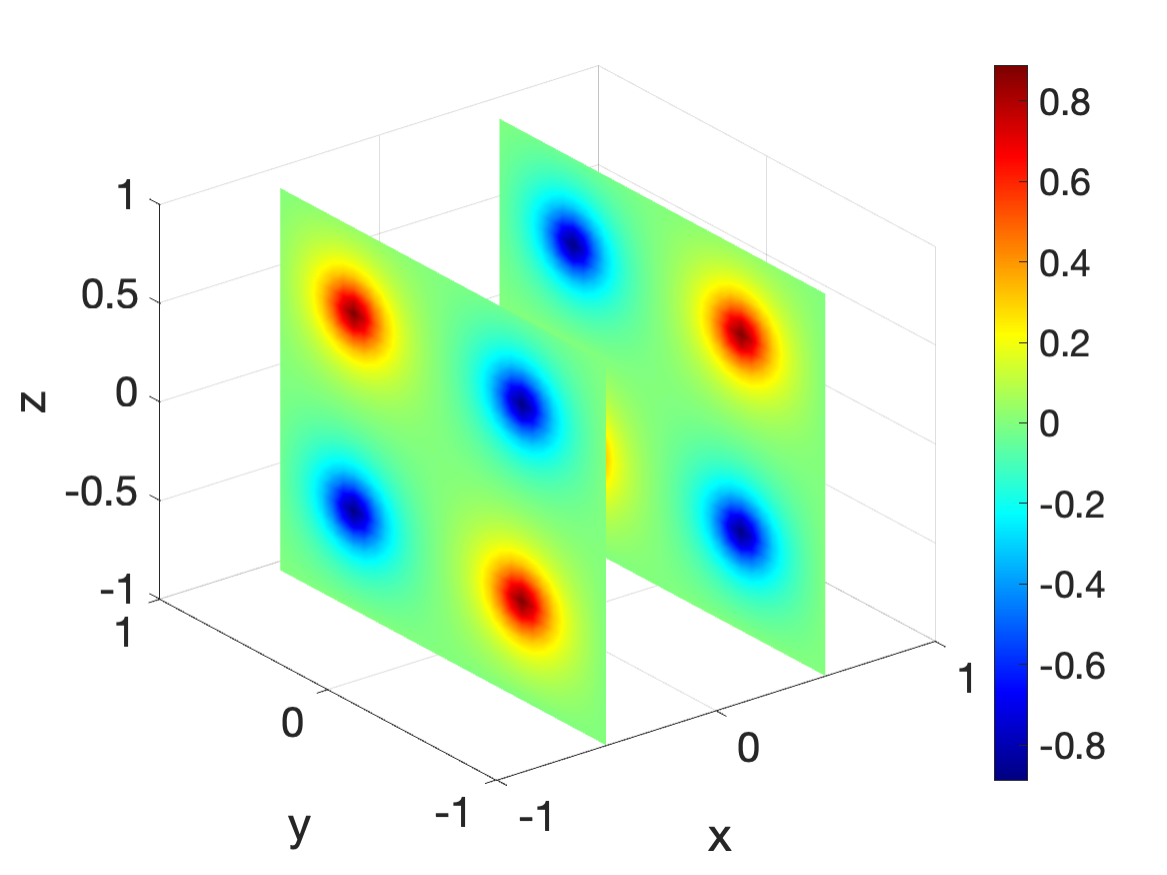}\hspace{-7pt}
 \includegraphics[width=0.25\textwidth]{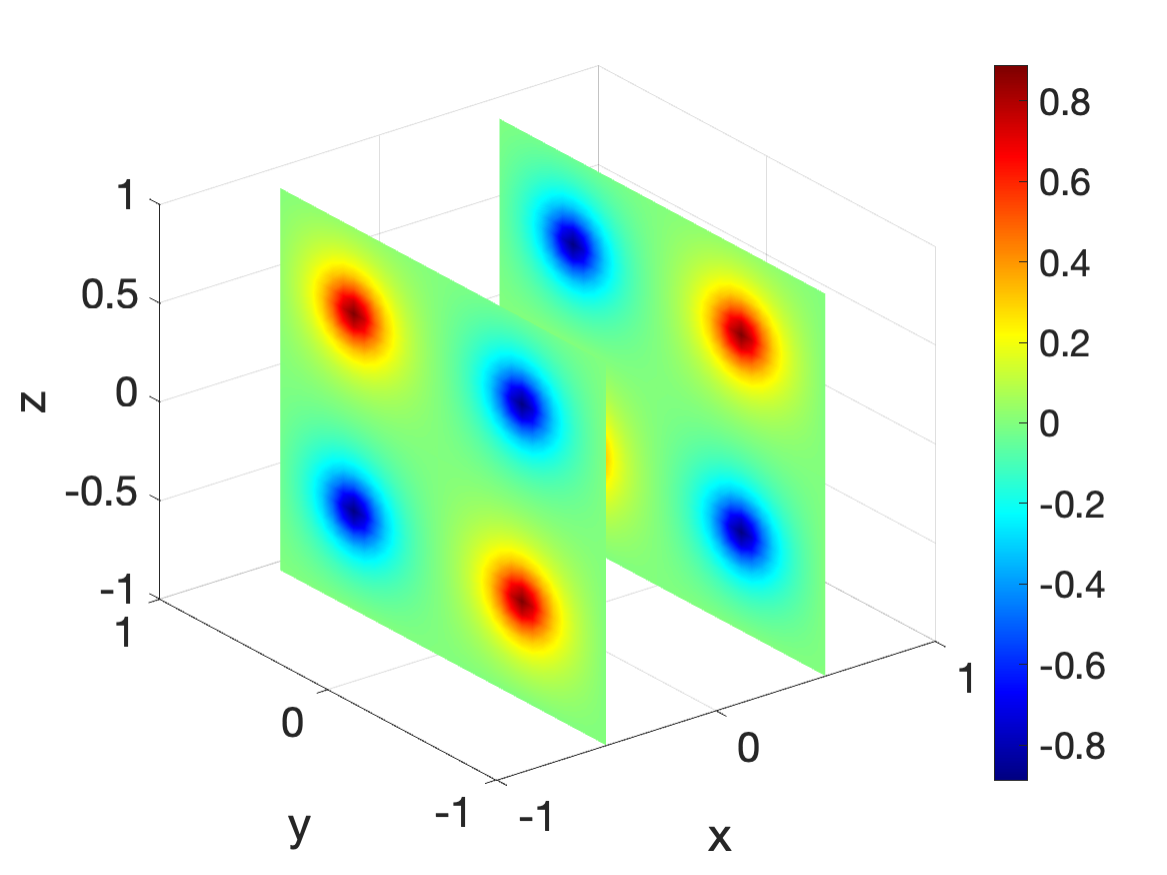}\hspace{-7pt}
 \includegraphics[width=0.25\textwidth]{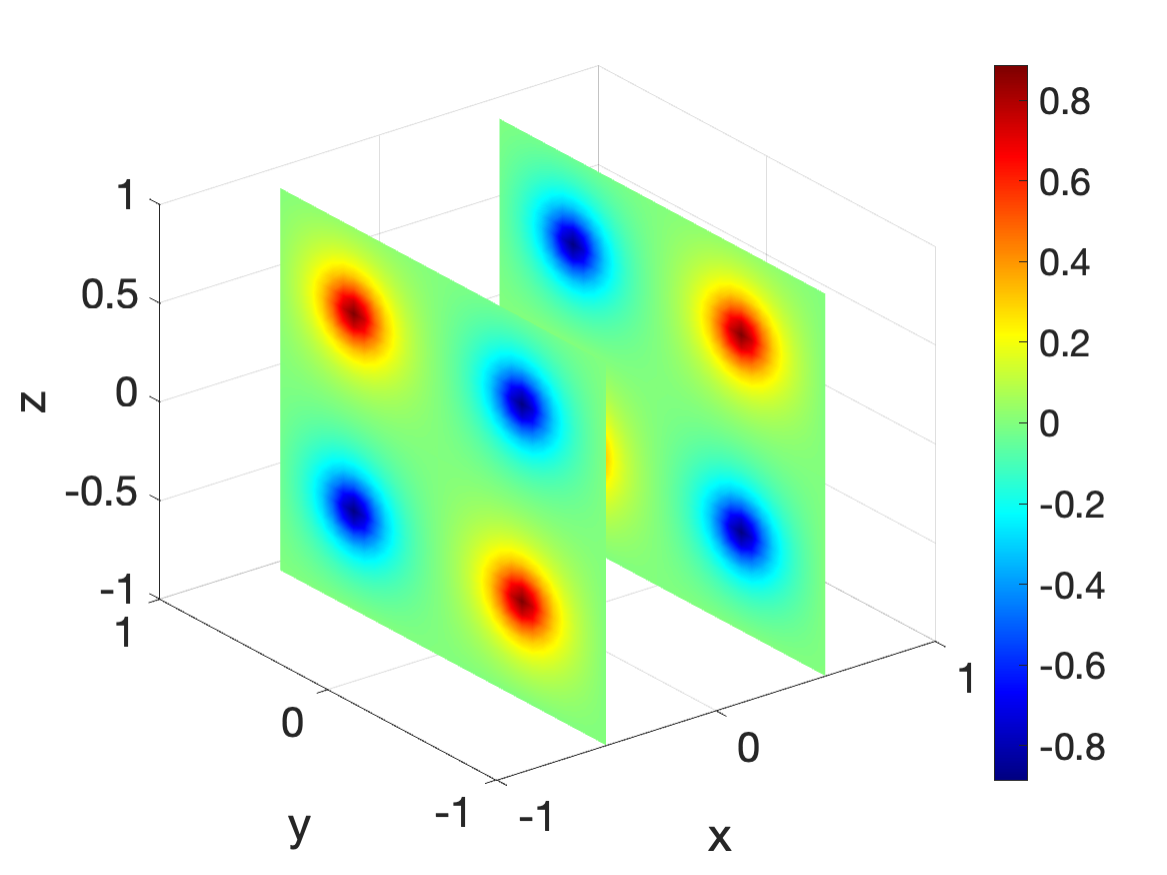}\hspace{-6pt}
 \includegraphics[width=0.25\textwidth]{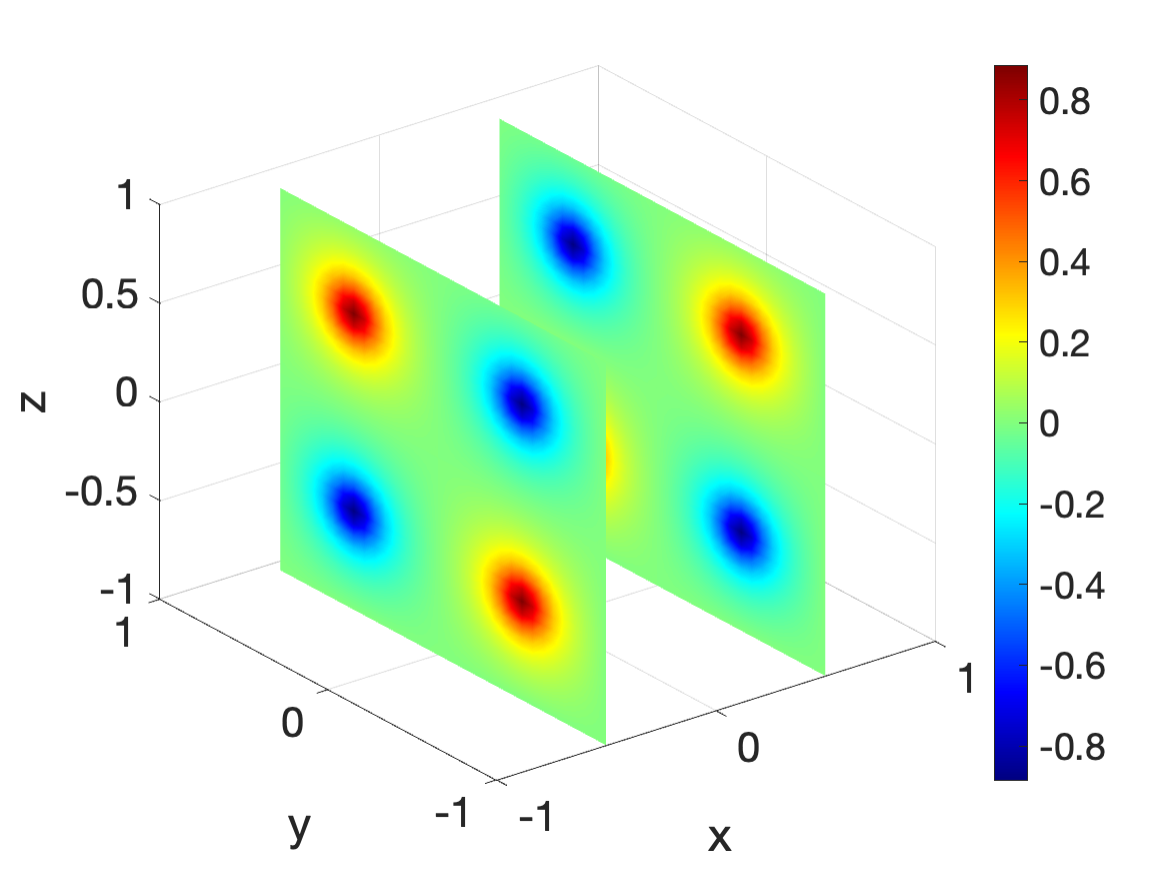}}
\caption{Example \ref{example4} with CNFDP.  Snapshots of $p$, $n$ and $\phi$  (from top to bottom) at $t = 0, 0.01, 0.05, 2.0$ (from left to right),  respectively. } 
 \label{fig:Ex3Dphi}
\end{figure}

\begin{figure}[!htp]
 \centering
 \begin{minipage}{0.3\linewidth}
  \centering
  \includegraphics[width=0.9\linewidth]{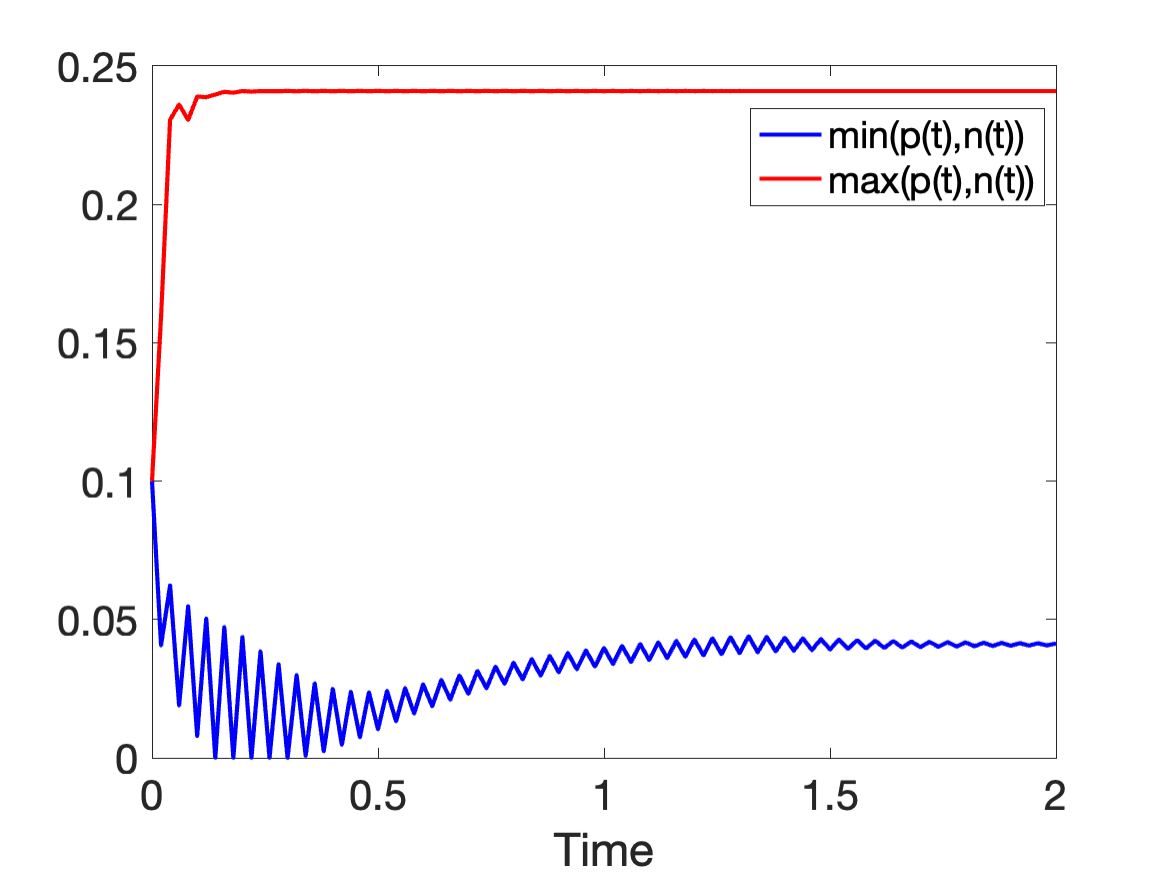}
 \end{minipage}
 \begin{minipage}{0.3\linewidth}
  \centering
  \includegraphics[width=0.9\linewidth]{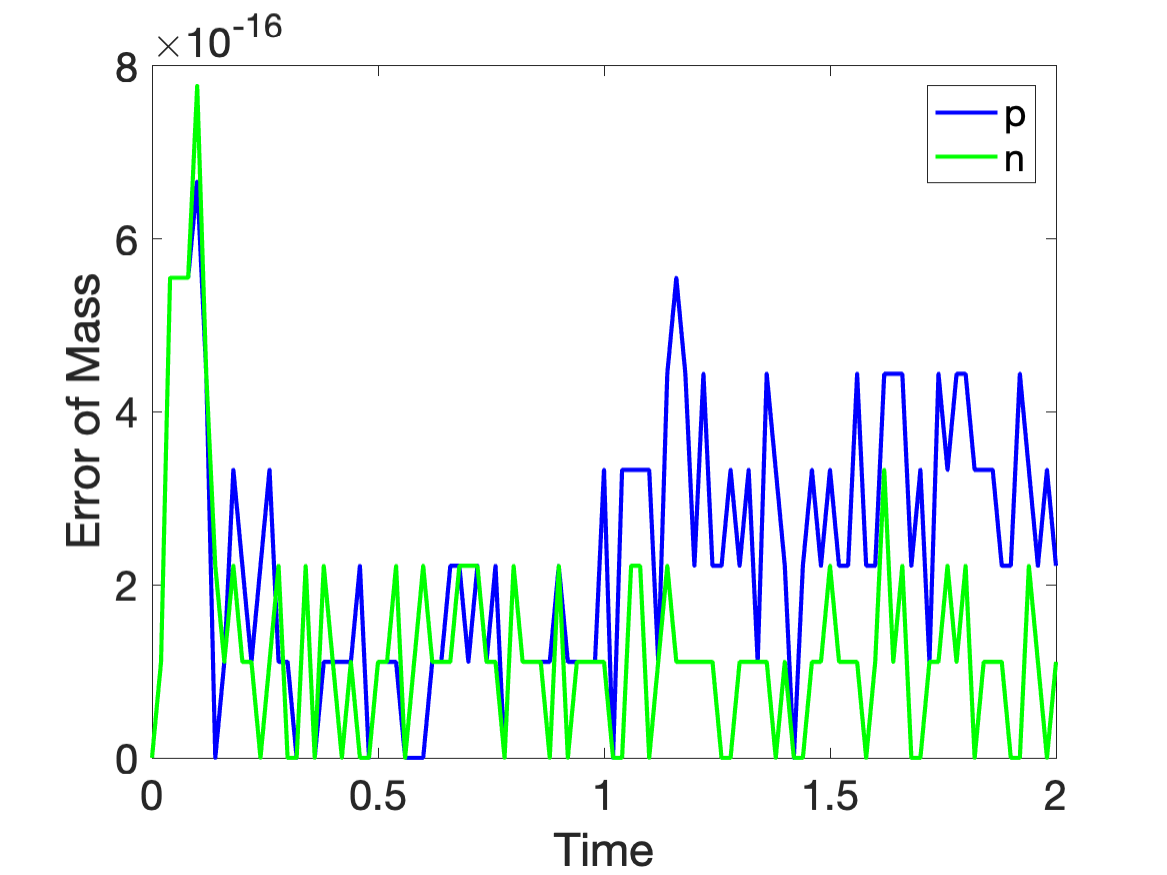}
 \end{minipage}
 \begin{minipage}{0.3\linewidth}
  \centering
  \includegraphics[width=0.9\linewidth]{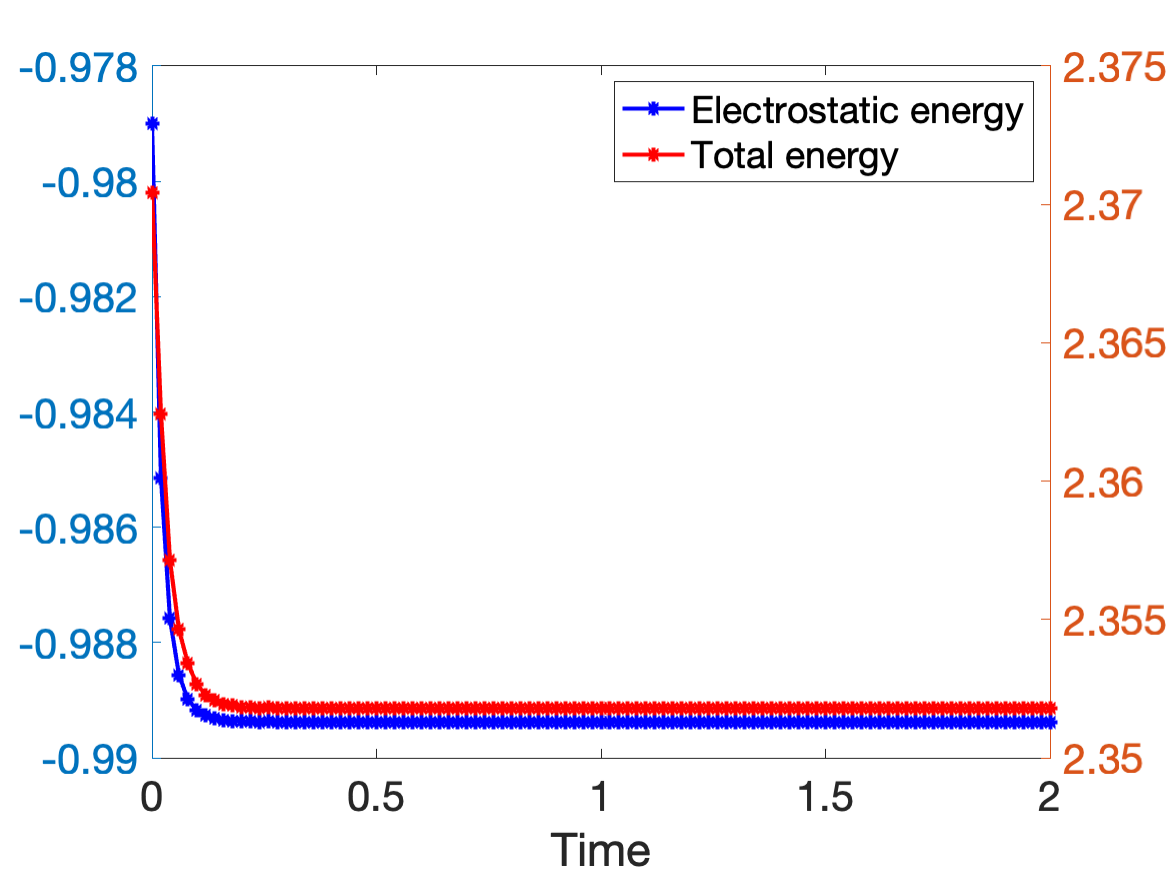}
 \end{minipage}
   \caption{Example \ref{example4}. Left:  lower and upper bounds of (n,p).  Middle: the error of discrete mass for $p$ and $n$.  Right: the total energy \eqref{energytotal} and  the electric potential energy  \eqref{energyel}.}
   \label{fig:EX3Dmass}
\end{figure}

\begin{example}\label{example4} 
The  PNP system \eqref{PDE} with  fixed charges on $\Omega=[-1,1]^3$ with periodic boundary conditions, i.e. the electric potential  equation is  replaced by
\begin{equation}
-\Delta \phi  = p - n + \rho,
\end{equation}
where the external fixed charge distribution $\rho$ is given by
\begin{equation*}
\begin{aligned}
 \rho(x,y,z)  &= 200 \sum\limits_{\varepsilon_{x},\varepsilon_{y},\varepsilon_{z}=\pm1}\varepsilon_{x}\varepsilon_{y}\varepsilon_{z}e^{-100[(x+\varepsilon_xx_i)^2+(y+\varepsilon_yy_i)^2+(z+\varepsilon_zz_i)^2]}.\end{aligned}
\end{equation*}
The  initial data is chosen as $p(x,y,z,0)  = 0.1$ and $n(x,y,z,0) =0.1$. 
\end{example}
We simulate Example \ref{example4} using CNFDP  with $\tau=2/100$ and $h=1/32$.
The initial distribution of electrostatic potential $\phi$ (or the distribution of the fixed charges'  potential) is obtained by solving the Poisson equation with initial concentrations.
Fig. \ref{fig:Ex3Dphi}   shows  the snapshots of the concentrations $p$ and $n$ and electric potential $\phi$ at different times.  As time evolves, the mobile ions are attracted by opposite fixed charges.  
Meanwhile,  the electrostatic potential of the fixed charges is gradually screened  by accumulated mobile ions of opposite signs.
From Fig. \ref{fig:EX3Dmass},  we can observe  the point-wise bounds of  $(n,p)$,  the discrete mass, the total energy \eqref{energytotal} and the electric potential energy \eqref{energyel}, which confirms that CNFDP is positivity preserving and mass conserving. In addition, the energies  \eqref{energytotal} and \eqref{energyel}  decrease monotonically in time.  

\section{Conclusion}
\label{sec:conclusion}
In this work, we developed efficient numerical schemes CNFDP and CNFDP2 to simulate the  PNP equation \eqref{PDE},  with positivity preserving and mass conservation at the discrete level. The key idea is to  solve a minimization problem to project (correct) the intermediate solution to satisfy the desired physically constraints.  By projection idea,  either in $L^2$ norm (CNFDP) or $H^1$  norm (CNFDP2),   the positivity and mass conservation of the numerical solutions can be guaranteed.  For CNFDP, error estimates were shown to be at $O(\tau^2+h^2)$ for the ion densities in $L^2$ norm and for the electric potential in $H^1$ norm. The proposed method with projection technique is a linearized
scheme, which can be implemented very efficiently, with  negligible computational cost in solving a nonlinear algebraic system. The method can be directly generalized to the case of multi-ions PNP system without any difficulty. Extension of the techniques and ideas to other more complex problems or higher order methods will be subject to  future investigation.

\section*{Appendix} 
\label{sec:appendix}
In this Appendix, we briefly describe the process about how to determine the  Lagrange multipliers $\lambda^{k+1}_h(\mathbf{x}),\eta^{k+1}_h(\mathbf{x})\in X$ and $\xi^{k+1},\gamma^{k+1}\in\mathbb{R}$. 

{\bf $L^2$ projection \eqref{projnonzero}}. From the  complementary relaxation conditions \eqref{projnonzero},   we  can obtain numerical solutions $p_{h}^{k+1}$  and $n_{h}^{k+1}$  expressed as
\begin{equation}\label{psolu}
(p_{h}^{k+1},\lambda_h^{k+1})=\left\{\begin{array}{ll}
\left[\tilde{p}_{h}^{k+1}-\xi^{k+1},0\right], & \text { if } \tilde{p}_{h}^{k+1}-\xi^{k+1}>0,\\
\left[0,-(\tilde{p}_{h}^{k+1}-\xi^{k+1})\right], & \text { otherwise},\end{array}\right.
\end{equation}
and
\begin{equation}\label{nsolu}
(n_{h}^{k+1},\eta_h^{k+1})=\left\{\begin{array}{ll}
\left[\tilde{n}_{h}^{k+1}-\gamma^{k+1},0\right], & \text { if } \tilde{n}_{h}^{k+1}-\gamma^{k+1}>0, \\
\left[0,-(\tilde{n}_{h}^{k+1}-\gamma^{k+1})\right], & \text { otherwise,}\end{array}\right.
\end{equation}
where  the Lagrange multipliers $\xi^{k+1},\gamma^{k+1}\in\mathbb{R}$ are determined by the mass constraints, e.g. for $\xi^{k+1}$,
\begin{equation*}
h^2 \sum_{0<\tilde{p}_{h}^{k+1} (\mathbf{x})- \xi^{k+1}} \left(\tilde{p}_{h}^{k+1}(\mathbf{x})- \xi^{k+1}\right) = \langle{p}^{0}_h,1\rangle.
\end{equation*}
Hence $\xi^{k+1}\in\mathbb{R}$ is the solution to the nonlinear algebraic equation
\begin{equation}\label{ksii}
F(\xi^{k+1})=h^2 \sum_{i,j=1}^N   \left( \tilde{p}_{h}^{k+1}(\mathbf{x}_{i,j})- \xi^{k+1}\right)^+ - \langle{p}^{0}_h,1\rangle=0,
\end{equation}
where $f^+=\max\{f,0\}$ for $f\in\mathbb{R}$.
Semi-smooth Newton methods \cite{Kojima1986,2003Finite} can be employed to solve \eqref{ksii}, i.e. for some initial guess $\xi_0\in\mathbb{R}$, find the root of $F(\xi)=0$ by updating
\begin{equation}\label{semismoothNewton}
\xi_{s+1}=\xi_{s}-V^{-1}(\xi_{s})F\left(\xi_{s}\right),\quad s=0,1,\cdots,
\end{equation}
where $V(\xi_s)$ is  a generalized derivate in semi-smooth Newtown methods as
\begin{equation*}
V(\xi_{s})=-h^2 \sum_{i,j=1}^N \text{sgn}( (\tilde{p}_{h}^{k+1}(\mathbf{x}_{i,j})- \xi_s)^+),
\end{equation*}
and $\text{sgn}(\cdot)$ is the sign function with $\text{sgn}(s)=1$ ($s>0$), $\text{sgn}(0)=0$ and $\text{sgn}(s)=-1$ ($s<0$).
Noticing $\xi^{k+1}$ is supposed to be of small magnitude   \cite{ChengQing1} and  $\xi^{k+1} \geq 0$ (Lemma \ref{lem:lampoistive}),  we can choose $\xi_{0}=0$  to start the semi-smooth Newton iterations. 
Once $\xi^{k+1}$ is known, we can update $\left(p_{h}^{k+1}, \lambda_{h}^{k+1}\right)$ according to \eqref{psolu}.
$\gamma^{k+1}$ and $n_h^{k+1}$ can be computed in the same way.  In all our numerical experiments,  \eqref{semismoothNewton} converges in only one iteration so that the cost of solving \eqref{ksii} is negligible.  The secant method can be also used to solve \eqref{ksii} \cite{ChengQing1}.

{\bf $H^1$ projection \eqref{hproj}-\eqref{hpproj}}. We present a semi-smooth Newton method to find the  Lagrange multipliers $\lambda^{k+1}_h(\mathbf{x}),\eta^{k+1}_h(\mathbf{x})\in X$ and $\xi^{k+1},\gamma^{k+1}\in\mathbb{R}$ under $H^1$ projection.
The KKT condition \eqref{hproj}-\eqref{hpproj} for $\lambda^{k+1}$ and $\xi^{k+1}$ can be transformed to an equivalent  system as follows
\begin{equation}\label{F1F2}
\begin{aligned}
  &F_1(U,\xi;\mathbf{x})= -\Delta_h {U}^{+} (\mathbf{x}) + U(\mathbf{x})+\xi^{k+1}- (I-\Delta_h)\tilde{p}^{k+1}_h(\mathbf{x})=0,\quad \mathbf{x}\in\Omega_h,\\
   & F_2(U,\xi)=  h^2 \sum_{i,j=1}^N {U}^{+}(\mathbf{x}_{i,j})-M_0=0,
\end{aligned}
\end{equation}
where  $ U(\mathbf{x})=U^+(\mathbf{x}) - U^-(\mathbf{x})$ , $U^+ = \max(U^{k+1}(\mathbf{x}),0)$,  $ U^- = \max(-U^{k+1}(\mathbf{x}),0)$,
and $U^+(\mathbf{x})={p}^{k+1}_h(\mathbf{x})$, $U^-(\mathbf{x})=\lambda_{h}^{k+1}(\mathbf{x})$. $\eta^{k+1}$ and $\gamma^{k+1}$ can be solved in the same way.

The semi-smooth Newton method can be applied to solve \eqref{F1F2}. We start with the generalized Jacobian to be used in the semi-smooth Newton iterations. In a simplified form, at $(U(\mathbf{x}),\xi)\in X\times\mathbb{R}$, the generalized Jacobian $J$ (acting on a vector $(V(\mathbf{x}),\zeta)\in X\times\mathbb{R}$) can be written as
\begin{equation}\label{F1F2D}
J\begin{pmatrix}
V(\mathbf{x})\\
\zeta
\end{pmatrix}=\begin{pmatrix}
   I - \Delta_h \left(\text{sgn}(U^+(\mathbf{x})) V(\mathbf{x})\right)+\zeta \\
   h^2 \sum  (\text{sgn} (U^+(\mathbf{x})) V(\mathbf{x})
\end{pmatrix}.
\end{equation}
By a careful computation based on \eqref{F1F2}-\eqref{F1F2D},  a semi-smooth Newton method for solving \eqref{F1F2} with initial $\xi_0$ and $U_0(\mathbf{x})$ can be designed as
follows:
\begin{align*}
\xi_{s+1} =& \xi_{s} + \frac{F_2(U_s,\xi_s)+h^2 \sum_{i,j=1}^N \text{sgn}(U_s^+(\mathbf{x}_{i,j}))V_1(\mathbf{x}_{i,j})}{h^2 \sum\limits_{i,j=1}^N (\text{sgn}(U^+_s(\mathbf{x}_{i,j}))V_2(\mathbf{x}_{i,j})},\\
U_{s+1} =& U_{s} +V_s-(\xi_{s+1} - \xi_{s})V_2,
\end{align*}
where $V_1(\mathbf{x})=( I - \Delta_h \text{sgn}(U_s^+) )^{-1}(-F_1(U_s,\xi_s))$ (by abuse of notation) and $V_2(\mathbf{x})=(   I - \Delta_h \text{sgn}(U_s^+) )^{-1}\mathbf{1}($ ($\mathbf{1}$ is the constant $1$ function on $\Omega_h$).  For solving $V_1$ and $V_2$, we can use  $(  I - \Delta_h )^{-1}$ as a pre-conditioner. 
When Newton iteration converges for some $s\ge0$,  we have $p^{k+1}(\mathbf{x})= U^+_{s+1}(\mathbf{x})$.  To start the semi-smooth Newton iteration, we may set $\xi_{0}=0$ and $U_0(\mathbf{x})=\tilde{p}^{k+1}_h(\mathbf{x})$.


\bibliographystyle{siamplain}
\bibliography{PNP_mp_arxiv}
\end{document}


\maketitle

\section{A detailed example}

Here we include some equations and theorem-like environments to show
how these are labeled in a supplement and can be referenced from the
main text.
Consider the following equation:
\begin{equation}
  \label{eq:suppa}
  a^2 + b^2 = c^2.
\end{equation}
You can also reference equations such as \cref{eq:matrices,eq:bb} 
from the main article in this supplement.

\lipsum[100-101]

\begin{theorem}
An example theorem.
\end{theorem}

\lipsum[102]
 
\begin{lemma}
An example lemma.
\end{lemma}

\lipsum[103-105]

Here is an example citation: \cite{KoMa14}.

\section[Proof of Thm]{Proof of \cref{thm:bigthm}}
\label{sec:proof}

\lipsum[106-112]

\section{Additional experimental results}
\Cref{tab:smfoo} shows additional
supporting evidence. 

\begin{table}[htbp]
\footnotesize
  \caption{Example table.}\label{tab:smfoo}
\begin{center}
  \begin{tabular}{|c|c|c|} \hline
   Species & \bf Mean & \bf Std.~Dev. \\ \hline
    1 & 3.4 & 1.2 \\
    2 & 5.4 & 0.6 \\ \hline
  \end{tabular}
\end{center}
\end{table}

\bibliographystyle{siamplain}
\bibliography{references}